\documentclass[12pt,a4paper]{article}
\usepackage{amsmath}
\usepackage{amssymb}\usepackage{bbm}
\usepackage{booktabs}
\usepackage{threeparttable}
\usepackage{multirow}
\usepackage{graphicx}
\usepackage{subfigure}
\usepackage[colorlinks,linkcolor=blue,citecolor=blue]{hyperref}
\usepackage{float}
\usepackage{indentfirst}
\usepackage[numbers,sort]{natbib}

\newcommand{\dx}{\Delta x}
\newcommand{\dy}{\Delta y}
\newcommand{\dt}{\Delta t}

\usepackage{array}
\newcolumntype{L}[1]{>{\raggedright\arraybackslash}p{#1}}

\usepackage{theorem}
\theorembodyfont{\normalfont}

\newtheorem{rem}{Remark}[section]

\newtheorem{example}{Example}[section]

\numberwithin{equation}{section}
\numberwithin{figure}{section}
\numberwithin{table}{section}

\newcommand{\nv}{{\nu}}
\newcommand{\betamaxplot}[1]{\includegraphics[width=7.5cm,angle=0]{#1}}
\usepackage{caption}
\begin{document}

\captionsetup[figure]{labelfont={bf},name={Fig.},labelsep=period}

\baselineskip=1.5pc
\begin{center}
	{\Large {\bf An efficient and robust fifth-order HWENO scheme with gradient reconstruction for compressible Navier--Stokes equations}}
\end{center}

\begin{center}
	Zijun Guo\footnote{School of Mathematical Sciences, Xiamen University, Xiamen, Fujian 361005, P.R. China. E-mail: guozj@stu.xmu.edu.cn.},
	Chuan Fan\footnote{Key Laboratory of Mathematical Modelling and High Performance Computing of Air Vehicles (NUAA), MIIT, Nanjing University of Aeronautics and Astronautics, Nanjing, Jiangsu 210016, China. E-mail: cfan\_math@nuaa.edu.cn.},
	Zhuang Zhao\footnote{Corresponding author. School of Mathematical Sciences and Fujian Provincial
Key Laboratory of Mathematical Modeling and High-Performance Scientific Computing, Xiamen University, Xiamen, Fujian 361005, P.R. China. E-mail: zzhao@xmu.edu.cn.}
\end{center}

\vspace{.05in}
\centerline
{\bf Abstract\ }
\bigskip
In this paper, we propose an efficient and robust fifth-order finite-volume Hermite weighted essentially non-oscillatory (HWENO) scheme with gradient reconstruction for the compressible Navier--Stokes equations. The key idea is to construct weak-derivative moments for each dissipative variable by taking the arithmetic average of two HWENO interface traces, which are then processed componentwise through the standard scalar nonlinear HWENO reconstruction formula, with only the input moments differing from those used in the convective reconstruction, thereby supplying the gradients needed for the viscous fluxes.  This reconstruction achieves fifth-order accuracy using only candidate polynomials of degree at most four, avoids the order reduction that would result from direct differentiation, and introduces negligible additional computational cost. To ensure robustness, especially in extreme test cases involving strong shocks or low densities, we apply a positivity-preserving limiter to the conservative states while retaining the nonlinear HWENO reconstruction for the viscous gradients. This combination is essential for stable computations in flows with strong gradients. Both the conservative and gradient reconstructions share the same compact stencils and candidate-polynomial structure, enabling a single implementation routine without introducing additional algorithmic complexity. Numerical results confirm that the proposed scheme delivers fifth-order accuracy and high resolution, offers competitive computational efficiency, and robustly resolves challenging finite-Reynolds-number Navier--Stokes flows in which the positivity-preserving limiter is actively engaged.

\vspace{.2in}
\vfill
\vfill {\bf Key Words:} Compressible Navier--Stokes equations, HWENO scheme, Gradient reconstruction, Positivity-preserving, Finite-volume method

{\bf AMS(MOS) subject classification}: 65M60, 35L65

\newpage
\baselineskip=2pc
\section{Introduction}

The compressible Navier--Stokes equations model viscous, heat-conducting gas flows. Their convective part is nonlinear and hyperbolic, whereas the viscous stress and heat flux introduce diffusion through velocity and internal-energy gradients. The relative importance of these mechanisms changes with the Reynolds number. At large Reynolds numbers, shocks, contact discontinuities, and thin shear layers require a low-dissipation treatment of convection; at finite Reynolds numbers, the viscous gradients must be resolved with comparable accuracy. A useful high-order discretization must therefore control oscillations in the convective part while retaining accuracy in the viscous part \cite{SjogreenYee2003}.

Many high-order methods have been developed for solving convection terms, such as the WENO schemes
\cite{js,bgsAO,hs,ZStetra,lpr,ZQd,ZSWENOMR}, the DG methods \cite{CockburnLinShu1989,CockburnShu1998,ZhangJCP2017,LinChanTomas2023,OFDG2,OEDG}, and the HWENO schemes \cite{QSHw1,ZhaoChenQiu2020HybridHWENO,LiShuQiu2021MRHWENO,FanZhaoXiongQiu2023,LuoBaum}. {Among them, HWENO schemes leverage the advantages of both WENO and DG methods.} {In the spatial reconstruction, they utilize not only the information of function values but also that of immediate neighboring cells, thereby avoiding excessively wide stencils while also preventing the use of excessive information within the target cell, which would otherwise lead to a reduced CFL number.} However, when dealing with convection-dominated problems, the numerical solutions may contain discontinuities. The first HWENO scheme \cite{QSHw1} has to employ two sets of candidate reconstruction stencils to circumvent the discontinuities, yet their robustness still falls short of that of the WENO scheme of the same type \cite{js}. To improve the robustness of the scheme, Zhao et al. \cite{ZhaoChenQiu2020HybridHWENO} drew inspiration from the limiter techniques in DG methods and proposed a hybrid HWENO scheme that modifies the first-order moments near discontinuities before performing the spatial reconstruction, which effectively enhances the robustness of the scheme. This idea was later extended to achieve unified reconstruction stencils \cite{FanQiuZhao2024}, further improving the computational efficiency and practicality of the HWENO scheme. Therefore, in this work, we adopt the HWENO scheme with this unified reconstruction stencil for the treatment of the convective terms.

{When the WENO schemes \cite{FanZhangQiu2022PPWENONS,TanChengShu2025}, DG methods \cite{ZhangJCP2017,LinChanTomas2023,LiuBuzzardZhang2024}, and HWENO method \cite{FanZhangQiu2021PPhybridHWENO} are extended to the compressible Navier--Stokes equations, the treatment of the viscous fluxes becomes a central difficulty.} One direct approach, adopted in the HWENO schemes of \cite{FanZhangQiu2021PPhybridHWENO}, differentiates the polynomial reconstructed for the convective fluxes. This approach is compact and reuses the nonlinear HWENO reconstruction, but direct differentiation reduces the formal order of accuracy of the gradient approximation by one in smooth regions. {For nonlinear degenerate parabolic equations, Liu et al. instead approximated the second derivative directly in conservative flux-difference form, obtaining sixth- and eighth-order finite-difference WENO schemes with high-order, nonoscillatory performance for nonsmooth solutions {\cite{LiuShuZhang2011}}.} Other finite-volume approaches have addressed this issue through compact high-order reconstruction of viscous fluxes \cite{CuetoFelguerosoColominasNogueira2007,ZhangLiuShu2012}, though these typically require wider stencils or a separate discretization of the derivative variables. Alternatively, gradients can be introduced as auxiliary variables following the LDG methodology \cite{TianXuKuertenVegt2015}. Such treatments can restore the required order of accuracy, but may demand additional reconstruction data. Moreover, while linear gradient approximations are often adequate for smooth flows, their robustness can deteriorate in low-density and low-pressure computations involving large gradients, even when positivity-preserving (PP) limiters are applied \cite{FanZhangQiu2022PPWENONS}. Despite these efforts, a viscous discretization that simultaneously preserves a compact stencil, retains full design order without resorting to higher-degree polynomials, and remains robust in severe flows is still a challenge.

In this work, we address this challenge by proposing a fifth-order finite-volume HWENO scheme with gradient reconstruction (HWENO-GR) for the compressible Navier--Stokes equations. The key idea is to construct weak-derivative moments for each dissipative variable by taking the arithmetic average of two HWENO interface traces, thereby supplying the gradients needed for the viscous fluxes. This construction is motivated by the weak-derivative formulation in the local discontinuous Galerkin (LDG) framework \cite{TianXuKuertenVegt2015}, as used for KdV-type equations in \cite{HybridLDGHWENO}. However, our scheme is not an LDG discretization: it neither evolves an auxiliary gradient system nor employs alternating numerical fluxes for the derivative variables. Instead, the weak-derivative moments are processed componentwise through the scalar nonlinear HWENO reconstruction formula, with only the input moments differing from those used in the convective reconstruction. Consequently, the gradient reconstruction introduces no additional algorithmic complexity. The resulting reconstruction uses candidate polynomials of degree at most four and supplies fifth-order viscous gradients at both face and volume quadrature points, avoiding the order loss inherent in direct differentiation.

For comparison, we also consider a direct-reconstruction formulation, denoted by the HWENO-DR scheme, in which the viscous gradients are obtained by differentiating the conservative candidates and reusing the convective nonlinear weights. This treatment is based on the direct-differentiation idea used in the hybrid HWENO scheme \cite{FanZhangQiu2021PPhybridHWENO}, with its implementation adapted to the present reconstruction. The resulting construction is efficient but suffers from the order loss described above. Furthermore, thanks to the componentwise processing of the weak-derivative moments in the HWENO-GR scheme, the computational cost of HWENO-GR is comparable to that of the HWENO-DR scheme. In fact, for some extreme test cases, the HWENO-GR scheme is even faster, primarily because its improved robustness allows for larger time steps. The proposed HWENO-GR scheme is efficient and robust, delivering full fifth-order accuracy without additional computational cost or algorithmic complexity. The scheme can also accommodate arbitrary positive linear weights that sum to unity and a compact reconstruction stencil, as in the HWENO scheme with unified stencils \cite{FanQiuZhao2024}. In contrast, the hybrid HWENO scheme \cite{FanZhangQiu2021PPhybridHWENO} is restricted to fixed linear weights, which limits it to fourth-order accuracy in two dimensions and makes it difficult to extend to unstructured meshes. For extreme flows involving strong shocks or expansions, which may create low-density or low-pressure states, we apply a PP-limiter \cite{FanZhangQiu2021PPhybridHWENO} to the conservative states before the numerical fluxes are evaluated, while retaining the nonlinear HWENO reconstruction for the viscous gradients. This combination ensures robust computations in flows with strong gradients. Numerical results confirm that the scheme delivers fifth-order accuracy and high resolution, offers competitive computational efficiency, and robustly resolves challenging finite-Reynolds-number Navier--Stokes flows in which the PP-limiter is actively engaged.

The remainder of this paper is organized as follows. Sections~\ref{sec2:hweno-gr-1d} and \ref{sec3:hweno-gr-2d} present the one- and two-dimensional HWENO-GR schemes, respectively, while their flowchart is provided in Section \ref{Flowchart}. Numerical results are reported in Section~ \ref{Numerical tests}, followed by conclusions in Section~ \ref{Concluding}.

\section{One-dimensional HWENO-GR scheme}\label{sec2:hweno-gr-1d}\label{sec:gradient-reconstruction}\label{sec:gradient-1d}

This section presents the one-dimensional formulation of the HWENO-GR scheme. Subsection 2.1 gives the semi-discrete moment formulation, the convective reconstruction, and the HWENO-DR scheme. Subsections 2.2 and 2.3 present the proposed gradient reconstruction and its accuracy analysis. Subsection 2.4 compares the two gradient treatments and gives the Runge--Kutta discretization.

\subsection{Convective discretization and semi-discrete formulation}\label{sec2:1dcase}

Consider the one-dimensional compressible Navier--Stokes equations in the form
\begin{equation}\label{sec2:1dHCLS}
	\begin{cases}
		\mathbf{u}_t+f(\mathbf{u},\mathbf{s})_x=0,\\
		\mathbf{u}(x,0)=\mathbf{u}_0(x),\\
		\mathbf{s}(x,t)=(\mu_x,e_x)^T,
	\end{cases}
\end{equation}
where \(\mathbf{u}=(\rho,m,E)^T\), \(m=\rho \mu\), \(\rho\) is the density, \(\mu\) is the velocity, and \(E\) is the total energy. The total flux is split as \(f=f^a-f^d\), with
\[
f^a=
\begin{pmatrix}
\rho \mu\\ \rho \mu^2+p\\ (E+p)\mu
\end{pmatrix},
\qquad
f^d=
\begin{pmatrix}
0\\ \tau\\ \tau \mu-q
\end{pmatrix}.
\]
For an ideal gas, the specific internal energy and pressure are
\(e=E/\rho-\mu^2/2\) and \(p=(\gamma-1)\rho e\), respectively. The viscous stress and heat flux are
\(\tau=\eta \mu_x/\mathrm{Re}\) and \(q=-\gamma e_x/(\mathrm{Pr}\,\mathrm{Re})\), with \(\eta=4/3\).
{Here \(\mathrm{Re}\) denotes the Reynolds number, and \(\mathrm{Pr}\) is the dimensionless Prandtl number. The symbol \(\mu\) denotes the streamwise velocity rather than the dynamic viscosity. Unless otherwise stated, \(\mathrm{Pr}=0.72\), corresponding to air, is used; Example~\ref{Example:ViscousShockTube} uses \(\mathrm{Pr}=0.73\).}
Thus the dissipative variables whose derivatives are needed in the viscous flux are \(\mu\) and \(e\).

Let the computational domain be partitioned into uniform cells
\(I_i=[x_{i-\frac12},x_{i+\frac12}]\), with cell center
\(x_i=(x_{i-\frac12}+x_{i+\frac12})/2\), cell size
\(\dx=x_{i+\frac12}-x_{i-\frac12}\), and scaled coordinate
\(\xi_i=(x-x_i)/\dx\). On this mesh, the finite-volume unknowns are the
cell average and the scaled first-order moment,
$
\bar{\mathbf{u}}_i=\frac1{\dx}\int_{I_i}\mathbf{u}(x,t)\,\mathrm{d}x \ \text{and} \ \bar{\mathbf{v}}_i=\frac1{\dx}\int_{I_i}\mathbf{u}(x,t)\xi_i\,\mathrm{d}x,
$
respectively. The projection of \eqref{sec2:1dHCLS} onto
\(\operatorname{span}\{1,\xi_i\}\), followed by integration by parts, gives
the two moment balances. Replacing the physical face fluxes by the numerical
flux \(\hat f_{i\pm\frac12}\) and evaluating the volume term with the
four-point Legendre--Gauss--Lobatto rule yields the conservative
semi-discrete scheme
\begin{equation}\label{sec2:1dSemi_HWENO}
	\left\{
	\begin{aligned}
		\frac{\mathrm{d}\bar{\mathbf{u}}_i}{\mathrm{d}t}
		&=-\frac1{\dx}\bigl(\hat f_{i+\frac12}-\hat f_{i-\frac12}\bigr)
		\triangleq\mathcal{F}^1_i(\bar{\mathbf{u}},\bar{\mathbf{v}}),\\
		\frac{\mathrm{d}\bar{\mathbf{v}}_i}{\mathrm{d}t}
		&=-\frac1{2\dx}\bigl(\hat f_{i-\frac12}+\hat f_{i+\frac12}\bigr)
		+\frac1{\dx}\sum_{G=1}^{4}\hat\omega_G
		f(\mathbf{u},\mathbf{s})\big|_{\hat x_i^G,t}
		\triangleq\mathcal{F}^2_i(\bar{\mathbf{u}},\bar{\mathbf{v}}).
	\end{aligned}
	\right.
\end{equation}
The first equation is the conservative balance for the cell average. The second evolves
the first-order moment and contains a volume-flux contribution because the derivative of the
test function is nonzero. At an interface, we use the PP local
Lax--Friedrichs flux
\begin{align}\label{sec2:lf-flux-1d}
&\hat f_{i+\frac12}
=\frac12\left[
f(\mathbf{u}_{i+\frac12}^-,\mathbf{s}_{i+\frac12}^-)
+f(\mathbf{u}_{i+\frac12}^+,\mathbf{s}_{i+\frac12}^+)
-\beta_{i+\frac12}
\bigl(\mathbf{u}_{i+\frac12}^+-\mathbf{u}_{i+\frac12}^-\bigr)
\right], \\
&\beta_{i+\frac{1}{2}} > \max\limits_{\mathbf{u}_{i+\frac{1}{2}}^\pm, \mathbf{s}_{i+\frac{1}{2}}^\pm}\left[ |\mu| + \frac{1}{2\rho^2 e} \left( \sqrt{\rho^2 q^2 + 2\rho^2 e|\tau - p|^2} + \rho|q| \right) \right] \notag
\end{align}
where \(\beta_{i+\frac12}\) is chosen no smaller than the PP wave-speed bound in \cite{ZhangJCP2017}. The superscripts \(+\) and \(-\) denote the right and left traces at \(x_{i+\frac12}\).
The four quadrature points in \eqref{sec2:1dSemi_HWENO} are
\[
\hat x_i^1=x_{i-\frac12},\qquad
\hat x_i^2=x_{i-\frac{\sqrt5}{10}},\qquad
\hat x_i^3=x_{i+\frac{\sqrt5}{10}},\qquad
\hat x_i^4=x_{i+\frac12},
\]
with normalized weights \(\hat\omega_1=\hat\omega_4=1/12\) and
\(\hat\omega_2=\hat\omega_3=5/12\). Thus the endpoint quadrature values are the
interface traces, whereas the two interior values are used only in the volume-flux
quadrature of the first-order moment equation.

The states entering the convective flux are obtained by local characteristic
reconstruction \cite{js,QSHw1}. For \(k=1,2,3\), let
\(\boldsymbol{\ell}^{k}\) denote the left eigenvector associated with the \(k\)-th
characteristic field of the convective-flux Jacobian. The projected cell average
and scaled first-order moment in \(I_i\) are
\(\bar u_i^k=(\boldsymbol{\ell}^{k})^{T}\bar{\mathbf u}_i\) and
\(\bar v_i^k=(\boldsymbol{\ell}^{k})^{T}\bar{\mathbf v}_i\), respectively.
The scalar reconstruction is performed independently for each \(k\), and the
reconstructed characteristic values are then transformed back to conservative variables.

For notational simplicity, the superscript \(k\) is omitted below.  Hence,
\(\bar u_i\) and \(\bar v_i\) denote the projected moments for a fixed characteristic
field. The HWENO reconstruction with unified stencils uses the values
\(\{\bar u_{i-1},\bar u_i,\bar u_{i+1},\bar v_{i-1},\bar v_{i+1}\}\) to define
three candidate polynomials \(p_m\), \(m=0,1,2\), which are then combined using the corresponding nonlinear weights
\(\omega_m\). The approximated polynomial and its modified first-order moment are reconstructed by
\begin{equation}\label{sec:solution-value-1d}
\left\{
\begin{aligned}
&u_i(x)
=\omega_0\left(
\frac1{\gamma_0}p_0(x)
-\frac{\gamma_1}{\gamma_0}p_1(x)
-\frac{\gamma_2}{\gamma_0}p_2(x)\right)
{}+\omega_1p_1(x)+\omega_2p_2(x),\\
&\widehat v_i
=\frac1{\dx}\int_{I_i}u_i(x)\frac{x-x_i}{\dx}\,\mathrm{d}x.
\end{aligned}
\right.
\end{equation}
The algorithm for constructing the polynomials \(p_m\) and computing the corresponding nonlinear weights \(\omega_m\).
is given in the Appendix \ref{app:nonlinear-hweno1}. In particular, when evaluating the second relation, the nonlinear weights \(\omega_m\) are constant over
\(I_i\) and thus remain outside the integral, leaving only the candidate polynomials \(p_m(x)\) to be integrated. The resulting modified moment \(\widehat v_i\) is then only used in the Runge--Kutta time discretization \eqref{sec2:1dRK3_HWENO}. Finally, applying \eqref{sec:solution-value-1d} in each characteristic field and transforming back yields the required conservative point values and modified moments. The polynomial reconstruction and moment modification share a single set of candidate stencils, while the positive linear weights \(\gamma_m\) maintain a unit sum. Following
\cite{ZhangSReview,FanZhangQiu2021PPhybridHWENO},
{a PP-limiter rescales the reconstructed conservative point values about the cell average at the four Legendre--Gauss--Lobatto points in \eqref{sec2:1dSemi_HWENO}. It preserves the cell average and positive density and pressure, and it does not act on the modified first-order moment \(\widehat v_i\) or on the gradient reconstruction described below.}

The comparison implementation denoted by the HWENO-DR scheme in this paper adapts the direct
differentiation idea of \cite{FanZhangQiu2021PPhybridHWENO} to the viscous
derivatives. It uses the same characteristic-wise reconstruction as that for the convective
fluxes, but obtains the viscous derivatives from differentiated solution polynomials.
In particular,  the nonlinear weights \(\omega_m\)
are reused. Let \(p'_m=\mathrm{d}p_m/\mathrm{d}x\).
For each characteristic field, the derivative is taken in the form
\begin{equation}\label{sec:dr-gradient-value-1d}
u^{\mathrm{DR}}_x(x)
=\omega_0\left(
\frac1{\gamma_0}p'_0(x)
-\frac{\gamma_1}{\gamma_0}p'_1(x)
-\frac{\gamma_2}{\gamma_0}p'_2(x)\right)
{}+\omega_1p'_1(x)
+\omega_2p'_2(x).
\end{equation}
After transformation back to conservative variables, the resulting values of
\((\rho_x,m_x,E_x)^T\) are converted to the gradients of the dissipative variables
\(\mu_x\) and \(e_x\). The HWENO-DR scheme provides a direct treatment of the viscous terms, where the reuse above
eliminates the separate evaluation of nonlinear weights for the differentiated reconstructed polynomials. Its main limitation, however, is that differentiating fifth-order solution polynomials generally produces only fourth-order gradients.

\subsection{Nonlinear HWENO reconstruction of the viscous gradients}\label{sec:gradient-reconstruction-1d}

To recover fifth-order accuracy, the HWENO-GR scheme reconstructs the viscous
derivatives from weak-derivative moments, rather than by differentiating the solution
polynomials. Let \(w\) denote either \(\mu\) or \(e\).
At an interface, the value used in the weak-derivative relation is the arithmetic average of the two HWENO traces,
\begin{equation}\label{sec:interface-average-1d}
\widehat w_{i+\frac12}
=\frac12\left(w_{i+\frac12}^-+w_{i+\frac12}^+\right).
\end{equation}
The auxiliary relation \(r=w_x\) is used to generate cell moments for the derivative.
Multiplying this relation by a test function \(\phi\in\mathbbm{P}^1(I_i)\), integrating over
\(I_i\), and replacing the boundary values of \(w\) by the interface averages in
\eqref{sec:interface-average-1d}, we obtain
\begin{equation}\label{sec:weak-derivative-1d}
\int_{I_i}r(x)\phi(x)\,\mathrm{d}x
=-\int_{I_i}w(x)\phi_x(x)\,\mathrm{d}x
{}+\widehat w_{i+\frac12}\phi(x_{i+\frac12}^-)
-\widehat w_{i-\frac12}\phi(x_{i-\frac12}^+),
\ \forall \phi\in\mathbbm{P}^1(I_i).
\end{equation}
The weak relation \eqref{sec:weak-derivative-1d} has the same local
integration-by-parts form as the auxiliary-gradient relation in LDG methods
\cite{TianXuKuertenVegt2015,TianXuYangVegt2025,HybridLDGHWENO}. The LDG
discretization rewrites the governing equation as a first-order system and approximates
the solution and its auxiliary gradient in discontinuous polynomial spaces coupled by
numerical fluxes. Here \eqref{sec:weak-derivative-1d} is used only to form cell moments
of \(r\) from the HWENO approximation of \(w\). The derivative \(r\) is not introduced
as an independent discrete unknown, nor is an LDG flux pair required. Instead, its moments are directly supplied to the finite-volume HWENO reconstruction.

Taking \(\phi=1\) and \(\phi=\xi_i\) in \eqref{sec:weak-derivative-1d} gives
\begin{equation}\label{sec:weak-moments-1d}
\bar r_i^{\,0}
=\frac{\widehat w_{i+\frac12}-\widehat w_{i-\frac12}}{\dx},
\qquad
\bar r_i^{\,1}
=\frac{\frac12(\widehat w_{i+\frac12}+\widehat w_{i-\frac12})-\bar w_i}{\dx}.
\end{equation}
In \eqref{sec:weak-moments-1d},
\(\bar w_i=\sum_{G=1}^{4}\hat\omega_G
w\bigl(\mathbf{u}_i(\hat x_i^G)\bigr)\) is computed from the conservative states supplied
by the convective HWENO reconstruction at the four quadrature points specified in
Subsection~\ref{sec2:1dcase}. The quantities \(\bar r_i^{\,0}\) and
\(\bar r_i^{\,1}\) are the zeroth-order and first-order weak moments of \(r\) on
\(I_i\), respectively. For \(w\in\{\mu,e\}\), the derivative at each quadrature
point \(\hat x_i^G\), \(G=1,\ldots,4\), is reconstructed as
\begin{equation}\label{sec:grad-value-1d}
w_x(x)
=\omega_0\left(
\frac1{\gamma_0}p_0(x)
-\frac{\gamma_1}{\gamma_0}p_1(x)
-\frac{\gamma_2}{\gamma_0}p_2(x)\right)
{}+\omega_1p_1(x)+\omega_2p_2(x).
\end{equation}
The candidate polynomials \(p_m\) and nonlinear weights \(\omega_m\) in
\eqref{sec:grad-value-1d} are obtained by applying the HWENO reconstruction in
Appendix~\ref{app:nonlinear-hweno1} to the weak-moment data
\(\bar r_k^{\,0}\) and \(\bar r_k^{\,1}\) on the corresponding stencils, where the positive linear weights
\(\gamma_m\) have unit sum.
Evaluating \eqref{sec:grad-value-1d} with \(w=\mu\) and \(w=e\) gives \(\mu_x\)
and \(e_x\), respectively, at the interface and volume quadrature points. From these derivatives, we obtain
\(\tau=\eta\mu_x/\mathrm{Re}\) and
\(q=-\gamma e_x/(\mathrm{Pr}\,\mathrm{Re})\) in the dissipative flux. Consequently, the convective and viscous fluxes are evaluated at the same locations, while the
viscous derivatives are reconstructed from weak-derivative moments rather than obtained
by differentiating the reconstructed  polynomials.

\subsection{Fifth-order accuracy of the compatible gradient reconstruction}
\label{sec:gradient-accuracy}

The order reduction of the HWENO-DR scheme stems from direct polynomial differentiation: if
\(p=u+O(\dx^5)\), then
\(p'=u_x+O(\dx^4)\), making direct differentiation of a fifth-order reconstruction only fourth-order accurate.
Fan et al.
\cite{FanZhangQiu2021PPhybridHWENO} avoid this loss in their hybrid scheme by using a
sixth-order approximation in the smooth case, whose derivative is fifth-order accurate, thus preserving overall fifth-order accuracy. Without this enhancement, direct differentiation suffers the same order reduction as the HWENO-DR scheme.

The HWENO-GR scheme does not differentiate the solution polynomial. Its weak-derivative formulas
divide by \(\dx\), so fifth-order derivative moments follow if the cell and interface
averages in their numerators are sixth-order accurate. For a smooth solution, let
\(\mathbf u_i(x)\) denote the conservative
reconstruction in \(I_i\). Conservation gives
\[
\frac1{\dx}\int_{I_i}
\bigl(\mathbf{u}_i(x)-\mathbf{u}(x)\bigr)\,\mathrm{d}x=0.
\]
Since \(\mathbf u_i-\mathbf u=O(\dx^5)\), expansion about \(x_i\) yields
\[
w(\mathbf{u}_i)-w(\mathbf{u})
=\nabla_{\mathbf{u}}w(\mathbf{u}(x_i))
\bigl(\mathbf{u}_i-\mathbf{u}\bigr)+O(\dx^6),
\]
and the leading term vanishes after cell averaging. With the degree-five quadrature,
\begin{equation}\label{sec:primitive-average-accuracy}
\bar w_i=\bar w_i^{\,\mathrm{exact}}+O(\dx^6),
\end{equation}
where ``exact'' denotes the exact cell average. At \(x_{i+\frac12}\), the two
large-stencil traces satisfy
\[
\begin{aligned}
p_0^-(x_{i+\frac12})
&=w(x_{i+\frac12})
-\frac{241}{95760}\dx^5w^{(5)}(x_{i+\frac12})+O(\dx^6),\\
p_0^+(x_{i+\frac12})
&=w(x_{i+\frac12})
+\frac{241}{95760}\dx^5w^{(5)}(x_{i+\frac12})+O(\dx^6),
\end{aligned}
\]
so the leading fifth-order errors cancel in their arithmetic average. The smooth-region
nonlinear weights preserve this cancellation \cite{FanQiuZhao2024}, giving
\begin{equation}\label{sec:interface-average-accuracy}
\widehat w_{i+\frac12}=\frac12\left(w_{i+\frac12}^-+w_{i+\frac12}^+\right)=w(x_{i+\frac12})+O(\dx^6).
\end{equation}
Substituting \eqref{sec:primitive-average-accuracy} and \eqref{sec:interface-average-accuracy} into \eqref{sec:weak-moments-1d} shows that both weak-derivative moments are fifth-order accurate. Consequently, their nonlinear reconstruction yields fifth-order derivative values at the interface and volume quadrature points. The sixth-order estimates mentioned above are intermediate consequences of conservation and interface-error cancellation; the final gradient approximation remains fifth-order accurate.

\subsection{Comparison and time discretization}\label{sec:time-discretization-1d}

\begin{figure}[H]
\centering
\includegraphics[width=0.77\textwidth]{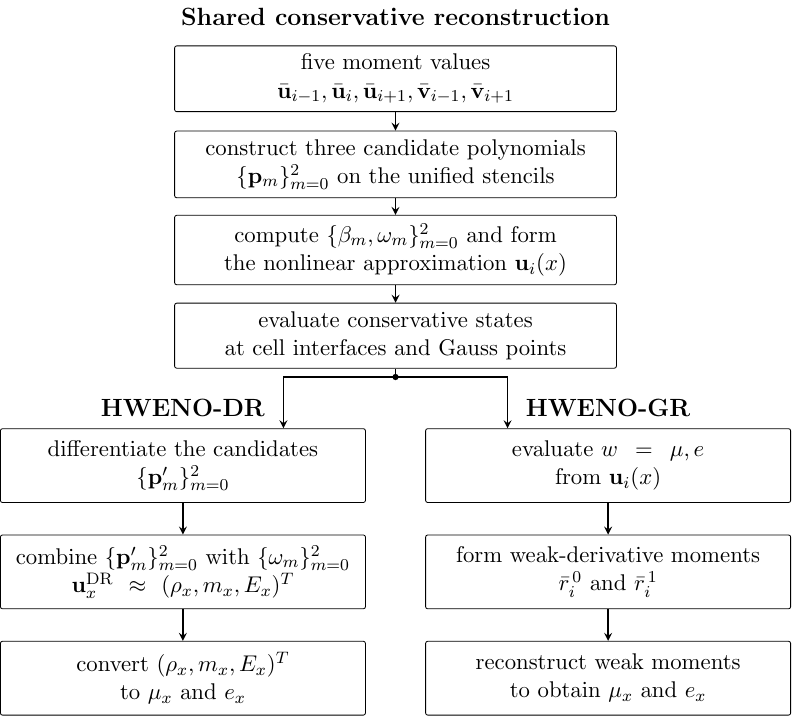}
\caption{Shared conservative reconstruction and alternative one-dimensional
viscous-gradient treatments in the HWENO-DR and HWENO-GR schemes.}
\label{fig:dr-gr-gradient-comparison}
\end{figure}

In the HWENO-DR scheme used here, the differentiated conservative
candidates are combined using the nonlinear weights from the convective reconstruction. The HWENO-GR scheme
instead forms weak moments of
\(w_x\), \(w=\mu,e\), from the conservative values in
\eqref{sec:solution-value-1d} and reconstructs these moments directly. The common
convective reconstruction and the two viscous-gradient treatments are summarized in
Fig.~\ref{fig:dr-gr-gradient-comparison}.

For time discretizations, we use the third-order SSP Runge--Kutta method with the modified
first-order moments
\begin{equation}\label{sec2:1dRK3_HWENO}
	\left\{
	\begin{aligned}
		\begin{bmatrix}
			\bar{\mathbf{u}}^{(1)}_i \\ \bar{\mathbf{v}}^{(1)}_i
		\end{bmatrix}
		&=
		\begin{bmatrix}
			\bar{\mathbf{u}}^{n}_i \\ \hat{\mathbf{v}}^{n}_i
		\end{bmatrix}
		+\dt
		\begin{bmatrix}
			\mathcal{F}^1_i(\bar{\mathbf{u}}^{n},\bar{\mathbf{v}}^{n}) \\
			\mathcal{F}^2_i(\bar{\mathbf{u}}^{n},\bar{\mathbf{v}}^{n})
		\end{bmatrix},\\
		\begin{bmatrix}
			\bar{\mathbf{u}}^{(2)}_i \\ \bar{\mathbf{v}}^{(2)}_i
		\end{bmatrix}
		&=\frac34
		\begin{bmatrix}
			\bar{\mathbf{u}}^{n}_i \\ \hat{\mathbf{v}}^{n}_i
		\end{bmatrix}
		+\frac14\left(
		\begin{bmatrix}
			\bar{\mathbf{u}}^{(1)}_i \\ \hat{\mathbf{v}}^{(1)}_i
		\end{bmatrix}
		+\dt
		\begin{bmatrix}
			\mathcal{F}^1_i(\bar{\mathbf{u}}^{(1)},\bar{\mathbf{v}}^{(1)}) \\
			\mathcal{F}^2_i(\bar{\mathbf{u}}^{(1)},\bar{\mathbf{v}}^{(1)})
		\end{bmatrix}\right),\\
		\begin{bmatrix}
			\bar{\mathbf{u}}^{n+1}_i \\ \bar{\mathbf{v}}^{n+1}_i
		\end{bmatrix}
		&=\frac13
		\begin{bmatrix}
			\bar{\mathbf{u}}^{n}_i \\ \hat{\mathbf{v}}^{n}_i
		\end{bmatrix}
		+\frac23\left(
		\begin{bmatrix}
			\bar{\mathbf{u}}^{(2)}_i \\ \hat{\mathbf{v}}^{(2)}_i
		\end{bmatrix}
		+\dt
		\begin{bmatrix}
			\mathcal{F}^1_i(\bar{\mathbf{u}}^{(2)},\bar{\mathbf{v}}^{(2)}) \\
			\mathcal{F}^2_i(\bar{\mathbf{u}}^{(2)},\bar{\mathbf{v}}^{(2)})
		\end{bmatrix}\right).
	\end{aligned}
	\right.
\end{equation}
Here \(\hat{\mathbf{v}}_i^n\), \(\hat{\mathbf{v}}_i^{(1)}\), and
\(\hat{\mathbf{v}}_i^{(2)}\) are the modified first-order moments supplied by the
HWENO reconstruction \eqref{sec:solution-value-1d} at the
corresponding stages.

\section{Two-dimensional HWENO-GR scheme}\label{sec3:hweno-gr-2d}\label{sec2:2dcase}\label{sec:gradient-2d}

\subsection{Fully discrete finite-volume formulation}\label{sec3:2d-fully-discrete}

The two-dimensional dimensionless compressible Navier--Stokes equations are written as
\begin{equation}\label{sec2:2dHCLS}
\begin{cases}
\mathbf{u}_t+f(\mathbf{u},\mathbf{s})_x+g(\mathbf{u},\mathbf{s})_y=0,\\
\mathbf{u}(x,y,0)=\mathbf{u}_0(x,y),\\
\mathbf{s}(x,y,t)=(\mu_x,\mu_y,\nv_x,\nv_y,e_x,e_y)^T,
\end{cases}
\end{equation}
where \(\mathbf{u}=(\rho,\rho \mu,\rho \nv,E)^T\). The fluxes are
\[
\begin{aligned}
f(\mathbf{u},\mathbf{s}) &=
\begin{pmatrix}
\rho \mu\\
\rho \mu^2+p-\tau_{xx}/\mathrm{Re}\\
\rho \mu \nv-\tau_{xy}/\mathrm{Re}\\
(E+p)\mu-\left(\tau_{xx}\mu+\tau_{xy}\nv+\gamma e_x/\mathrm{Pr}\right)/\mathrm{Re}
\end{pmatrix},\\
g(\mathbf{u},\mathbf{s}) &=
\begin{pmatrix}
\rho \nv\\
\rho \mu \nv-\tau_{xy}/\mathrm{Re}\\
\rho \nv^2+p-\tau_{yy}/\mathrm{Re}\\
(E+p)\nv-\left(\tau_{xy}\mu+\tau_{yy}\nv+\gamma e_y/\mathrm{Pr}\right)/\mathrm{Re}
\end{pmatrix}.
\end{aligned}
\]
{Here \(\mu\) and \(\nu\) are the velocity components in the \(x\)- and \(y\)-directions, respectively.}
\(e=E/\rho-(\mu^2+\nv^2)/2\), and \(p=(\gamma-1)\rho e\). The viscous stresses are
\(\tau_{xx}=4\mu_x/3-2\nv_y/3\), \(\tau_{xy}=\mu_y+\nv_x\), and
\(\tau_{yy}=4\nv_y/3-2\mu_x/3\). Thus the dissipative variables reconstructed below are
\(\mu\), \(\nv\), and \(e\).

Let \(I_i=[x_{i-\frac12},x_{i+\frac12}]\) and
\(I_j=[y_{j-\frac12},y_{j+\frac12}]\) denote uniform cells in the \(x\)- and
\(y\)-directions, respectively, and set \(I_{i,j}=I_i\times I_j\),
with cell center \((x_i,y_j)\), mesh sizes
\(\dx=x_{i+\frac12}-x_{i-\frac12}\) and
\(\dy=y_{j+\frac12}-y_{j-\frac12}\), and scaled coordinates
\(\xi_i=(x-x_i)/\dx\) and \(\eta_j=(y-y_j)/\dy\).
On \(I_{i,j}\), the finite-volume unknowns are the cell average and the two
scaled first-order moments,
$
\bar{\mathbf{u}}_{i,j}
=\frac1{\dx\dy}\iint_{I_{i,j}}\mathbf{u}\,\mathrm{d}x\mathrm{d}y,\
\bar{\mathbf{v}}_{i,j}
=\frac1{\dx\dy}\iint_{I_{i,j}}\mathbf{u}\xi_i\,\mathrm{d}x\mathrm{d}y,\
\bar{\mathbf{w}}_{i,j}
=\frac1{\dx\dy}\iint_{I_{i,j}}\mathbf{u}\eta_j\,\mathrm{d}x\mathrm{d}y.
$
The projection of \eqref{sec2:2dHCLS} onto
\(\operatorname{span}\{1,\xi_i,\eta_j\}\), followed by integration by parts,
gives the three moment balances. Replacing the physical face fluxes by
numerical fluxes and applying the three-point Gauss rule to the face and volume
integrals yields the conservative semi-discrete system
\begin{equation}\label{sec2:2dSemi_HWENO}
	\left\{
	\begin{aligned}
		\frac{{\rm d} \bar{\mathbf{u}}_{i,j}}{{\rm d} t}={}&
		-\frac{1}{\dx}\sum_{G=1}^3\hat{\omega}_G
		(\hat{f}_{i+\frac12,G}-\hat{f}_{i-\frac12,G})
		-\frac{1}{\dy}\sum_{G=1}^3\hat{\omega}_G
		(\hat{g}_{G,j+\frac12}-\hat{g}_{G,j-\frac12}), \\
		\frac{{\rm d} \bar{\mathbf{v}}_{i,j}}{{\rm d} t}={}&
		-\frac{1}{2\dx}\sum_{G=1}^3\hat{\omega}_G
		(\hat{f}_{i-\frac12,G}+\hat{f}_{i+\frac12,G})
		-\frac{1}{\dy}\sum_{G=1}^3\hat{\omega}_G
		\frac{\hat{x}_i^G-x_i}{\dx}
		(\hat{g}_{G,j+\frac12}-\hat{g}_{G,j-\frac12})\\
		&+\frac{1}{\dx}\sum_{G,H=1}^3\hat{\omega}_G\hat{\omega}_H
		f(\mathbf{u},\mathbf{s})(\hat{x}_i^{G},\hat{y}_j^H), \\
		\frac{{\rm d} \bar{\mathbf{w}}_{i,j}}{{\rm d}t}={}&
		-\frac{1}{\dx}\sum_{G=1}^3\hat{\omega}_G
		\frac{\hat{y}_j^G-y_j}{\dy}
		(\hat{f}_{i+\frac12,G}-\hat{f}_{i-\frac12,G})
		-\frac{1}{2\dy}\sum_{G=1}^3\hat{\omega}_G
		(\hat{g}_{G,j-\frac12}+\hat{g}_{G,j+\frac12})\\
		&+\frac{1}{\dy}\sum_{G,H=1}^3\hat{\omega}_G\hat{\omega}_H
		g(\mathbf{u},\mathbf{s})(\hat{x}_i^{G},\hat{y}_j^H).
	\end{aligned}
	\right.
\end{equation}
The face fluxes
are evaluated by the PP local Lax--Friedrichs flux. At the vertical
and horizontal faces, respectively, they are
\begin{align}\label{sec2:lf-flux-2d}
&\hat f_{i+\frac12,G}
=\frac12\left[
f(\mathbf{u}_{i+\frac12,G}^-,\mathbf{s}_{i+\frac12,G}^-)
+f(\mathbf{u}_{i+\frac12,G}^+,\mathbf{s}_{i+\frac12,G}^+)
-\beta_{i+\frac12,j}
\left(\mathbf{u}_{i+\frac12,G}^+-\mathbf{u}_{i+\frac12,G}^-\right)
\right],\\
&\hat g_{G,j+\frac12}
=\frac12\left[
g(\mathbf{u}_{G,j+\frac12}^-,\mathbf{s}_{G,j+\frac12}^-)
+g(\mathbf{u}_{G,j+\frac12}^+,\mathbf{s}_{G,j+\frac12}^+)
-\beta_{i,j+\frac12}
\left(\mathbf{u}_{G,j+\frac12}^+-\mathbf{u}_{G,j+\frac12}^-\right)
\right],\notag\\
&\beta_{i+\frac12,j}>
\max_{\mathbf{u}_{i+\frac12,G}^{\pm},\,\mathbf{s}_{i+\frac12,G}^{\pm}}
\left[|\mu|+\frac{1}{2\rho^2e}\left(
\sqrt{\rho^2q_x^2+2\rho^2e
\left|\left(\frac{\tau_{xx}}{\mathrm{Re}}-p,\frac{\tau_{xy}}{\mathrm{Re}}\right)\right|^2}
+\rho|q_x|\right)\right],\notag\\
&\beta_{i,j+\frac12}>
\max_{\mathbf{u}_{G,j+\frac12}^{\pm},\,\mathbf{s}_{G,j+\frac12}^{\pm}}
\left[|\nv|+\frac{1}{2\rho^2e}\left(
\sqrt{\rho^2q_y^2+2\rho^2e
\left|\left(\frac{\tau_{xy}}{\mathrm{Re}},\frac{\tau_{yy}}{\mathrm{Re}}-p\right)\right|^2}
+\rho|q_y|\right)\right].\notag
\end{align}
Here \(q_x=\gamma e_x/(\mathrm{Pr}\,\mathrm{Re})\) and
\(q_y=\gamma e_y/(\mathrm{Pr}\,\mathrm{Re})\). Following
\cite{ZhangJCP2017}, each maximum is taken over the reconstructed states on both sides
of the three face quadrature points.
The three Gauss points in each coordinate direction are
\[
\begin{aligned}
\hat x_i^1&=x_{i-\frac{\sqrt{15}}{10}},\qquad
\hat x_i^2=x_i,\qquad
\hat x_i^3=x_{i+\frac{\sqrt{15}}{10}},\\
\hat y_j^1&=y_{j-\frac{\sqrt{15}}{10}},\qquad
\hat y_j^2=y_j,\qquad
\hat y_j^3=y_{j+\frac{\sqrt{15}}{10}},
\end{aligned}
\]
with normalized weights
\(\hat\omega_1=\hat\omega_3=5/18\) and \(\hat\omega_2=4/9\).
Accordingly, \(\hat f_{i\pm\frac12,G}\) is evaluated at
\((x_{i\pm\frac12},\hat y_j^G)\), \(\hat g_{G,j\pm\frac12}\) at
\((\hat x_i^G,y_{j\pm\frac12})\), and the volume fluxes at the tensor-product points
\((\hat x_i^G,\hat y_j^H)\).

To close \eqref{sec2:2dSemi_HWENO}, the face and volume states required by the
convective terms are obtained by local characteristic reconstruction \cite{js,QSHw}.
For \(k=1,2,3,4\), let \(\boldsymbol{\ell}^{k}\) denote the left eigenvector
associated with the \(k\)-th characteristic field of the relevant convective-flux
Jacobian. The projected cell average and scaled first-order moments in cell \(I_{i,j}\) are
$
\bar u_{i,j}^k=(\boldsymbol{\ell}^{k})^{T}\bar{\mathbf u}_{i,j},
\bar v_{i,j}^k=(\boldsymbol{\ell}^{k})^{T}\bar{\mathbf v}_{i,j},
\bar w_{i,j}^k=(\boldsymbol{\ell}^{k})^{T}\bar{\mathbf w}_{i,j}.
$
The scalar reconstruction is performed independently for each \(k\), and the
reconstructed characteristic values are then transformed back to conservative variables.

For notational simplicity, the superscript \(k\) is omitted below. Hence,
\(\bar u_{i,j}\), \(\bar v_{i,j}\), and \(\bar w_{i,j}\) denote, respectively,
the projected cell average and the scaled first-order moments in the \(x\)- and
\(y\)-directions for a fixed characteristic field. The HWENO reconstruction with the unified stencils uses the
cell averages
\(\{\bar u_{i-1,j-1},\allowbreak\bar u_{i,j-1},\allowbreak
\bar u_{i+1,j-1},
\bar u_{i-1,j},\allowbreak\bar u_{i,j},\allowbreak
\bar u_{i+1,j},
\bar u_{i-1,j+1},\allowbreak\bar u_{i,j+1},\allowbreak
\bar u_{i+1,j+1}\}\), the first-order moments \(\{\bar v_{i-1,j},\allowbreak\bar v_{i+1,j},\allowbreak
\bar v_{i,j-1},\allowbreak\bar v_{i,j+1},\bar w_{i-1,j},\allowbreak\bar w_{i+1,j},\allowbreak
\bar w_{i,j-1},\allowbreak\bar w_{i,j+1}\}\) to define a
quartic candidate \(p_0(x,y)\) and four linear candidates \(p_m(x,y)\),
\(m=1,\ldots,4\), respectively, which are then combined using the corresponding nonlinear weights
\(\omega_m\). The approximated polynomial and its modified first-order moments are reconstructed by
\begin{equation}\label{sec:solution-value-2d}
\left\{
\begin{aligned}
&u_{i,j}(x,y)
=\omega_0\left(
\frac1{\gamma_0}p_0(x,y)
-\sum_{m=1}^4\frac{\gamma_m}{\gamma_0}p_m(x,y)\right)
{}+\sum_{m=1}^4\omega_m p_m(x,y),\\
&\widehat v_{i,j}
=\frac1{\dx\dy}\iint_{I_{i,j}}u_{i,j}(x,y)
\frac{x-x_i}{\dx}\,\mathrm{d}x\mathrm{d}y,\\
&\widehat w_{i,j}
=\frac1{\dx\dy}\iint_{I_{i,j}}u_{i,j}(x,y)
\frac{y-y_j}{\dy}\,\mathrm{d}x\mathrm{d}y.
\end{aligned}
\right.
\end{equation}
The algorithm for constructing the polynomials \(p_m(x,y)\) and computing the corresponding nonlinear weights \(\omega_m\).
is given in the Appendix \ref{app:nonlinear-hweno2}. In particular, when evaluating the second relation, the nonlinear weights \(\omega_m\) are constant over
\(I_i\) and thus remain outside the integral, leaving only the candidate polynomials \(p_m(x,y)\) to be integrated. The resulting modified moments \(\widehat v_{i,j}\) and
\(\widehat w_{i,j}\)  are then only used in the Runge--Kutta time discretization. Finally, applying \eqref{sec:solution-value-2d} in each characteristic field and transforming back yields the required conservative point values and modified moments. Also, the polynomial  reconstruction and the moment modification share a single set of candidate stencils, while the positive linear weights \(\gamma_m\)  maintain a unit sum. Following
\cite{ZhangSReview,FanZhangQiu2021PPhybridHWENO}, {the reconstructed conservative states are checked at the face Gauss points and the volume tensor-product Gauss points. If limiting is required, only the reconstructed values at the face Gauss points are modified before the numerical fluxes are evaluated. The volume values are checked but not modified, and the modified first-order moments and viscous gradients remain unchanged.}

As in the one-dimensional case, the HWENO-DR scheme
differentiates the candidate polynomials from the characteristic-wise convective
reconstruction. In particular, the nonlinear
weights \(\omega_m\) are reused, and for each characteristic field, the derivatives are
\begin{equation}\label{sec:dr-gradient-value-2d}
\left\{
\begin{aligned}
u_x^{\mathrm{DR}}(x,y)
&=\omega_0\left(
\frac1{\gamma_0}\partial_x p_0(x,y)
-\sum_{m=1}^4\frac{\gamma_m}{\gamma_0}\partial_x p_m(x,y)\right)
{}+\sum_{m=1}^4\omega_m\partial_x p_m(x,y),\\
u_y^{\mathrm{DR}}(x,y)
&=\omega_0\left(
\frac1{\gamma_0}\partial_y p_0(x,y)
-\sum_{m=1}^4\frac{\gamma_m}{\gamma_0}\partial_y p_m(x,y)\right)
{}+\sum_{m=1}^4\omega_m\partial_y p_m(x,y).
\end{aligned}
\right.
\end{equation}
After transformation back to conservative variables, the results are then converted to
$\mu_x$, $\mu_y$, $\nv_x$, $\nv_y$, $e_x$, and $e_y$.

For time discretizations, the semi-discrete system is advanced by the two-dimensional analogue of
\eqref{sec2:1dRK3_HWENO}, using the evolved cell average and first-order moments together
with the reconstructed modified moments. The HWENO-DR scheme reuses the reconstructed polynomials and
nonlinear weights of the convective reconstruction and is therefore inexpensive. However, its
limitation is that differentiating a fifth-order solution reconstruction generally produces
only fourth-order gradients.

\subsection{Nonlinear HWENO reconstruction of the viscous gradients}\label{sec:gradient-reconstruction-2d}

To recover fifth-order accuracy, the HWENO-GR scheme reconstructs the viscous
derivatives from weak-derivative moments rather than differentiating the solution
polynomials. We now construct the derivatives in
\(\mathbf{s}=(\mu_x,\mu_y,\nv_x,\nv_y,e_x,e_y)^T\), which are needed both in the
numerical face fluxes and in the volume quadrature terms of
\eqref{sec2:2dSemi_HWENO}. Let \(w\in\{\mu,\nv,e\}\) denote a dissipative variable.
The face averages below are used as boundary data for the weak-derivative relation.
At a vertical or horizontal face, this boundary trace is taken as the arithmetic
average of the two HWENO traces:
\[
\widehat w_{i+\frac12,j}(y)
=\frac12\left(w_{i+\frac12,j}^-(y)+w_{i+\frac12,j}^+(y)\right),
\qquad
\widehat w_{i,j+\frac12}(x)
=\frac12\left(w_{i,j+\frac12}^-(x)+w_{i,j+\frac12}^+(x)\right).
\]
The values at the internal tensor-product Gauss points are obtained from the
cellwise HWENO reconstruction:
\[
\widehat w_{i,j}^{G,H}
=w(\hat x_i^G,\hat y_j^H),
\qquad G,H=1,2,3.
\]
Denote the face moments of \(w\) by \(\mathcal{W}^{\ell,k}\), where
\(\ell=x,y\) identifies the face-normal direction and \(k=0,1\) is the moment order.
The zeroth-order and transverse first-order moments are
\begin{equation}\label{sec:2d-face-moments}
\begin{aligned}
\mathcal{W}^{x,0}_{i+\frac12,j}
&=\frac1{\dy}\int_{I_j}\widehat w_{i+\frac12,j}(y)\,\mathrm{d}y,
&
\mathcal{W}^{x,1}_{i+\frac12,j}
&=\frac1{\dy}\int_{I_j}
\widehat w_{i+\frac12,j}(y)\eta_j\,\mathrm{d}y,\\
\mathcal{W}^{y,0}_{i,j+\frac12}
&=\frac1{\dx}\int_{I_i}\widehat w_{i,j+\frac12}(x)\,\mathrm{d}x,
&
\mathcal{W}^{y,1}_{i,j+\frac12}
&=\frac1{\dx}\int_{I_i}
\widehat w_{i,j+\frac12}(x)\xi_i\,\mathrm{d}x .
\end{aligned}
\end{equation}
The required quantities in \eqref{sec:2d-face-moments} are evaluated from the same
three face quadrature points used in \eqref{sec2:2dSemi_HWENO}. To obtain cell data
for the derivatives, introduce the auxiliary relations \(r^x=w_x\) and \(r^y=w_y\).
Multiplying them by a test function \(\phi\in\mathbbm{P}^1(I_{i,j})\), integrating
over \(I_{i,j}\), and replacing the boundary traces by the interface averages above
gives
\[
\begin{aligned}
\iint_{I_{i,j}}r^x\phi\,\mathrm{d}x\mathrm{d}y
&=
\int_{I_j}\widehat w_{i+\frac12,j}\phi(x_{i+\frac12},y)\,\mathrm{d}y
-\int_{I_j}\widehat w_{i-\frac12,j}\phi(x_{i-\frac12},y)\,\mathrm{d}y
-\iint_{I_{i,j}}w\phi_x\,\mathrm{d}x\mathrm{d}y,\\
\iint_{I_{i,j}}r^y\phi\,\mathrm{d}x\mathrm{d}y
&=
\int_{I_i}\widehat w_{i,j+\frac12}\phi(x,y_{j+\frac12})\,\mathrm{d}x
-\int_{I_i}\widehat w_{i,j-\frac12}\phi(x,y_{j-\frac12})\,\mathrm{d}x
-\iint_{I_{i,j}}w\phi_y\,\mathrm{d}x\mathrm{d}y,
\end{aligned}
\]
for all \(\phi\in\mathbbm{P}^1(I_{i,j})\). As in the one-dimensional case, these weak
relations are used only to generate cell moments for the subsequent HWENO gradient
reconstruction, while no auxiliary LDG system or LDG flux pair is required.

Taking the three test functions
\(\phi=1,\xi_i,\eta_j\) produces the zeroth-order weak moment and the two directional
first-order weak moments of each derivative. For \(r^x\), these moments are
\begin{equation}\label{sec:2d-x-weak-moments}
\begin{aligned}
\bar r_{ij}^{x,0}
&=\frac{\mathcal{W}^{x,0}_{i+\frac12,j}
-\mathcal{W}^{x,0}_{i-\frac12,j}}{\dx},\\
\bar r_{ij}^{x,\xi}
&=\frac{\frac12\left(\mathcal{W}^{x,0}_{i+\frac12,j}
+\mathcal{W}^{x,0}_{i-\frac12,j}\right)-\bar w_{ij}}{\dx},\\
\bar r_{ij}^{x,\eta}
&=\frac{\mathcal{W}^{x,1}_{i+\frac12,j}-\mathcal{W}^{x,1}_{i-\frac12,j}}{\dx},
\end{aligned}
\end{equation}
where the superscripts \(0\), \(\xi\), and \(\eta\) identify the corresponding test
functions. Interchanging the coordinate directions gives the moments of \(r^y\):
\begin{equation}\label{sec:2d-y-weak-moments}
\begin{aligned}
\bar r_{ij}^{y,0}
&=\frac{\mathcal{W}^{y,0}_{i,j+\frac12}
-\mathcal{W}^{y,0}_{i,j-\frac12}}{\dy},\\
\bar r_{ij}^{y,\xi}
&=\frac{\mathcal{W}^{y,1}_{i,j+\frac12}-\mathcal{W}^{y,1}_{i,j-\frac12}}{\dy},\\
\bar r_{ij}^{y,\eta}
&=\frac{\frac12\left(\mathcal{W}^{y,0}_{i,j+\frac12}
+\mathcal{W}^{y,0}_{i,j-\frac12}\right)-\bar w_{ij}}{\dy}.
\end{aligned}
\end{equation}
In \eqref{sec:2d-x-weak-moments} and \eqref{sec:2d-y-weak-moments},
\(\bar w_{ij}=\sum_{G,H=1}^3\hat\omega_G\hat\omega_H
\widehat w_{i,j}^{G,H}\) is computed from the conservative HWENO reconstruction
at the tensor-product quadrature points.
Equations \eqref{sec:2d-x-weak-moments} and \eqref{sec:2d-y-weak-moments}
therefore convert the HWENO interface traces of \(w\) into cell moments of its two
derivatives. These weak-derivative moments serve as the input data for the nonlinear
HWENO reconstruction below.

For each derivative direction \(\ell=x,y\), these weak moments are passed to the
same two-dimensional nonlinear HWENO reconstruction used for the convective variables.
The reconstruction is cellwise. After it is formed, it is evaluated at all quadrature
points required by \eqref{sec2:2dSemi_HWENO}. At such a point, the reconstructed
directional derivative is written as
\begin{equation}\label{sec:grad-value-2d}
w_{\ell}(x,y)
=\omega_0\left(
\frac1{\gamma_0}p_0(x,y)
-\sum_{m=1}^4\frac{\gamma_m}{\gamma_0}p_m(x,y)\right)
{}+\sum_{m=1}^4\omega_m p_m(x,y),
\qquad \ell=x,y.
\end{equation}
The candidate polynomials \(p_m(x,y)\) and nonlinear weights \(\omega_m\) in
\eqref{sec:grad-value-2d} are obtained by applying the HWENO reconstruction in
Appendix~\ref{app:nonlinear-hweno2} to the weak-moment data
\(\bar r_k^{\ell,0}\), \(\bar r_k^{\ell,\xi}\), and
\(\bar r_k^{\ell,\eta}\) on the corresponding stencils, where the positive linear weights
\(\gamma_m\) have unit sum.
Finally, taking \(w=\mu,\nv,e\)
supplies all velocity and internal-energy derivatives entering the stress tensor and
heat flux.

The accuracy proof is the same as in Subsection~\ref{sec:gradient-accuracy}. In
particular, the reconstructed cell average satisfies
\begin{equation}\label{sec:primitive-average-accuracy-2d}
\bar w_{i,j}=\bar w_{i,j}^{\,\mathrm{exact}}+O(\dx^6).
\end{equation}
The same weak-moment argument then gives
\begin{equation}\label{sec:gradient-consistency-2d}
w_{\ell}(x,y)=\partial_{\ell}w(x,y)+O(\dx^5),
\qquad \ell=x,y,\quad w\in\{\mu,\nv,e\}.
\end{equation}
Thus the velocity and internal-energy derivatives used in the face and volume flux
quadratures are fifth-order accurate at the corresponding quadrature points.

\begin{rem}\label{rem:nonlinear-gradient-reconstruction}
The nonlinear HWENO reconstruction of the weak-derivative moments in \eqref{sec:grad-value-2d} is retained for
strongly compressible flows. Its linear counterpart is adequate for smooth solutions, but it may fail to be reliable in certain extreme-flow tests considered here. For example, in the two-dimensional Sedov problem, the linear gradient reconstruction produced very large PP wave-speed bounds \(\beta_{i+\frac12,j}\) and \(\beta_{i,j+\frac12}\). As a result, the admissibility restarts
reduced the time step to approximately \(10^{-10}\), making it practically impossible to reach the prescribed final time.
This diagnostic is reported in Fig.~\ref{sec3:Sedov2DBeta}. In contrast, the nonlinear HWENO reconstruction substantially reduces these bounds and is therefore used for the viscous gradients throughout the computations.
\end{rem}

\section{Flowchart of the HWENO-GR scheme}
\label{Flowchart}

\begin{figure}[H]
\centering
\includegraphics[width=0.92\textwidth]{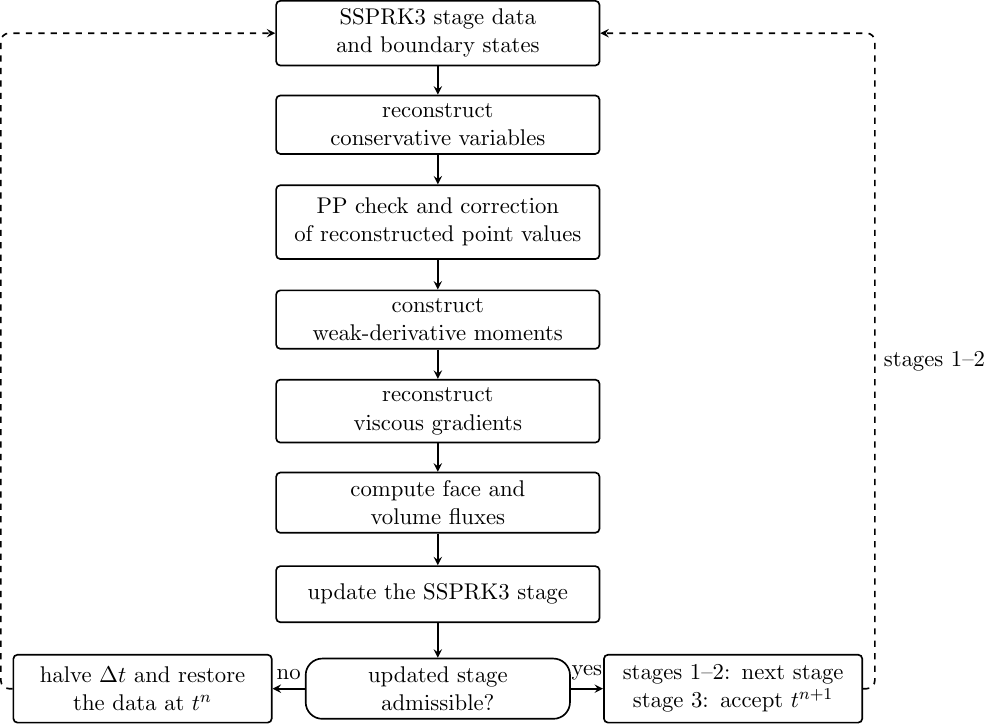}
\caption{{Stage-wise SSPRK3 procedure with PP correction and time-step control.}}
\label{fig:ssprk3-workflow}
\end{figure}

{In this section, we present the flowchart of the HWENO-GR scheme with PP limiters in the  one- and two-dimensional cases, as shown in Fig.~\ref{fig:ssprk3-workflow}, and the boundary treatment is given in Remark~\ref{rem:boundary-treatment}. The reconstructions described above define the spatial operator evaluated at each SSPRK3 stage. At a given stage, the conservative variables are reconstructed and positivity is checked at the required face and volume points. If correction is required, the PP limiter changes only the reconstructed face Gauss values. The weak-derivative moments and viscous gradients are then reconstructed, followed by flux evaluation and the SSPRK3 stage update. An inadmissible updated stage is rejected and restarted from \(t^n\) with \(\Delta t/2\); otherwise, stages 1 and 2 feed the next stage, while stage 3 yields \(t^{n+1}\).
}

\begin{rem}\label{rem:boundary-treatment}
{Boundary data are assigned in ghost cells. At a stationary reflective or slip wall, the normal velocity is extended oddly, while the tangential velocity and \(\rho\), \(p\), and \(e\) are extended evenly. At a stationary adiabatic no-slip wall, both velocity components are odd and \(\rho\), \(p\), and \(e\), and hence the temperature, are even; the total energy is then recomputed. Symmetry boundaries use the reflective parity, inflow states are prescribed, and outflow states are copied from the adjacent interior cells.}
\end{rem}

\section{Numerical tests}
\label{Numerical tests}

In this section, we examine the fifth-order accuracy, computational cost, and robustness
of the HWENO-GR scheme using smooth and non-smooth finite-Reynolds-number tests. The
HWENO-GR scheme and the HWENO-DR scheme use the same convective discretization and differ
only in the evaluation of the viscous gradients: the former reconstructs the gradients
from weak-derivative moments, whereas the latter differentiates the reconstructed polynomials.
All other discretization components and parameters are identical, so the comparison
isolates the viscous discretization. Their Euler-limit solutions agree up to round-off,
and the non-smooth comparisons are restricted to finite Reynolds numbers.

The convective and diffusive time-step estimates follow
\cite{ZhangJCP2017,FanZhangQiu2021PPhybridHWENO}. The PP-limiter uses
\(\varepsilon_{\mathrm{PP}}=10^{-12}\), and the CFL number is \(0.6\).
On a uniform mesh, let \(h=\dx\) in
one dimension and \(h=\min\{\dx,\dy\}\) in two dimensions. With
\(c=\sqrt{\gamma p/\rho}\), \(\alpha^*\) is the face maximum of \(|\mu|+c\) in one
dimension and of \(\max\{|\mu|+c,|\nv|+c\}\) in two dimensions. We take
\(\dt_{\mathrm{trial}}=\min\{\mathrm{CFL}\,h/\alpha^*,\,b\,\mathrm{Re}\,h^2\}\),
with \(b=0.001\). If a Runge--Kutta stage contains nonpositive density or pressure after
limiting, the step is rejected and repeated with half the trial step. The viscous gradients
are recomputed from the accepted stage data without separate limiting.

The Fortran 95 computations were run on an Intel(R) Xeon(R) Gold 6130 CPU at 2.10 GHz
with identical compiler settings. The CPU times represent relative costs for this
implementation. The reported PP-limiter percentage is the fraction of cell-stage checks
that require a nontrivial correction.

\subsection{Accuracy tests}\label{sec3:AccuracyTests}
The accuracy tests verify the fifth-order convergence of the HWENO-GR scheme and compare its errors and CPU times with those of the HWENO-DR scheme. The HWENO-GR scheme is more efficient and attains the designed fifth-order accuracy, whereas the HWENO-DR scheme exhibits only fourth-order accuracy.

\begin{example}\label{Example:NavierStokes1DTestOrder}
		We consider the one-dimensional Navier--Stokes system with \(Re=100\). A manufactured source term is added:
		\begin{equation*}
			\frac{\partial}{\partial{t}}\begin{bmatrix}	\rho \\ \rho\mu\\E \end{bmatrix} +
			\frac{\partial}{\partial{x}}\begin{bmatrix}	\rho\mu \\ \rho\mu^2+p-\tau\\ \mu(E+p)-(\mu \tau-q) \end{bmatrix} = \mathbf{S}(x,t),
		\end{equation*}
		where $\rho$ is the density, $\mu$ is the velocity, $E$ is the total energy, $p$ is the pressure, $\tau$ is the shear stress tensor, and $q$ denotes the heat diffusion flux. The source term $\mathbf{S}(x,t)$ is obtained analytically by substituting the following manufactured solution into the governing equations:
		\begin{equation*}
			\begin{aligned}
				\rho(x,t)&=2+0.2\exp(-t)\sin(2\pi x),\\
				\mu(x,t)&=1+0.2\exp(-t)\cos(2\pi x),\\
				e(x,t)&=1+0.2\exp(-t)\sin(2\pi x),
			\end{aligned}
		\end{equation*}
		with $E=\rho(e+\frac{1}{2}\mu^2)$ and $\gamma=1.4$. The computation is performed on the domain $[0,1]$ with periodic boundary conditions up to $T=0.1$. The numerical errors and CPU times are presented in Table \ref{sec3:1dNavierStokesTest}. The HWENO-GR scheme attains fifth-order accuracy on the refined meshes. The HWENO-DR scheme shows an observed rate close to four in this test. The corresponding error--CPU-time comparison is shown in Fig.~\ref{sec3:1dCPU_LinftyL1_1dNavierStokes}.
		\begin{table}[t]
			\centering
			\small
			\begin{threeparttable}
				\caption{Example \ref{Example:NavierStokes1DTestOrder}. One-dimensional Navier--Stokes equations: $L^\infty$ and $L^1$ errors, orders and CPU times of the HWENO-DR scheme and the HWENO-GR scheme.}
				\label{sec3:1dNavierStokesTest}
				{\begin{tabular}{cccccccccccccc}
						\toprule
						{Meshes} & ${L^\infty}$ error&Order&${L^1}$ error&Order &CPU \cr
						\midrule
						{HWENO-DR scheme}\\
						 40&	1.47E-06&   $- $&    5.50E-07&   $- $	&	1.56E-02	\\
						 80&	8.20E-09&   7.49&    3.35E-09&   7.36	&	4.69E-02	\\
						120&	1.29E-09&   4.56&    5.24E-10&   4.58	&	1.41E-01	\\
						160&	3.62E-10&   4.42&    1.47E-10&   4.42	&	2.97E-01	\\
						200&	1.37E-10&   4.35&    5.67E-11&   4.27	&	5.47E-01	\\
						240&	6.27E-11&   4.29&    2.61E-11&   4.26	&	8.91E-01	\\
						280&	3.25E-11&   4.26&    1.36E-11&   4.23	&	1.36E+00	\\
						\midrule
						{HWENO-GR scheme}\\
						 40&	1.66E-06&   $- $&    5.15E-07&   $- $	&	9.37E-03	\\
						 80&	8.26E-09&   7.65&    2.57E-09&   7.65	&	5.62E-02	\\
						120&	8.68E-10&   5.56&    3.25E-10&   5.11	&	1.59E-01	\\
						160&	2.08E-10&   4.97&    7.58E-11&   5.06	&	3.37E-01	\\
						200&	6.78E-11&   5.02&    2.47E-11&   5.03	&	6.11E-01	\\
						240&	2.71E-11&   5.02&    9.88E-12&   5.02	&	9.94E-01	\\
						280&	1.25E-11&   5.03&    4.53E-12&   5.05	&	1.48E+00	\\
						\bottomrule	
				\end{tabular}}
			\end{threeparttable}
		\end{table}
		\begin{figure}[!htbp]
			\centering
			\subfigure{\includegraphics[width=0.41\textwidth,angle=0]{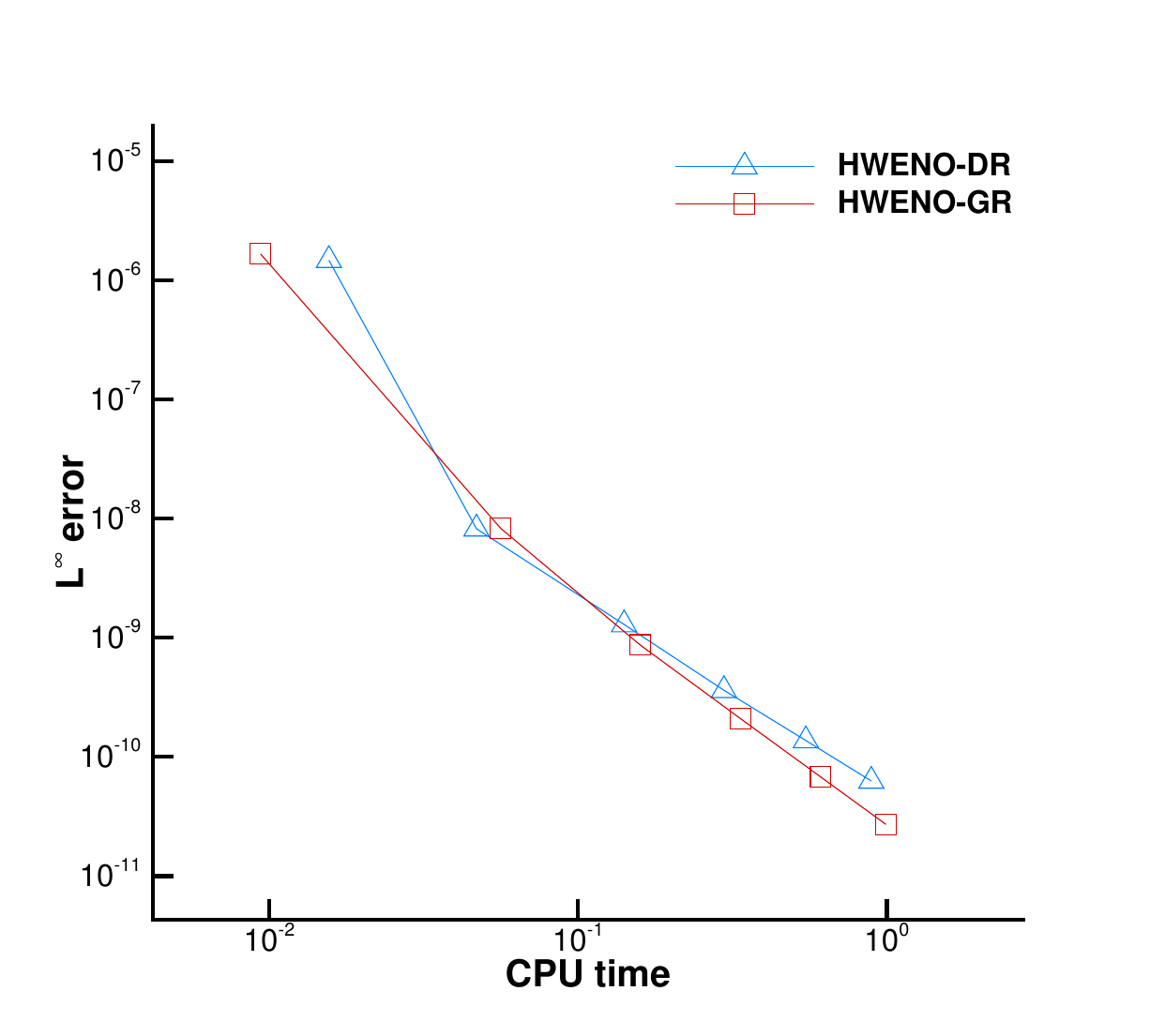}}
			\subfigure{\includegraphics[width=0.41\textwidth,angle=0]{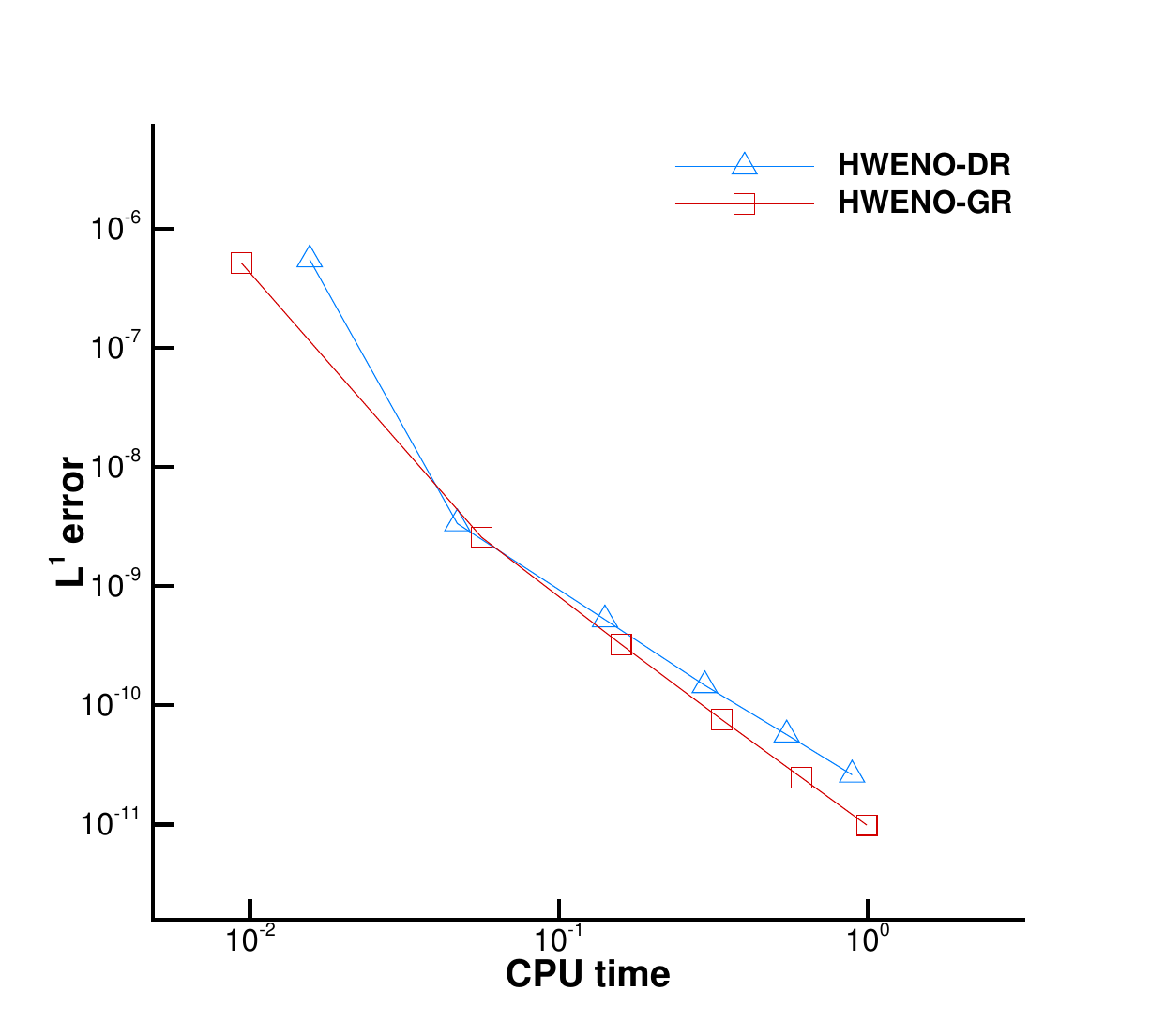}}
			\caption{Comparison of $L^\infty$, $L^1$ errors and CPU time for Example \ref{Example:NavierStokes1DTestOrder}.}
			\label{sec3:1dCPU_LinftyL1_1dNavierStokes}
		\end{figure}
\end{example}

\begin{example}\label{Example:NavierStokes2DTestOrder}
		We solve two-dimensional compressible Navier--Stokes equations with a manufactured source term and {\(\mathrm{Re}=100\)}:	
		\begin{equation*}
			\frac{\partial}{\partial{t}}\begin{bmatrix}	\rho \\ \rho\mu \\ \rho \nv  \\E \end{bmatrix} +
			\frac{\partial}{\partial{x}}\begin{bmatrix}	\rho\mu \\ \rho\mu^2+p-\frac{\tau_{xx}}{{\mathrm{Re}}}\\ \rho\mu \nv-\frac{\tau_{yx}}{{\mathrm{Re}}}\\ \mu(E+p)-\frac{\tau_{xx}\mu+\tau_{yx}\nv+\frac{\gamma}{{\mathrm{Pr}}}e_x}{{\mathrm{Re}}} \end{bmatrix}+
			\frac{\partial}{\partial{y}}\begin{bmatrix}	\rho \nv \\ \rho\mu \nv-\frac{\tau_{xy}}{{\mathrm{Re}}} \\ \rho \nv^2+p-\frac{\tau_{yy}}{{\mathrm{Re}}}\\ \nv(E+p)-\frac{\tau_{xy}\mu+\tau_{yy}\nv+\frac{\gamma}{{\mathrm{Pr}}}e_y}{{\mathrm{Re}}} \end{bmatrix} = \mathbf{S}(x,y,t),
		\end{equation*}
		where $\rho$ is the density, $\mu$ and \(\nv\) are the velocities in the $x$- and $y$-directions, respectively, $E$ is the total energy, $e$ is the internal energy, and $p$ is the pressure. The source term $\mathbf{S}(x,y,t)$ is obtained analytically by substituting the following manufactured solution into the governing equations:
		\begin{equation*}
			\begin{aligned}
				\rho(x,y,t)&=2+0.1\exp(-t)\Big[\sin(2\pi(2x-y))+\cos(2\pi(-x+2y))\Big],\\
				\mu(x,y,t)&=1+0.1\exp(-t)\Big[\sin(2\pi (-x+y))+\cos(2\pi (x+y))\Big],\\
				\nv(x,y,t)&=2+0.1\exp(-t)\Big[\sin(2\pi (x+y))+\cos(2\pi (x-y))\Big],\\
				e(x,y,t)&=2+0.1\exp(-t)\Big[\sin(2\pi (x-y))+\cos(2\pi (x+y))\Big],
			\end{aligned}
		\end{equation*}
		The manufactured solution deliberately combines nonseparable mixed-direction modes
		and uses different modal compositions for the primitive variables. This choice reduces
		the possibility that grid symmetry or cancellation of leading errors produces an
		artificially high observed convergence rate.
		Here, $E=\rho(e+\frac{1}{2}(\mu^2+\nv^2))$ and $\gamma=1.4$. The computation is performed on $[0,1]\times[0,1]$ with periodic boundary conditions in both directions up to $T=0.1$.
		
		The $L^\infty$ and $L^1$ errors, together with CPU times, are reported in Table \ref{sec3:2dNavierStokesTest}. On the refined meshes, the HWENO-GR scheme shows fifth-order convergence in both norms. For this mixed-direction manufactured solution, the results provide numerical evidence that the compatible gradient reconstruction retains fifth-order accuracy for the viscous terms in the complete two-dimensional Navier--Stokes discretization. For the HWENO-DR scheme, the observed order falls below five once the mesh is sufficiently refined and tends toward fourth-order accuracy, in agreement with the direct-differentiation accuracy discussed above. The error--CPU-time comparison is shown in Fig.~\ref{sec3:1dCPU_LinftyL1_2dNavierStokes}.
		\begin{table}[ht]
			\centering
			\small
			\begin{threeparttable}
				\caption{Example \ref{Example:NavierStokes2DTestOrder}. Two-dimensional Navier--Stokes equations: $L^\infty$ and $L^1$ errors, orders and CPU times of the HWENO-DR scheme and the HWENO-GR scheme.}
				\label{sec3:2dNavierStokesTest}
				{\begin{tabular}{ccccccccccccc}
						\toprule
						{Meshes} & ${L^\infty}$ error&Order&${L^1}$ error&Order &CPU \cr
						\midrule
						{HWENO-DR scheme}\\
						  $40\times40$&	4.12E-04&   $- $&    3.36E-05&   $- $	&	1.31E+01	\\
						  $80\times80$&	1.50E-06&   8.10&    1.11E-07&   8.24	&	1.51E+02	\\
						$120\times120$&	5.72E-08&   8.06&    8.76E-09&   6.26	&	6.50E+02	\\
						$160\times160$&	7.73E-09&   6.96&    1.64E-09&   5.82	&	1.88E+03	\\
						$200\times200$&	2.12E-09&   5.80&    5.06E-10&   5.27	&	4.25E+03	\\
						$240\times240$&	9.13E-10&   4.62&    2.12E-10&   4.77	&	8.25E+03	\\
						$280\times280$&	4.53E-10&   4.55&    1.07E-10&   4.44	&	1.46E+04	\\
						\midrule
						{HWENO-GR scheme}\\
						  $40\times40$&	3.96E-06&   $- $&    6.23E-07&   $- $	&	1.08E+01	\\
						  $80\times80$&	8.09E-08&   5.61&    1.86E-08&   5.07	&	1.22E+02	\\
						$120\times120$&	1.05E-08&   5.04&    2.45E-09&   5.00	&	5.44E+02	\\
						$160\times160$&	2.47E-09&   5.03&    5.83E-10&   5.00	&	1.56E+03	\\
						$200\times200$&	8.01E-10&   5.05&    1.91E-10&   5.00	&	3.56E+03	\\
						$240\times240$&	3.21E-10&   5.02&    7.67E-11&   5.00	&	6.96E+03	\\
						$280\times280$&	1.48E-10&   5.01&    3.55E-11&   5.00	&	1.24E+04	\\
						\bottomrule
				\end{tabular}}
			\end{threeparttable}
		\end{table}
		\begin{figure}[!htbp]
			\centering
			\subfigure{\includegraphics[width=0.41\textwidth,angle=0]{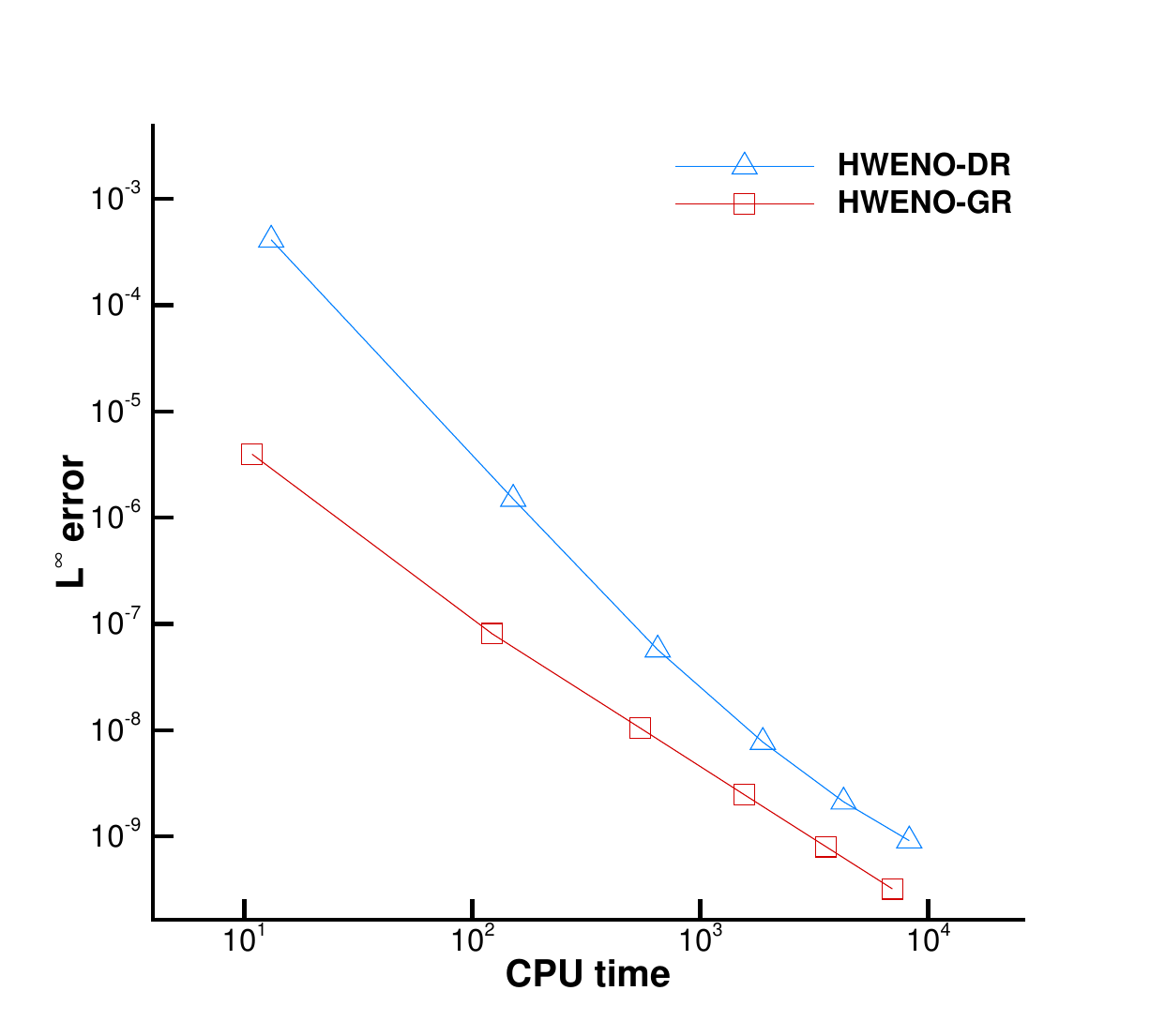}}
			\subfigure{\includegraphics[width=0.41\textwidth,angle=0]{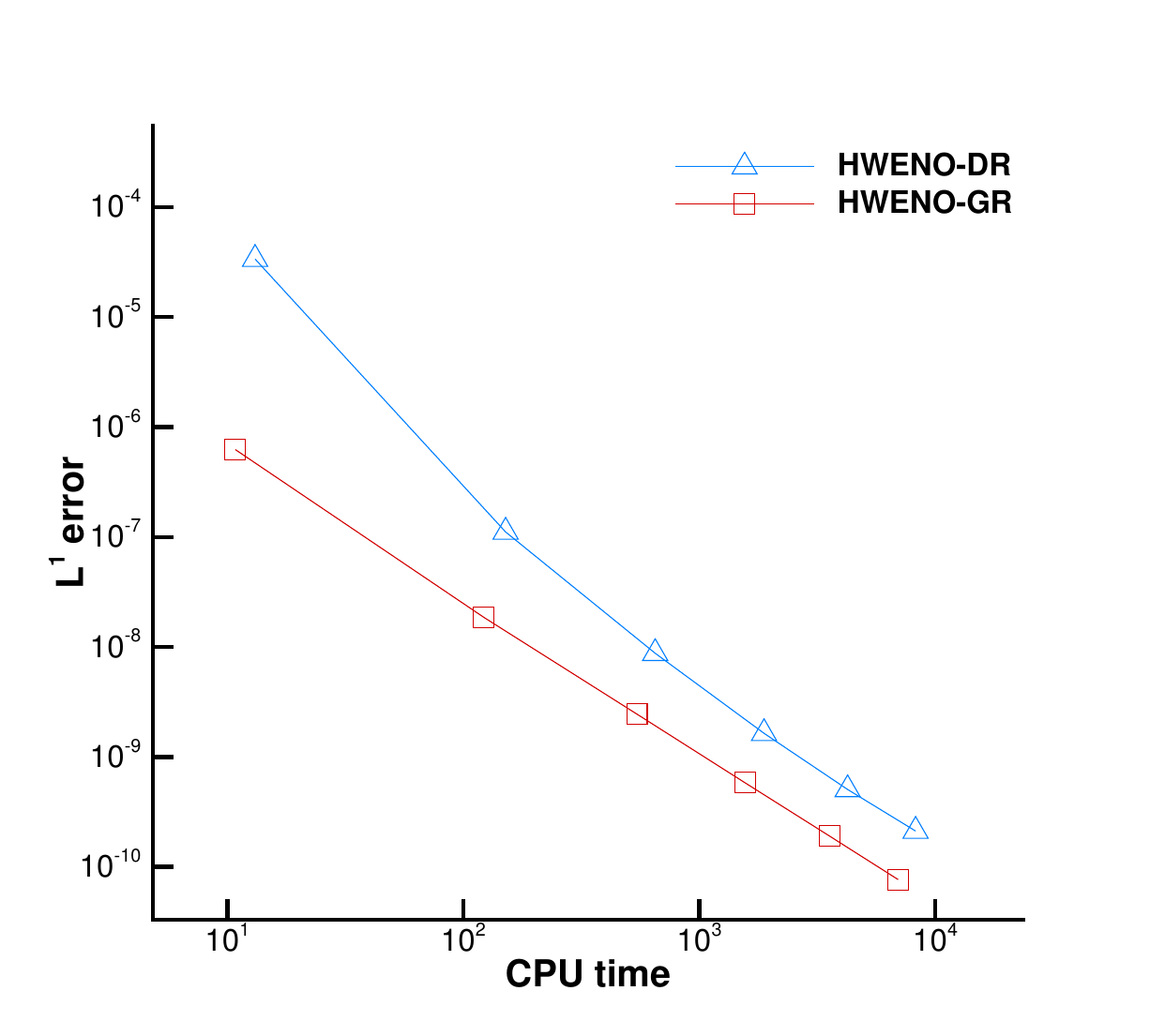}}
			\caption{Comparison of $L^\infty$, $L^1$ errors and CPU time for Example \ref{Example:NavierStokes2DTestOrder}.}
			\label{sec3:1dCPU_LinftyL1_2dNavierStokes}
		\end{figure}
\end{example}

\subsection{Non-smooth tests}\label{sec3:NonsmoothTests}
The HWENO-DR scheme and the HWENO-GR scheme are compared on finite-Reynolds-number
benchmark problems. Both schemes use the same convective discretization and quadrature rules,
differing only in the treatment of the viscous gradients. Specifically, the HWENO-DR scheme differentiates the reconstructed polynomials and reuses the convective nonlinear weights, whereas the HWENO-GR scheme reconstructs the gradients from weak-derivative moments.
Table~\ref{sec3:CPUtime} reports the total CPU time, the percentage of PP-limiter
activations, and \(T_{\mathrm{GR}}/T_{\mathrm{DR}}\), for which a value below unity
indicates a lower cost for the HWENO-GR scheme. The two schemes have comparable costs: the
HWENO-GR scheme is faster in most \(Re=1000\) tests and remains close to the HWENO-DR
scheme in the other viscous cases.

Across the non-smooth tests, the two schemes give similar numerical solutions.
The numerical tests show that both schemes remain stable for extreme
flows. However, the HWENO-GR scheme recovers the designed fifth-order accuracy in smooth regions without increasing algorithmic complexity, whereas the HWENO-DR scheme, which directly differentiates the reconstructed solution polynomials, achieves only fourth-order accuracy and does not exploit the full high-order accuracy of the underlying reconstruction.
\begin{table}[!htbp]
	\centering
	\small
	\setlength{\tabcolsep}{2.5pt}
	\begin{threeparttable}
	 \caption{Total CPU time in seconds, percentage of PP-limiter activations, and CPU-time ratio \(T_{\mathrm{GR}}/T_{\mathrm{DR}}\) for the finite-Reynolds-number non-smooth examples.}
		\label{sec3:CPUtime}
		{\begin{tabular}{L{0.35\textwidth}ccccc}
			\toprule
			\multirow{2}{*}{ Numerical example}&
			\multicolumn{2}{c}{HWENO-DR scheme}&
			\multicolumn{2}{c}{HWENO-GR scheme}&
			\multirow{2}{*}{\(T_{\mathrm{GR}}/T_{\mathrm{DR}}\)}\cr
			\cmidrule(lr){2-3} \cmidrule(lr){4-5}
			& CPU time & PP-limiter & CPU time & PP-limiter & \\
			\midrule
			\multicolumn{6}{l}{$Re = 1000$} \\
			\ref{Example:Lax1D} 1D Lax
			&2.19E-01 &0.00E+00\% &1.87E-01  &0.00E+00\% &0.85 \\
			\ref{Example:DoubleRare} 1D Double Rarefaction
			&7.81E-02 &1.18E-02\% &4.69E-02 &8.69E-04\% &0.60 \\
			\ref{Example:Sedov1D} 1D Sedov
			&2.09E+00 &3.00E+00\% &1.95E+00 &3.00E+00\% &0.93 \\
			\ref{Example:Leblanc1D} 1D Leblanc
			&3.19E+01 &2.50E-01\% &3.14E+01 &2.50E-01\% &0.98 \\
			\ref{Example:Sedov2D} 2D Sedov
			&6.33E+03 &9.93E-02\% &3.18E+03 &3.25E-01\% &0.50\\
			\ref{Example:ShockDiffraction} 2D shock-diffraction
			&3.18E+03 &4.29E-03\% &3.37E+03 &6.54E-04\% &1.06\\
			\ref{Example:HM2000} 2D Mach-2000
			&1.45E+04 &3.84E-01\% &1.17E+04 &4.32E+00\% &0.80\\
			\ref{Example:Mach10shock} 2D Mach-10 shock
			&9.37E+03 &8.57E-03\% &8.67E+03 &2.40E-03\% &0.93\\
			\midrule
			\multicolumn{6}{l}{Other Reynolds numbers} \\
			\ref{Example:ViscousShockMixing} 2D viscous shock--mixing-layer problem
			&2.49E+03 &0.00E+00\% &2.62E+03 &0.00E+00\% &1.05\\
			\ref{Example:ViscousShockTube} 2D viscous shock-tube
			&1.25E+05 &0.00E+00\% &1.36E+05 &0.00E+00\% &1.09\\
			\bottomrule
		\end{tabular}}
	\end{threeparttable}
	\label{table-cpu}
\end{table}
\begin{example}\label{Example:Lax1D}
We solve the Lax problem for one-dimensional Navier--Stokes equations with the initial conditions:
\begin{equation*}
	(\rho,\mu,p,\gamma)^\mathsf{T}=\begin{cases}
		(0.445, 0.698, 3.528, 1.4)^\mathsf{T}, \ -5\le x<0,
		\\(0.5, 0, 0.571, 1.4)^\mathsf{T},\ ~~~~~~~~~~~0\le x\le 5.
	\end{cases}
\end{equation*}	
The final time is $T = 1.3$, and outflow boundary conditions are imposed on all boundaries. Figure~\ref{sec3:Lax1D_WENOHWENO} shows the density profiles for several finite Reynolds numbers. For all four Reynolds numbers, the profiles from the two schemes are similar and agree well with the exact solution.

\begin{figure}[!ht]
	\centering
	\subfigure[$Re = 100$]{\includegraphics[width=7.5cm,angle=0]{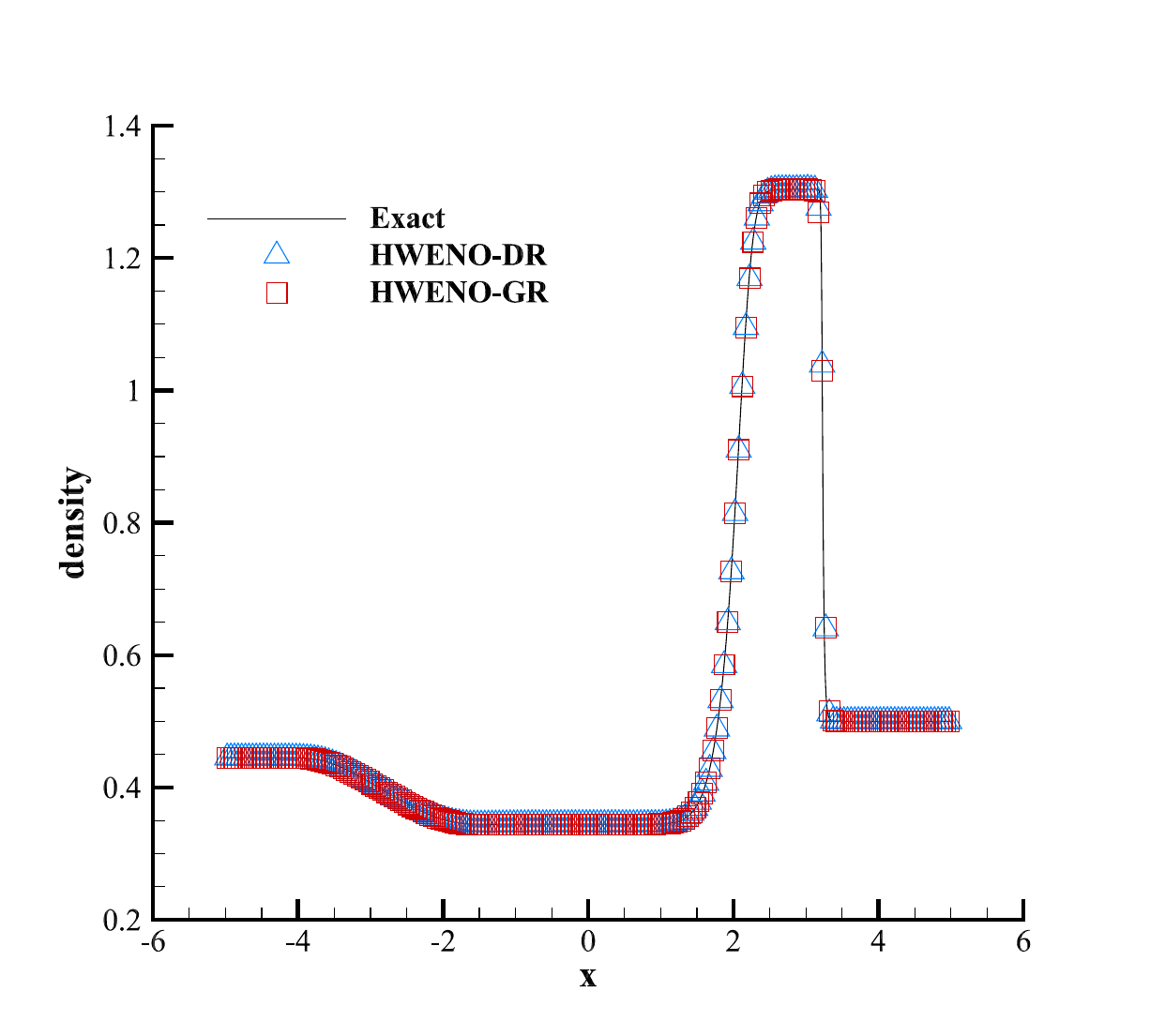}}
	\subfigure[$Re = 1000$]{\includegraphics[width=7.5cm,angle=0]{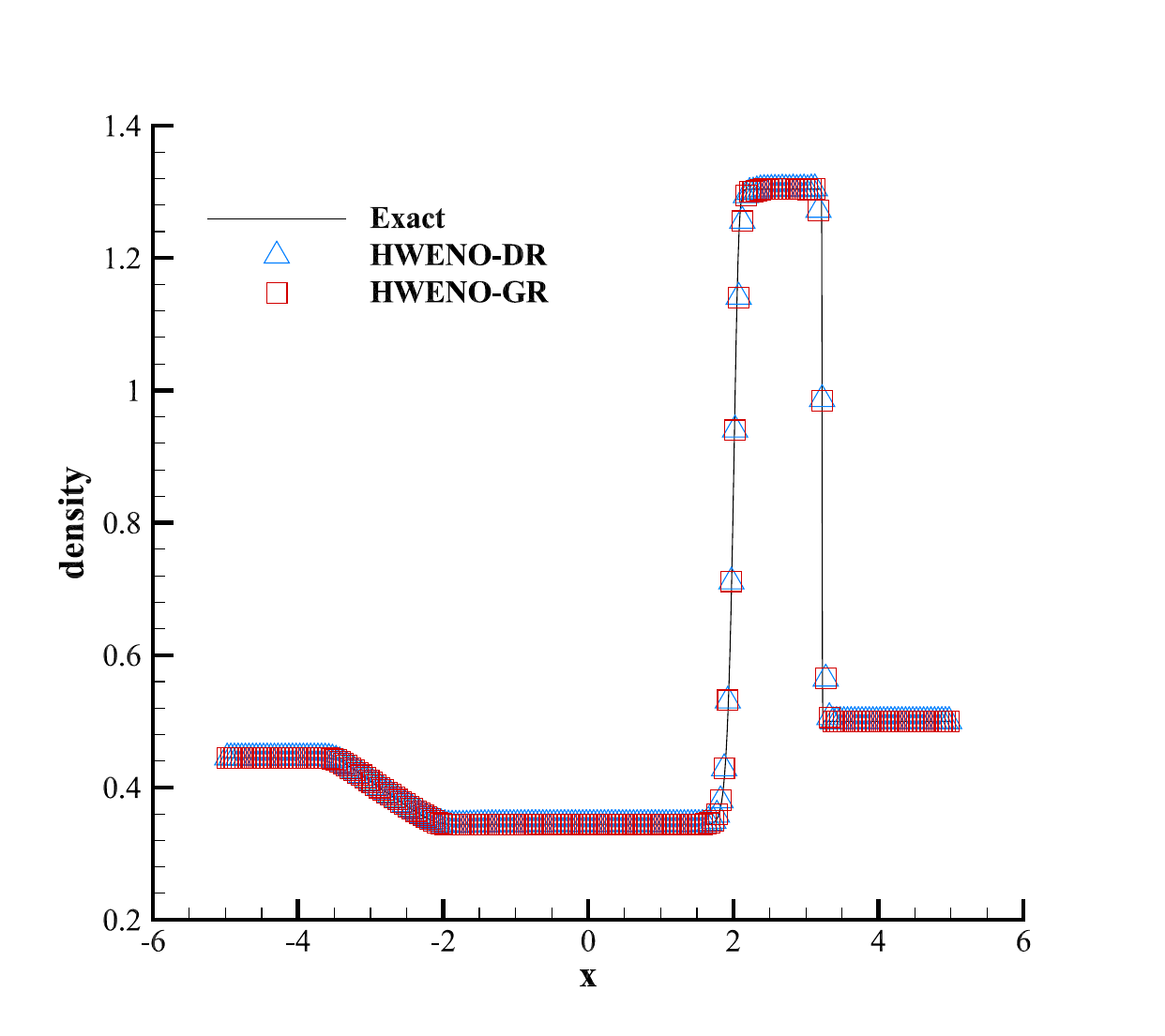}}
	\subfigure[$Re = 10000$]{\includegraphics[width=7.5cm,angle=0]{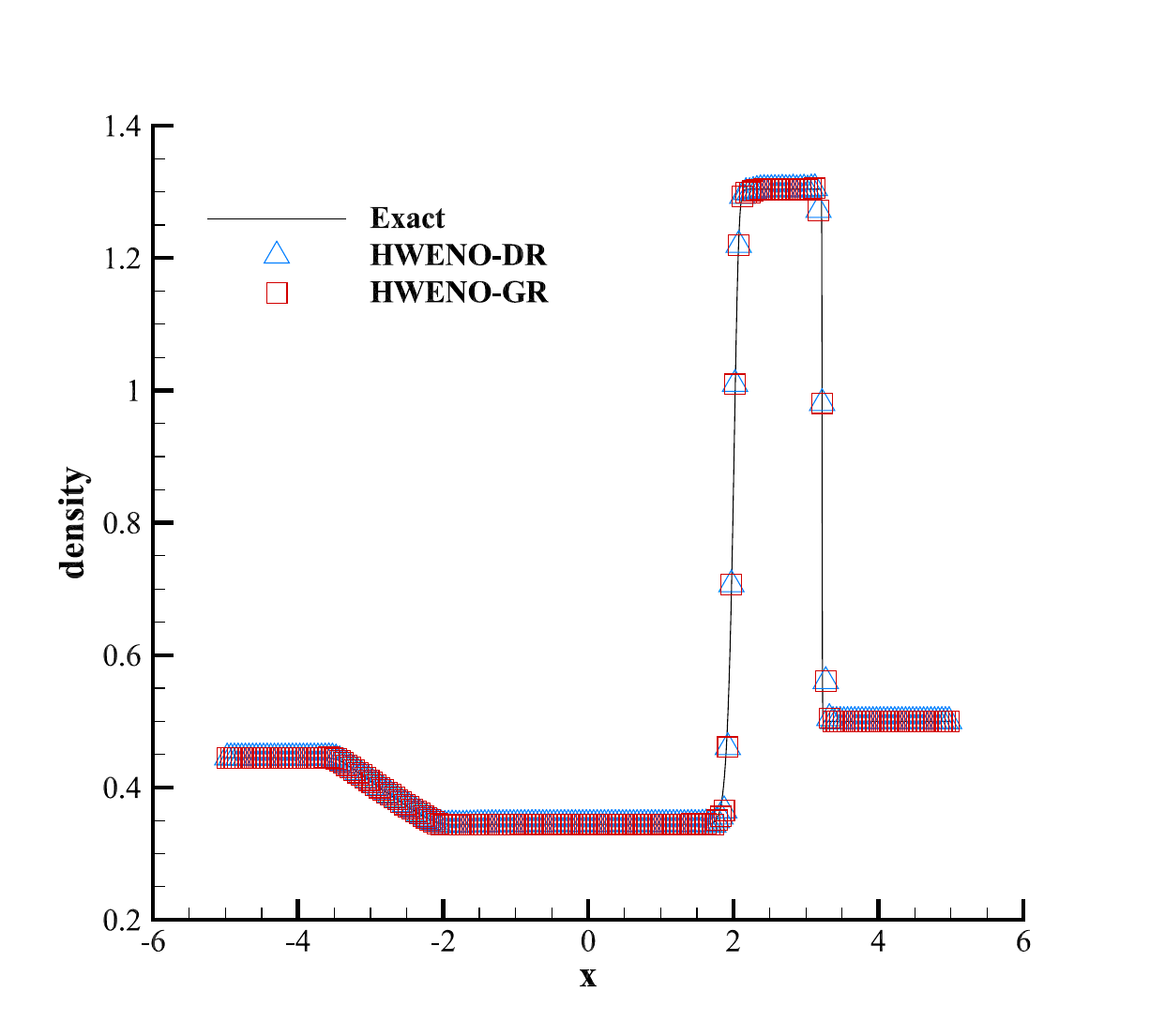}}
	\subfigure[$Re = 100000$]{\includegraphics[width=7.5cm,angle=0]{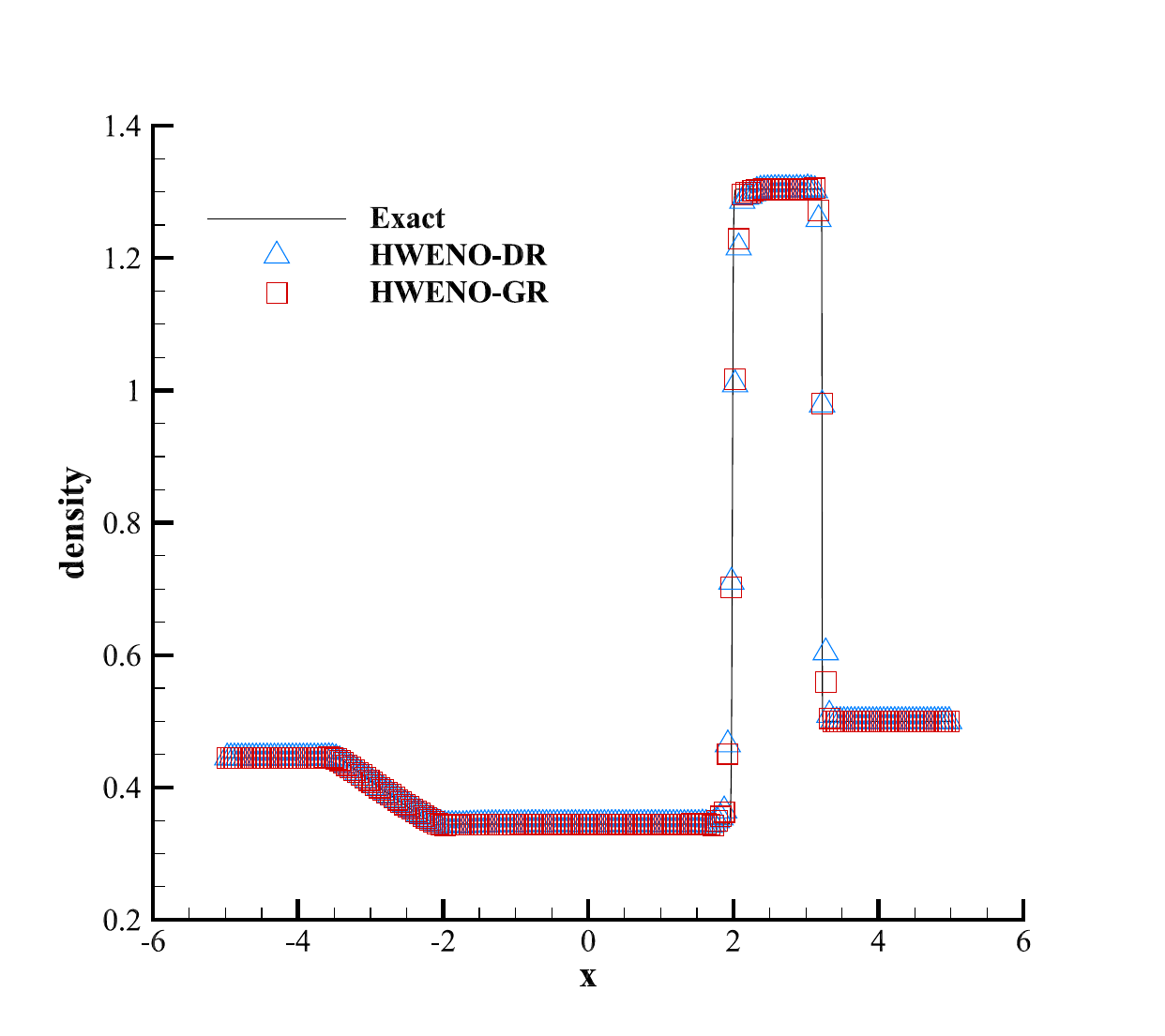}}
	\caption{Example \ref{Example:Lax1D}.
		Solutions computed by the HWENO-GR scheme and the HWENO-DR scheme on a uniform mesh with 200 cells.}
	\label{sec3:Lax1D_WENOHWENO}
\end{figure}
\end{example}

\begin{example}\label{Example:DoubleRare}
		We solve the double rarefaction wave problem for one-dimensional Navier--Stokes equations with \(Re=1000\) and the initial condition
		\begin{equation*}
			(\rho,\mu,p,\gamma)=\begin{cases}
				(7,-1,0.2,1.4),~-1<x<0,
				\\(7,1,0.2,1.4),~~~~~~~0<x<1.
			\end{cases}
		\end{equation*}	
		This test contains low-pressure and low-density regions, where many high-order schemes may produce negative density or pressure and fail. The final time is $T = 0.6$. The left and right boundary conditions are inflow and outflow, respectively. Figure~\ref{sec3:DoubleRare} gives the density, velocity, and pressure profiles. Similar profiles are obtained with both schemes, with good agreement with the exact solution.

		\begin{figure}[!htpb]
			\centering
			\subfigure{
			{\includegraphics[width=5.35cm,angle=0]{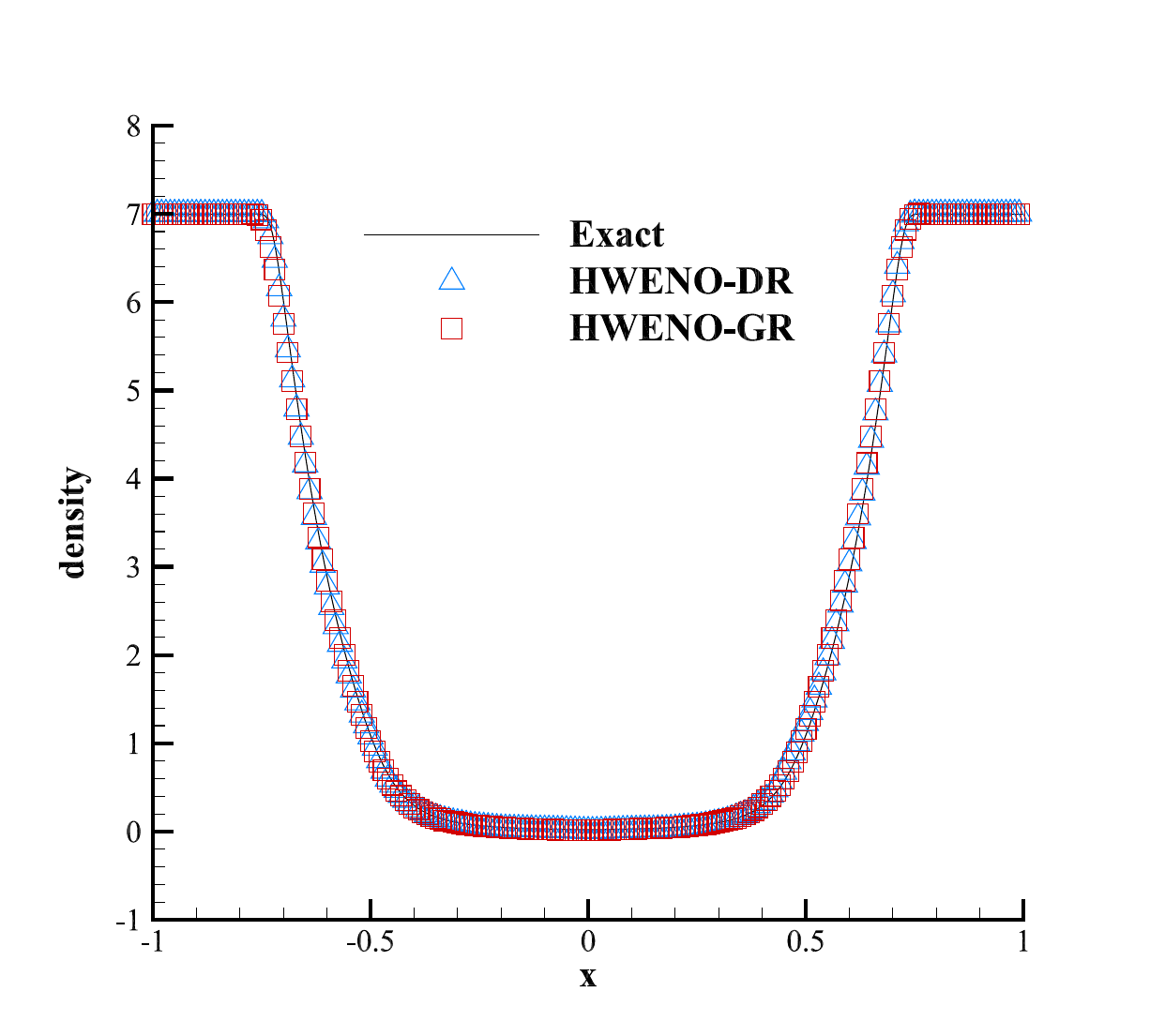}}
			{\includegraphics[width=5.35cm,angle=0]{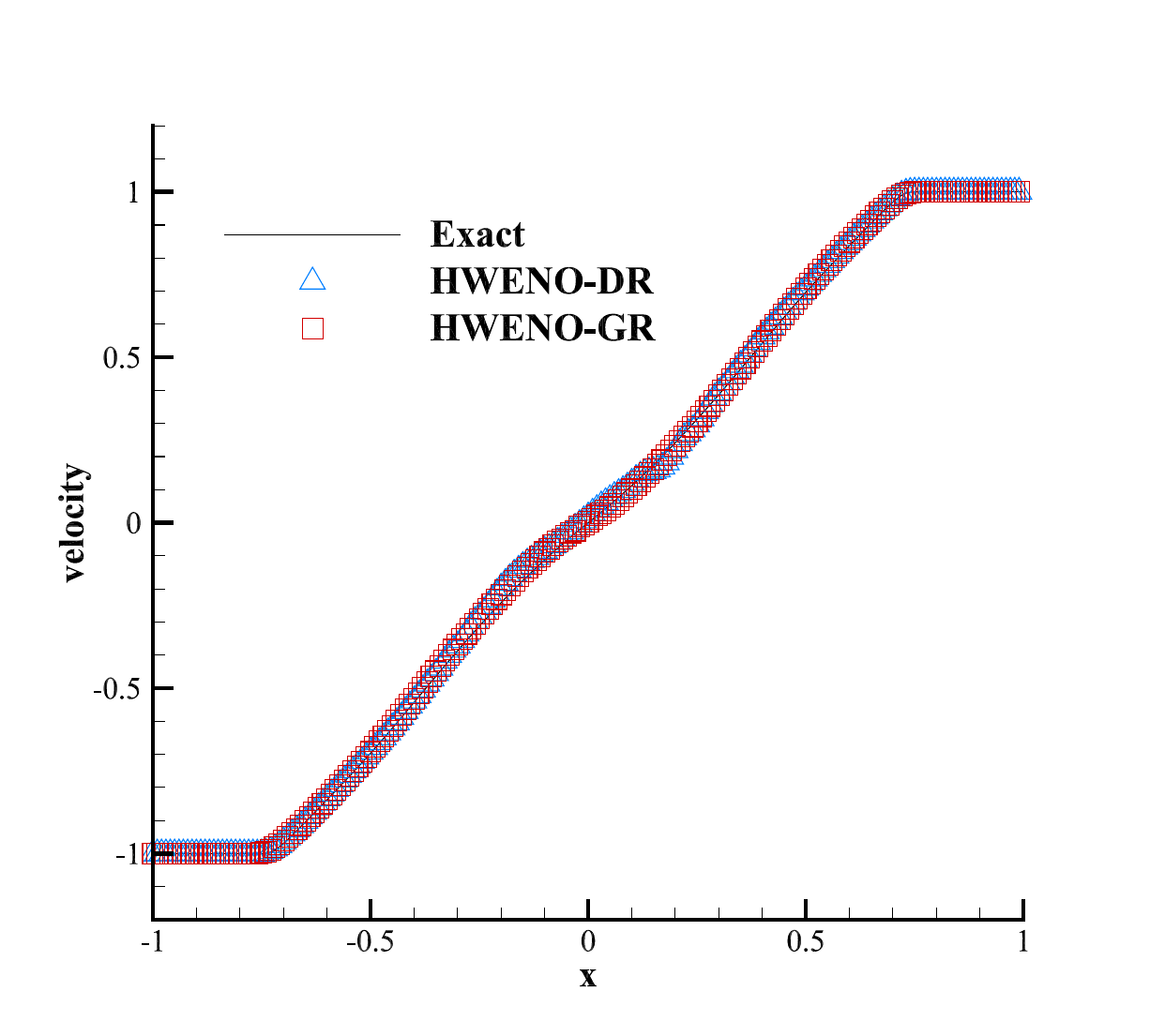}}
			{\includegraphics[width=5.35cm,angle=0]{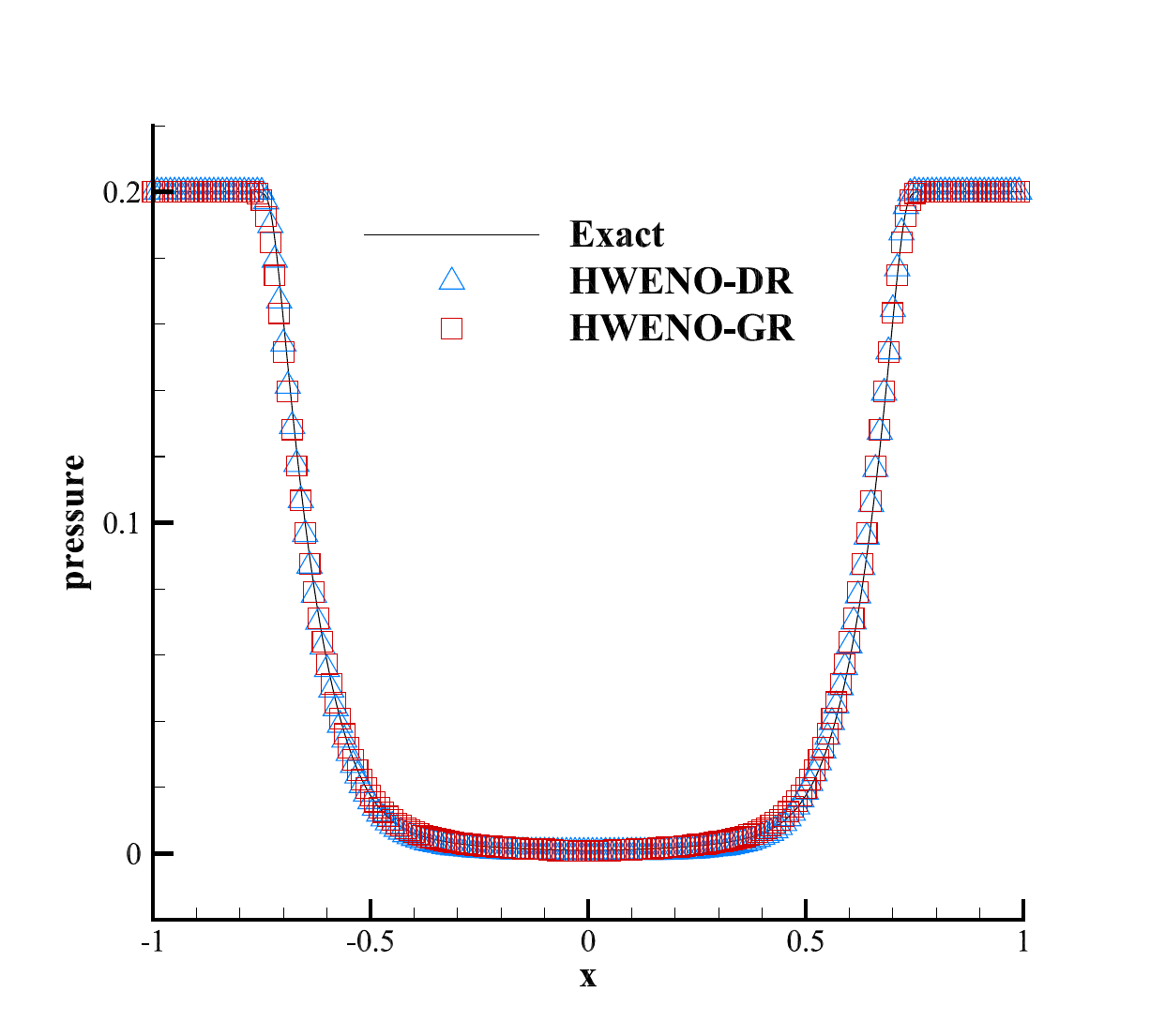}}}
			\caption{Example \ref{Example:DoubleRare}. Double rarefaction wave problem with \(Re=1000\) and $200$ cells.}
			\label{sec3:DoubleRare}
		\end{figure}
\end{example}
\begin{example}\label{Example:Sedov1D}
		We solve the Sedov blast wave problem for one-dimensional Navier--Stokes equations with the initial condition
		\begin{equation*}
			(\rho,\mu,E,\gamma)=\begin{cases}
				(1,0,10^{-12},1.4),~~x\in [-2,2] \setminus  \mbox{the center cell},
				\\(1,0,\frac{3200000}{\dx},1.4),~x\in \mbox{the center cell}.
			\end{cases}
		\end{equation*}	
		The final time is $T = 0.001$. The inlet and outlet conditions are imposed on the left and right boundaries, respectively. Figure~\ref{sec3:Sedov1D} shows the computed density, velocity, and pressure. The two sets of numerical profiles are similar and remain in good agreement with the exact solution.

		\begin{figure}[t]
			\centering
			\subfigure{
			{\includegraphics[width=5.35cm,angle=0]{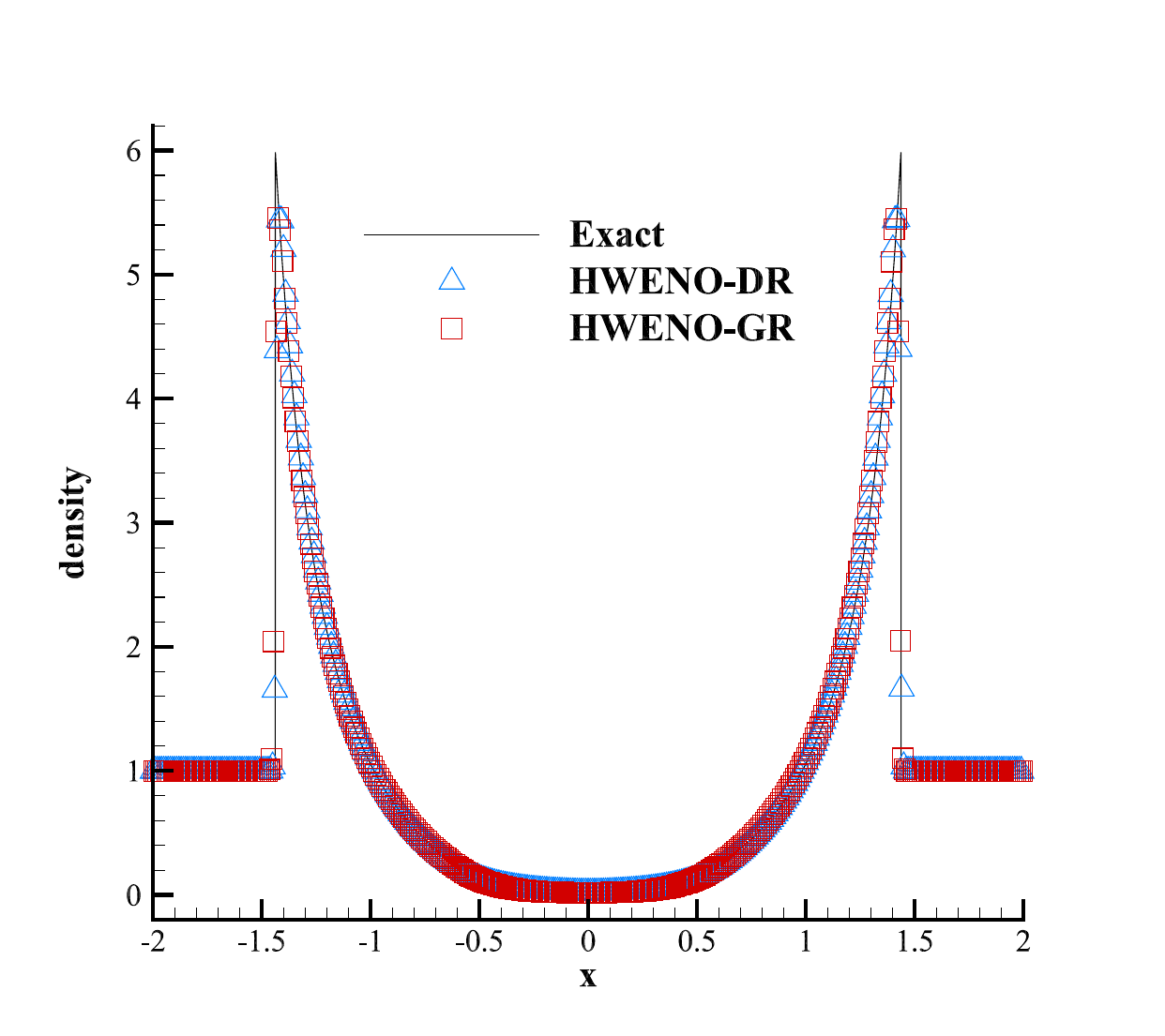}}
			{\includegraphics[width=5.35cm,angle=0]{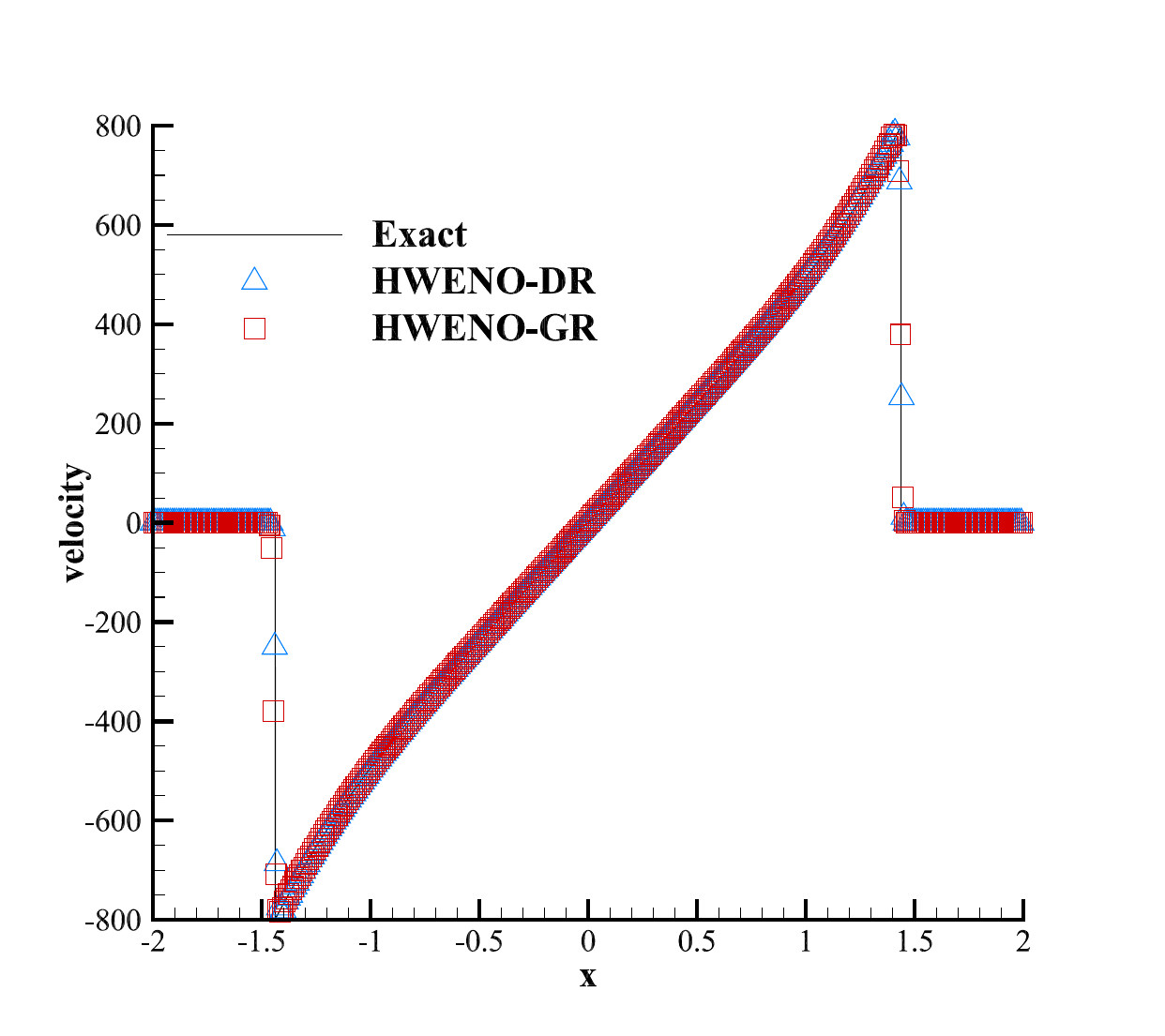}}
			{\includegraphics[width=5.35cm,angle=0]{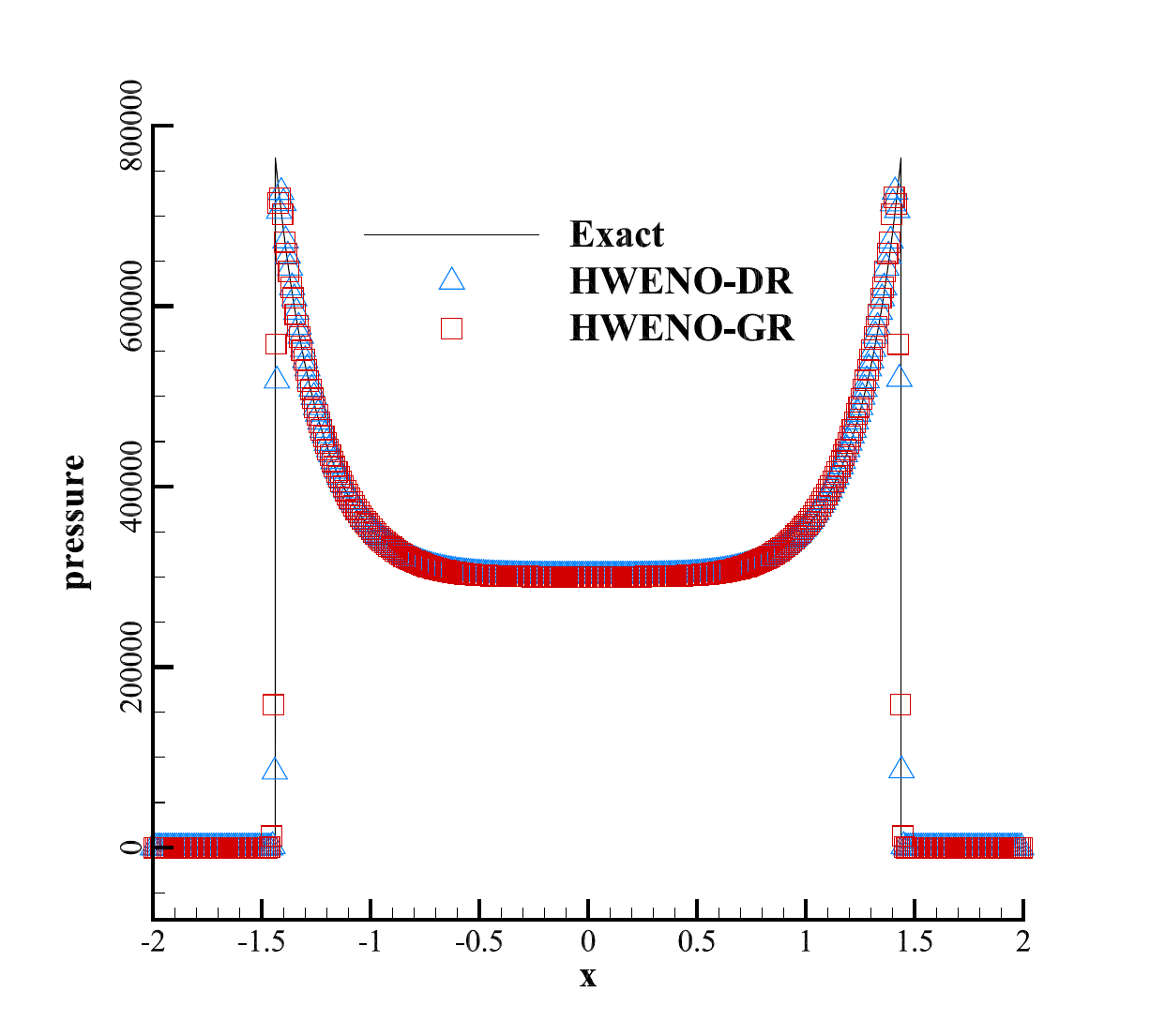}}}
			\caption{Example \ref{Example:Sedov1D}. One-dimensional Sedov problem with \(Re=1000\) and $400$ cells.}
			\label{sec3:Sedov1D}
		\end{figure}
\end{example}
\begin{example}\label{Example:Leblanc1D}
	We solve the Leblanc shock-tube problem for one-dimensional Navier--Stokes equations with the initial condition
	\begin{equation*}
		(\rho,\mu,p,\gamma)=\begin{cases}
			(2,0,10^{9},1.4),~~x\in [-10,0),
			\\(0.001,0,1,1.4),~x\in [0,10].
		\end{cases}
	\end{equation*}
	This problem contains extremely large jumps in density and pressure. The final time is $T = 0.0001$. The inlet and outlet conditions are imposed on the left and right boundaries, respectively. Figure~\ref{sec3:Leblanc1D} presents the density, velocity, and pressure profiles. Both schemes give similar results that agree well with the exact solution.

	\begin{figure}[t]
		\centering
		\subfigure{
		{\includegraphics[width=5.35cm,angle=0]{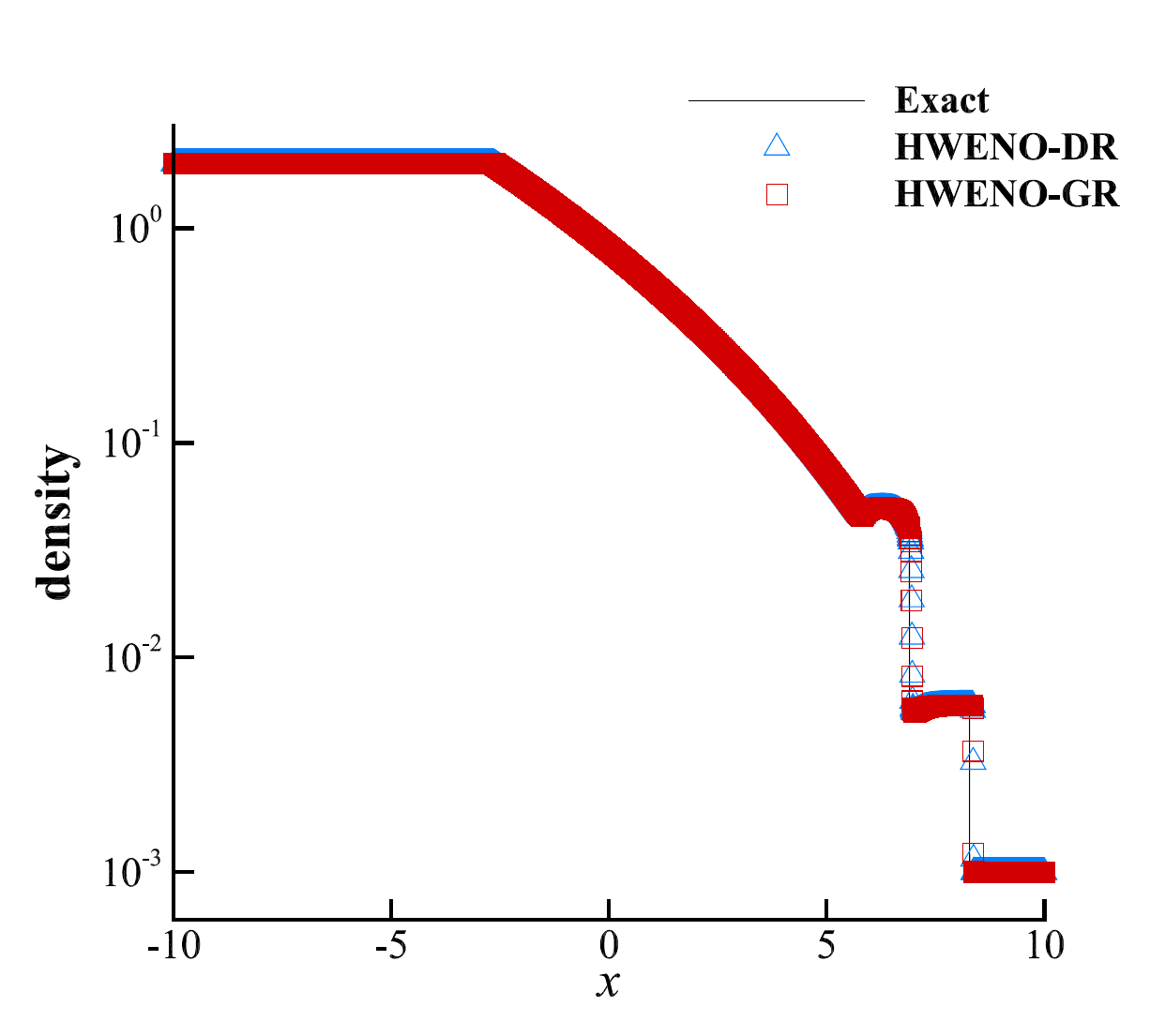}}
		{\includegraphics[width=5.35cm,angle=0]{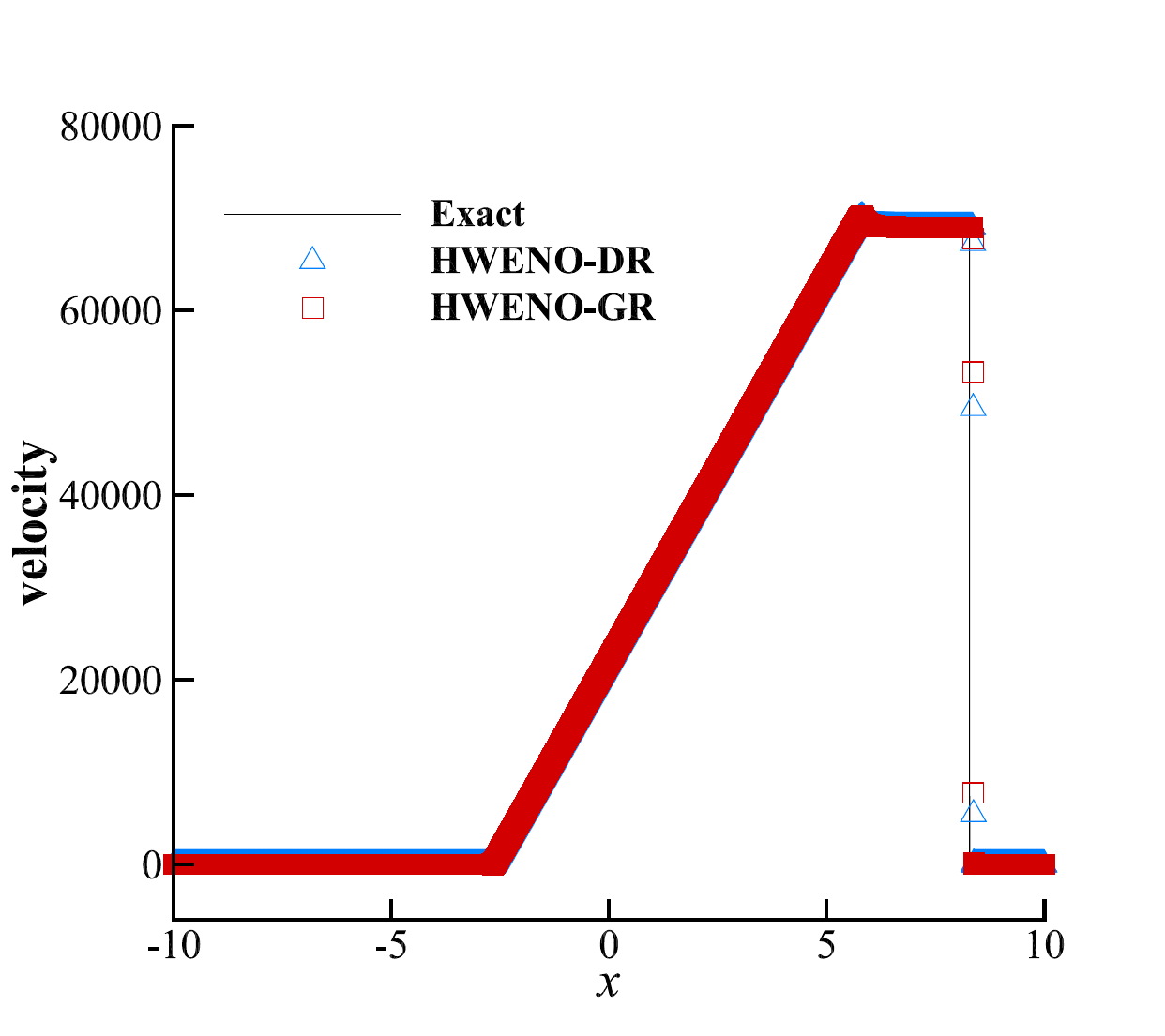}}
		{\includegraphics[width=5.35cm,angle=0]{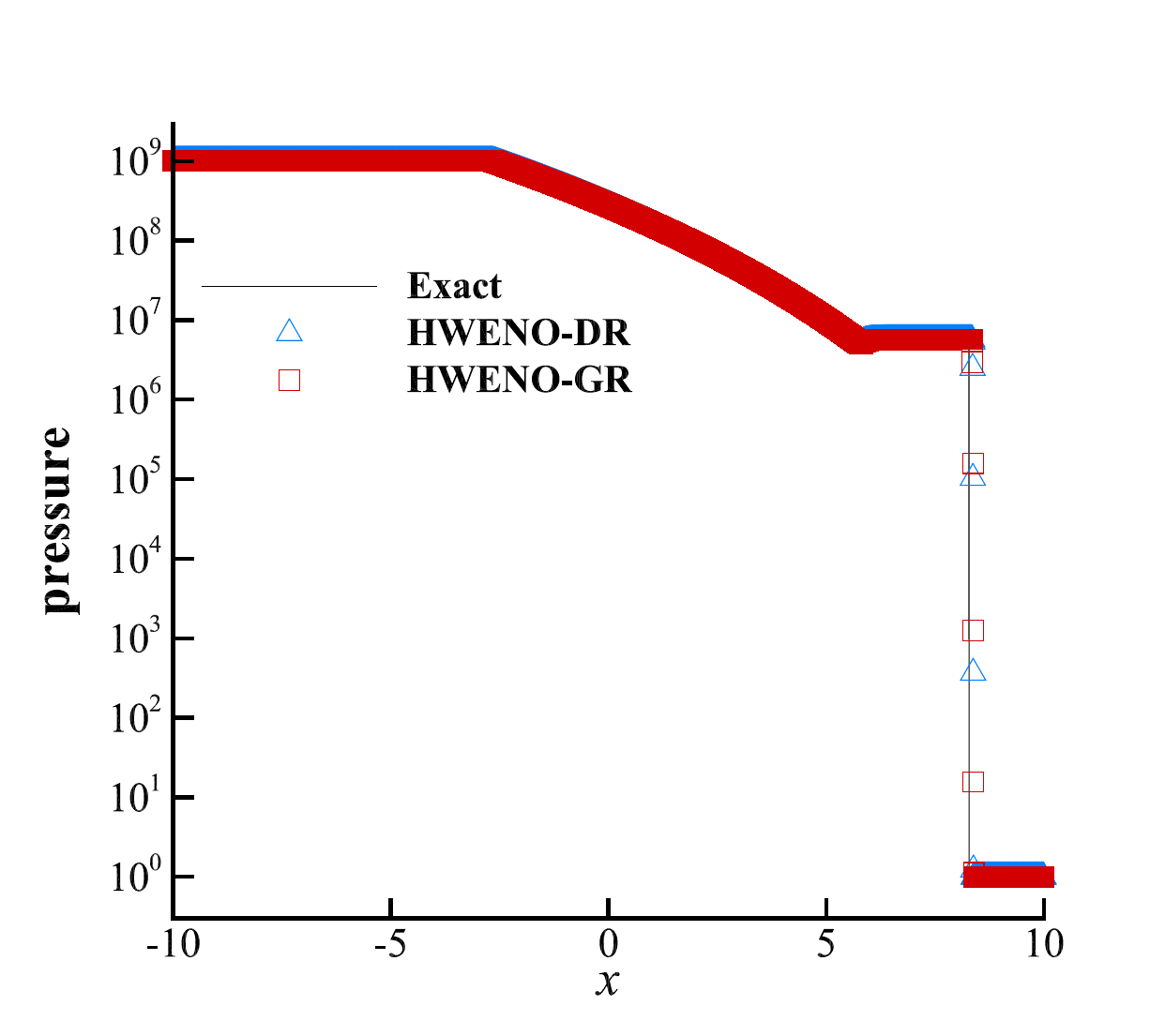}}}
		\caption{Example \ref{Example:Leblanc1D}. One-dimensional Leblanc problem with \(Re=1000\) and $3200$ cells.}
		\label{sec3:Leblanc1D}
	\end{figure}
\end{example}
\begin{example}\label{Example:Sedov2D}
		We solve a Sedov blast wave problem \cite{sedov1959similarity} for two-dimensional Navier--Stokes equations.	 The computational domain is $[0,1.1]\times[0,1.1]$ and the initial condition is
		\begin{equation*}
			(\rho,\mu,\nv,E,\gamma)=\begin{cases}
				(1,0,0,\frac{0.244816}{\dx\dy},1.4),~(x,y)\in[0,\dx]\times[0,\dy],
				\\(1,0,0,10^{-12},1.4),~~~~~\mbox{otherwise}.
			\end{cases}
		\end{equation*}	
		This test contains a strong shock and a large density variation. The PP-limiter is required; without it, negative density or pressure may occur. Reflective boundary conditions are employed on the left and bottom, while outflow conditions are applied on the right and upper boundaries. Figure~\ref{sec3:Sedov2D} presents the density contours at the final time $T=1$ for the HWENO-GR scheme and the HWENO-DR scheme with \(Re=1000\). The density contours obtained with the two schemes are similar.

		\begin{figure}[!htpb]
			\centering
			\subfigure[HWENO-DR scheme]{\includegraphics[width=7.5cm,height=7.5cm,angle=0]{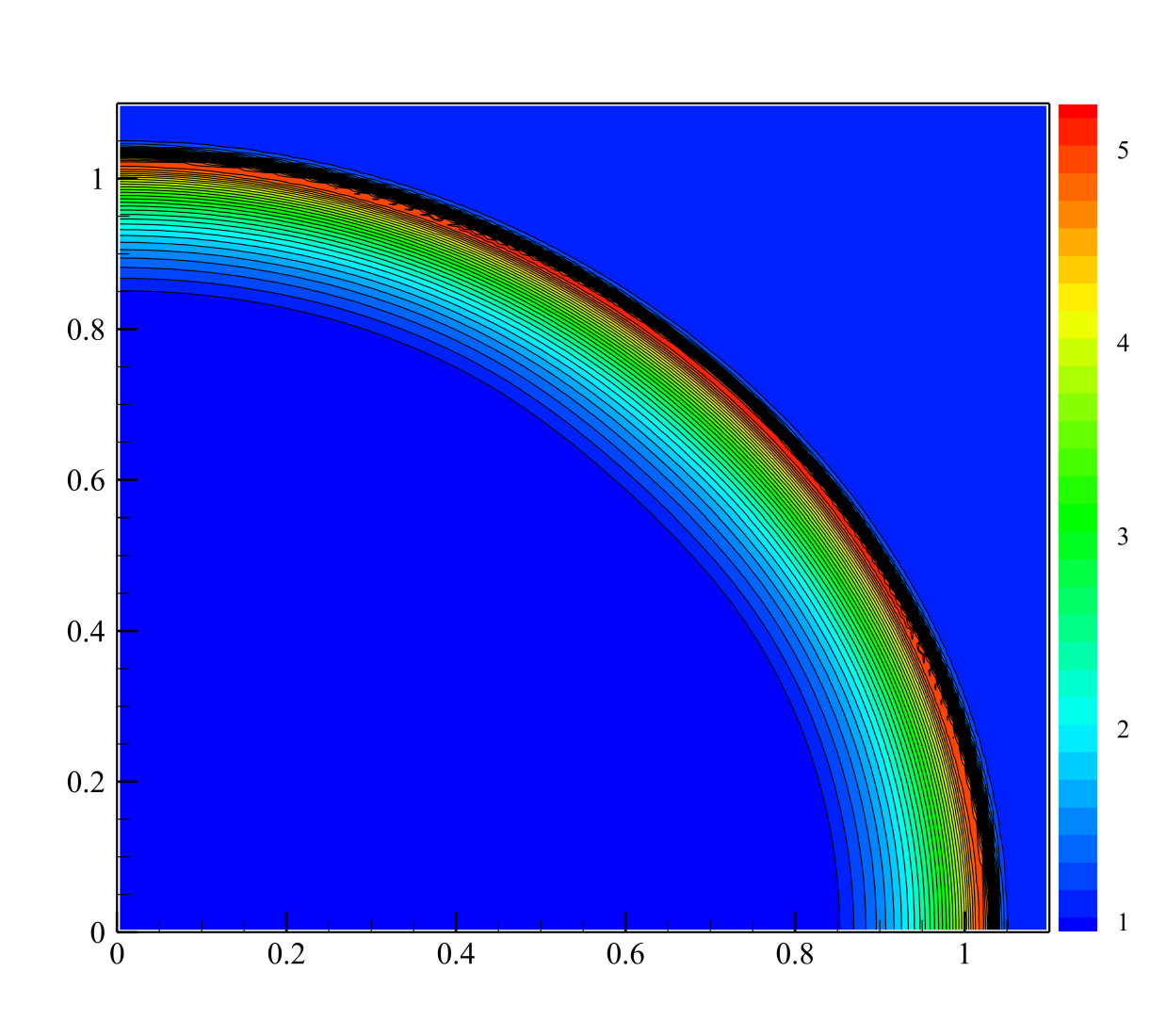}}
			\subfigure[HWENO-GR scheme]{\includegraphics[width=7.5cm,height=7.5cm,angle=0]{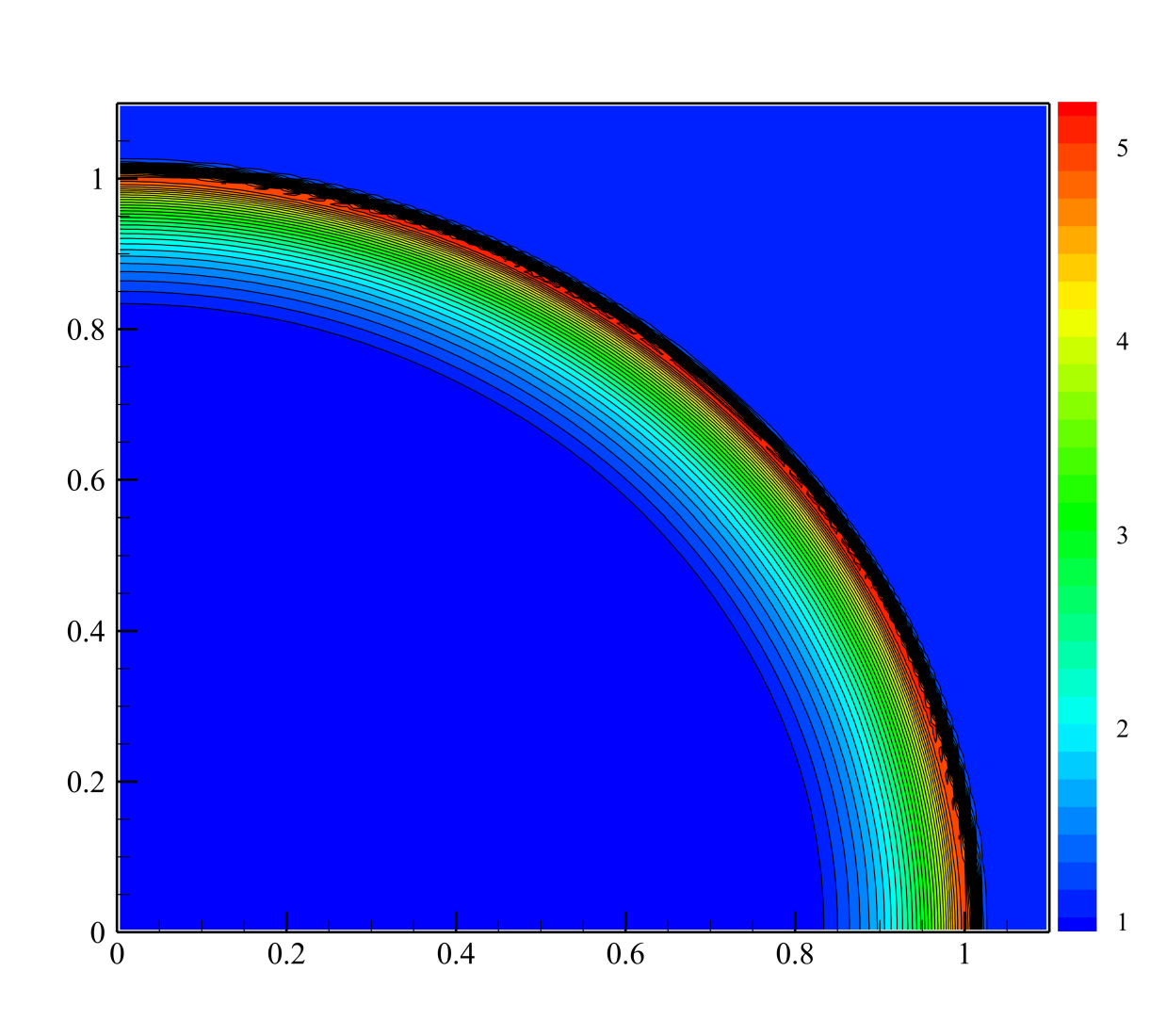}}
			\caption{Example \ref{Example:Sedov2D}.
				Two-dimensional Sedov problem with \(Re=1000\).  Contour plots of density with 30 equally spaced lines from 0.95 to 5.	Uniform meshes: $160\times160$.}
			\label{sec3:Sedov2D}
		\end{figure}
		To examine the reconstruction choice discussed in Remark~\ref{rem:nonlinear-gradient-reconstruction}, we monitor the PP wave-speed bounds during the viscous computation with \(Re=1000\). At reconstruction step \(n\), the recorded quantity is
		\[
		\beta_{\max}^{(n)}=max_{i,j}\left\{
		\beta_{i+\frac12}^{(n)},\,\beta_{j+\frac12}^{(n)}
		\right\},
		\]
		where \(\beta_{i+\frac12}\) and \(\beta_{j+\frac12}\) are the directional bounds in \eqref{sec2:lf-flux-2d}. Figure~\ref{sec3:Sedov2DBeta} plots \(\beta_{\max}^{(n)}\) on a logarithmic scale. With linear reconstruction of the weak-derivative moments, the bound remains between approximately \(10^{13}\) and \(10^{16}\) over the sampled steps. These large values lead to repeated inadmissible stages and step rejection, reducing the accepted time step to about \(10^{-10}\) and preventing the computation from reaching \(T=1\) in practice. With nonlinear HWENO reconstruction, the bound drops rapidly after the first few reconstruction steps and subsequently lies mainly between \(10^2\) and \(10^8\).
		\begin{figure}[!htpb]
			\centering
			\subfigure[Linear gradient reconstruction]{\betamaxplot{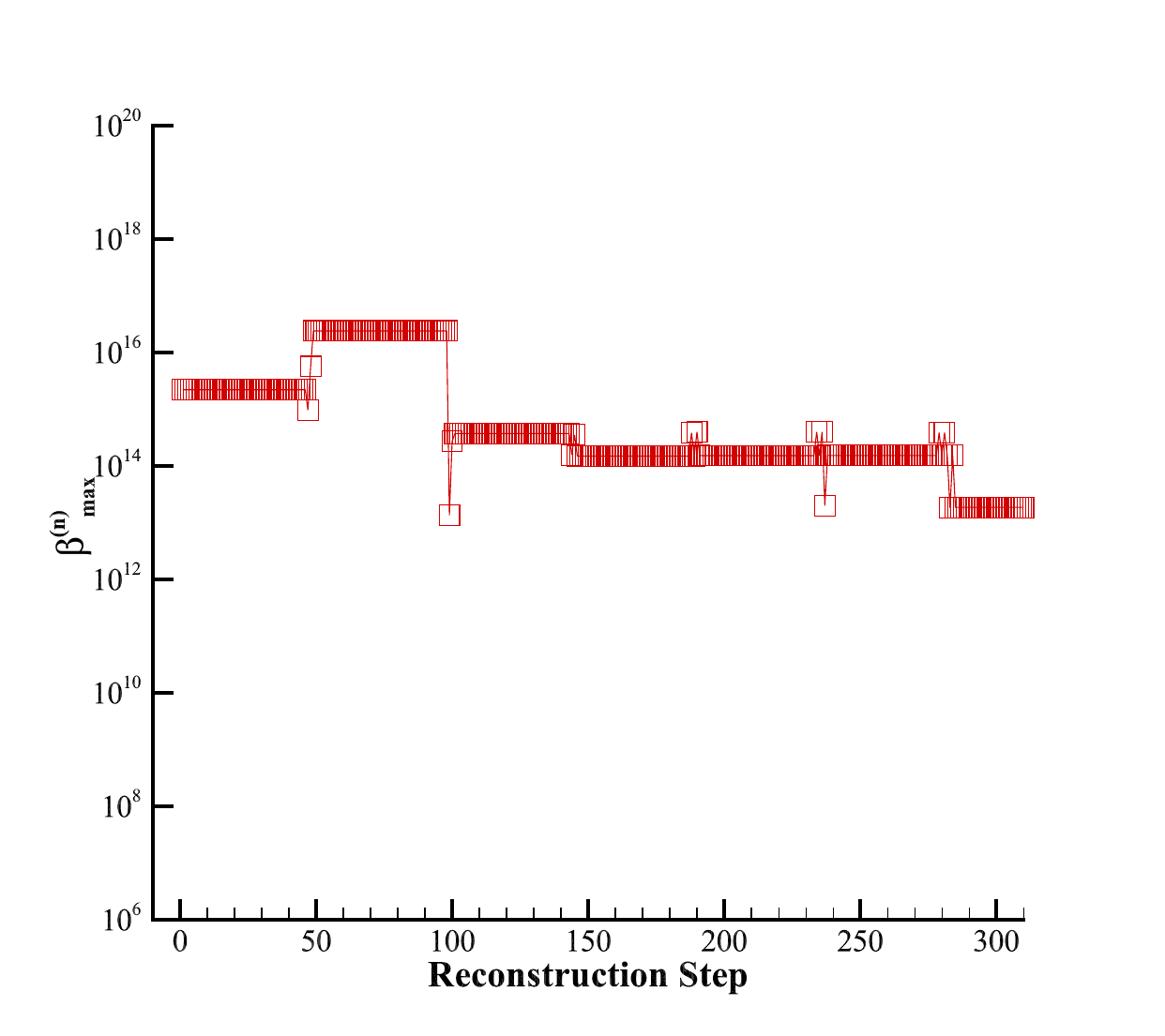}}
\subfigure[Nonlinear HWENO gradient reconstruction]{\betamaxplot{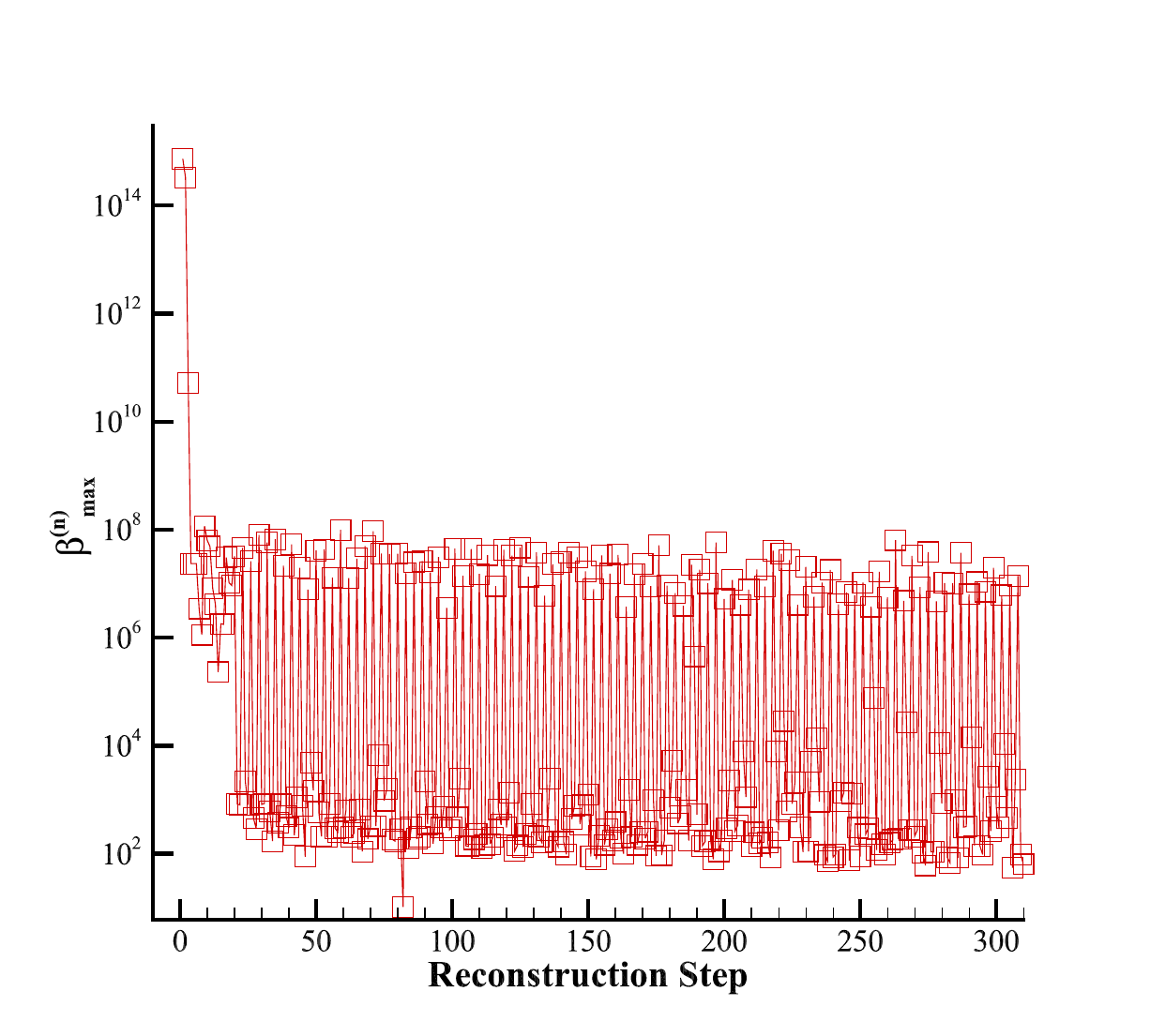}}
			\caption{Example \ref{Example:Sedov2D}. Evolution of the global maximum PP wave-speed bound \(\beta_{\max}^{(n)}\) for the two-dimensional viscous Sedov problem with \(Re=1000\). {Here, \(\beta_{\max}^{(n)}\) denotes the maximum value of \(\beta\) at reconstruction step \(n\).}}
			\label{sec3:Sedov2DBeta}
		\end{figure}
\end{example}
\begin{example}\label{Example:ShockDiffraction}
	We solve the shock-diffraction problem for the two-dimensional Navier--Stokes equations with \(Re=1000\). The computational domain is the union of $[0,1]\times[6,11]$ and $[1,13]\times[0,11]$. The initial condition is a right-moving Mach-5.09 shock, initially located at $x=0.5$ and $6\le y \le 11$, moving into undisturbed air with density 1.4 and pressure 1. The boundary conditions are inflow at $x=0$, $6 \le y\le 11$, outflow at $x=13$, $0\le y \le 11$, $1\le x\le 13$, $y=0$ and $0\le x\le 13$, $y=11$, and reflective at the walls $0\le x\le 1$, $y=6$ and $x=1$, $0\le y \le 6$. The final time is $T=2.3$. The density contours computed by the HWENO-GR scheme and the HWENO-DR scheme are shown in Fig. \ref{sec3:ShockDiffraction}. Similar density patterns are obtained with the two schemes.

	\begin{figure}[!htpb]
		\centering
		\subfigure[HWENO-DR scheme] {\includegraphics[width=6.5cm,height=5.5cm,angle=0]{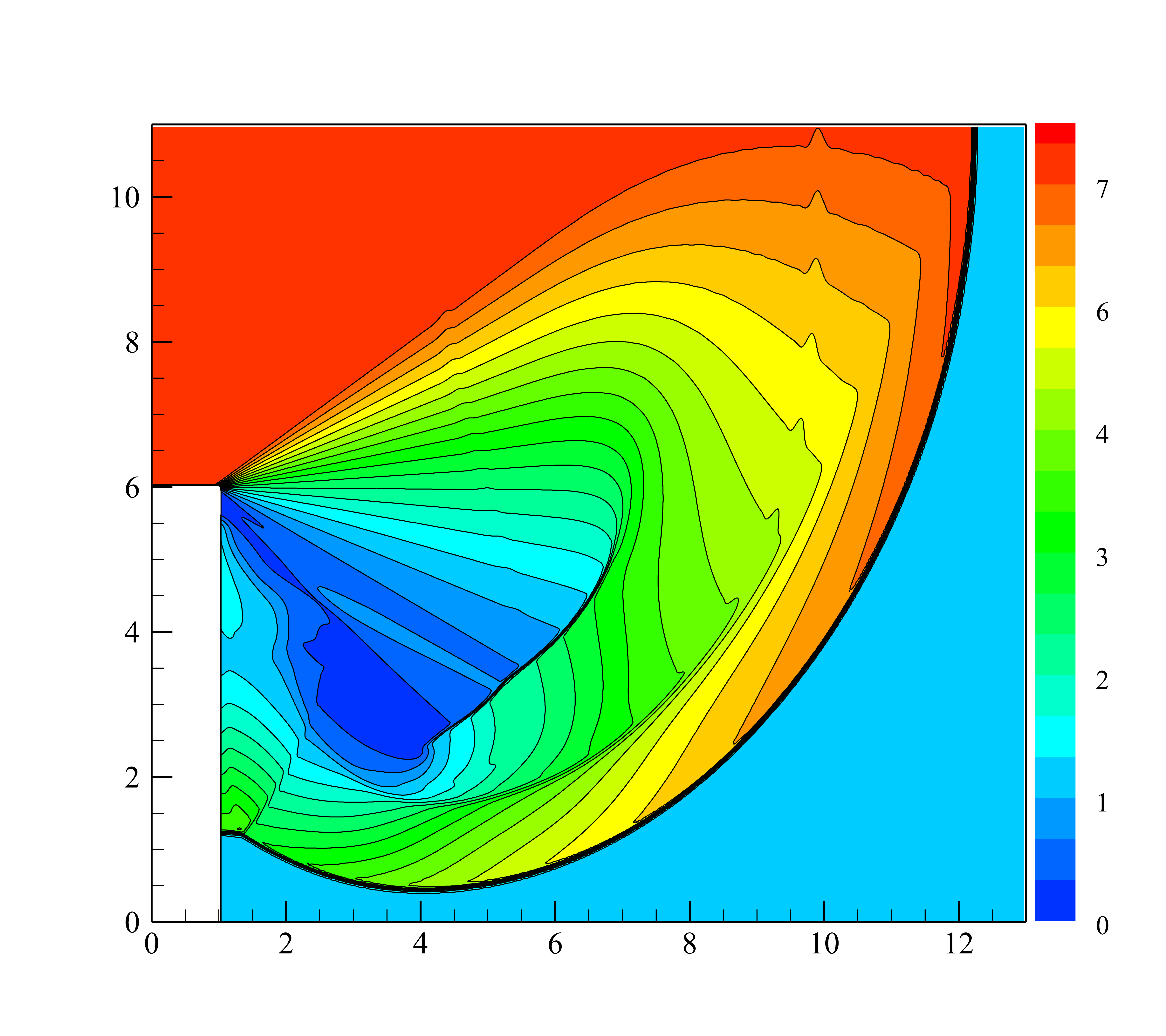}}
		\subfigure[HWENO-GR scheme] {\includegraphics[width=6.5cm,height=5.5cm,angle=0]{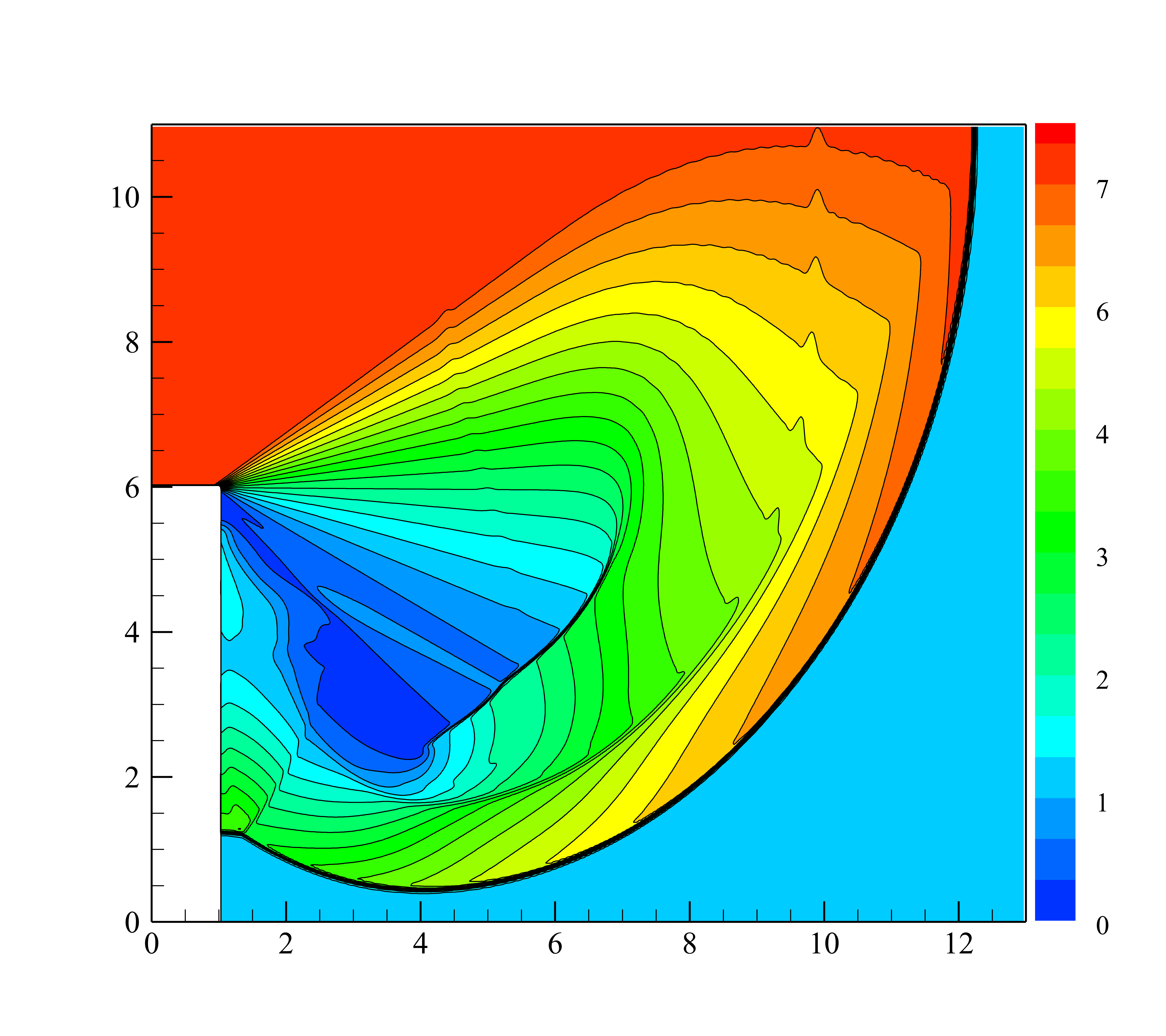}}
		\caption{Example \ref{Example:ShockDiffraction}. Shock-diffraction problem with \(Re=1000\). Contour plots of density with 20 equally spaced lines from 0.066227 to 7.0668. Uniform meshes: $64\times64$.}
		\label{sec3:ShockDiffraction}
	\end{figure}
\end{example}
\begin{example}\label{Example:HM2000}
	\begin{figure}[!htpb]
		\centering
		\subfigure[HWENO-DR scheme] {\includegraphics[width=8cm,height=4cm,angle=0]{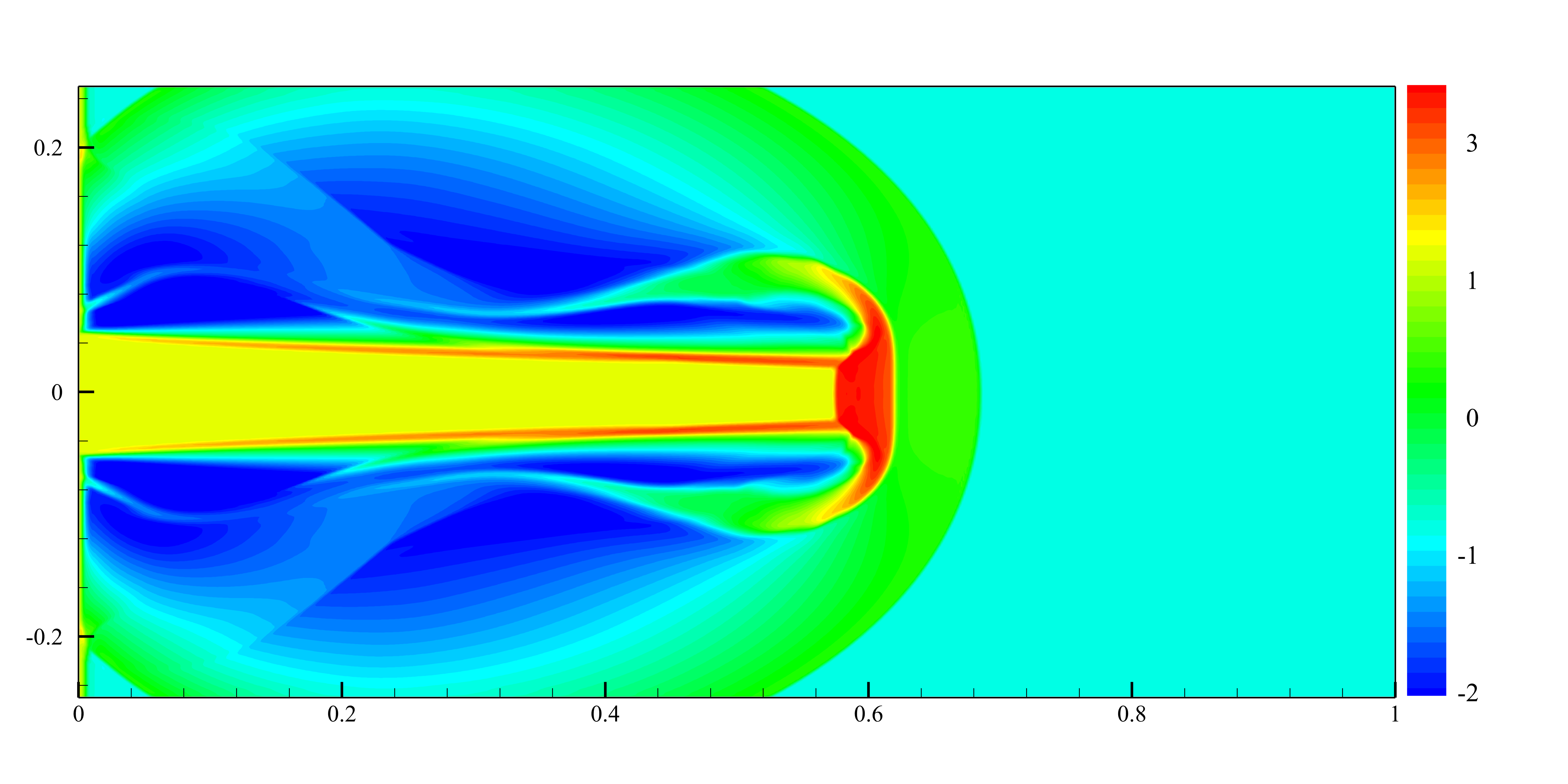}}
		\subfigure[HWENO-GR scheme] {\includegraphics[width=8cm,height=4cm,angle=0]{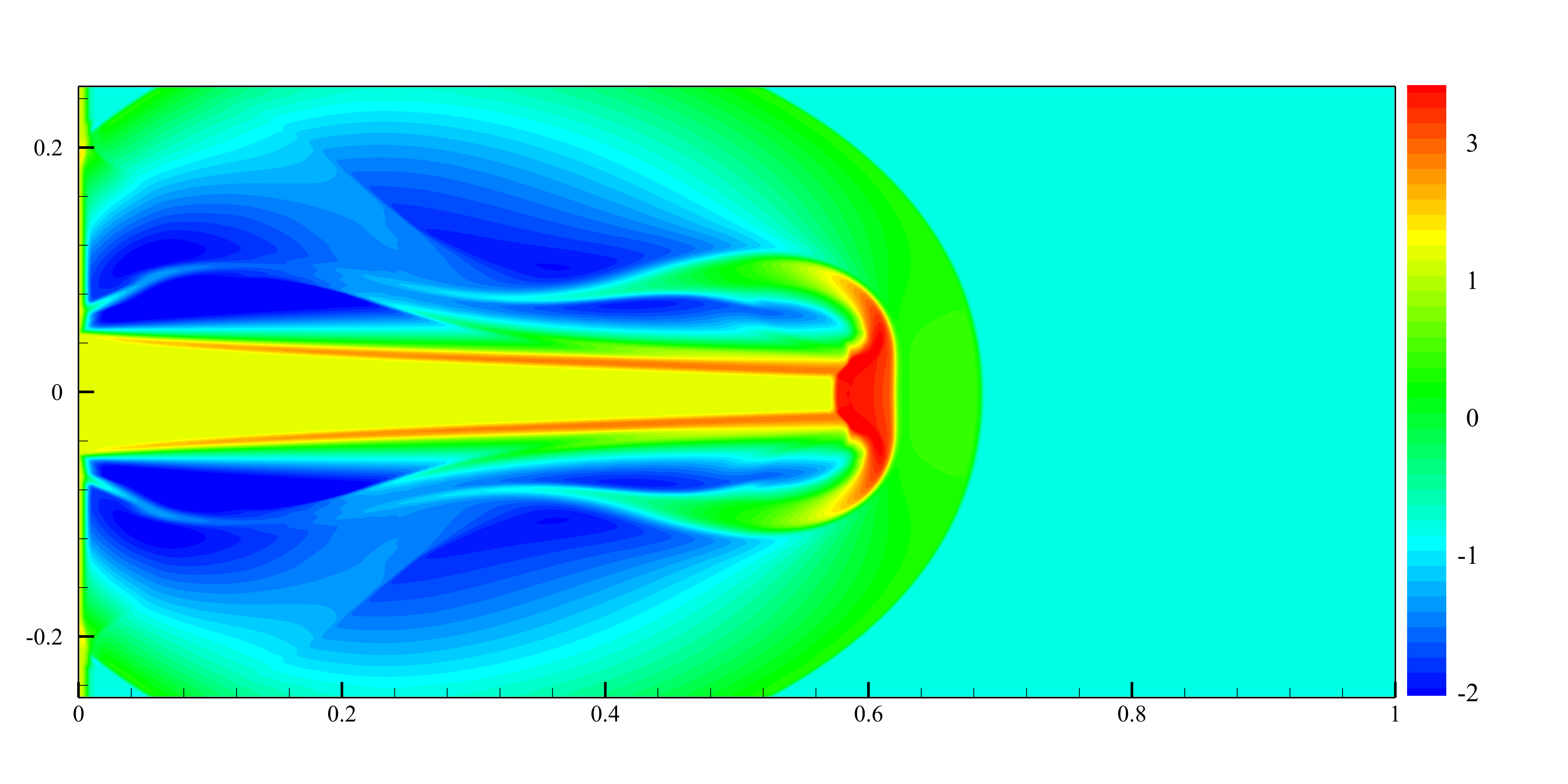}}
		\caption{Example \ref{Example:HM2000}. Mach-2000 problem with \(Re=1000\). Contour plots of density with 40 equally spaced lines from -2 to 3 and scales are logarithmic. Uniform meshes: $640\times320$.}
		\label{sec3:HM2000}
	\end{figure}
		Motivated by high-Mach astrophysical jet simulations \cite{ha2005numerical}, we consider the Mach-2000 benchmark without radiative cooling \cite{ZhangS,ZhangS2012} at \(Re=1000\). The computational domain is $[0,1]\times[-0.25,0.25]$. The domain is initially filled with ambient gas with $(\rho,\mu,\nv,p,\gamma)=(0.5,0,0,0.4127,\frac53)$. Outflow boundary conditions are imposed on the right, top, and bottom. On the left boundary, $(\rho, \mu, \nv, p, \gamma) = (5, 800, 0, 0.4127, \frac{5}{3})$ is prescribed for $|y|<0.05$, and $(5, 0, 0, 0.4127, \frac{5}{3})$ is prescribed elsewhere.
		Figure \ref{sec3:HM2000} shows the results obtained by the HWENO-GR scheme and the HWENO-DR scheme at the final time $T=0.001$. The two density fields are similar and consistent with those reported in \cite{ZhangS,ZhangS2012}.
\end{example}
\begin{example}\label{Example:Mach10shock}
	We solve the Mach-10 shock problem with shock reflection and diffraction for the two-dimensional Navier--Stokes equations with \(Re=1000\). The computational domain is the union of $[0,1]\times[0,1]$ and $[-1,1]\times[1,3]$. The initial condition is
	\begin{equation*}
		(\rho,\mu,\nv,p,\gamma)=\begin{cases}
			(8,\frac{33}{4}\sin(\frac{\pi}{3}),-\frac{33}{4}\cos(\frac{\pi}{3}),116.5,1.4),~x<\frac{1}{6}+\frac{y}{\sqrt{3}},
			\\(1.4,0,0,1,1.4),\quad\quad\quad\quad\quad\quad\quad\quad\quad\mbox{otherwise}.
		\end{cases}
	\end{equation*}
	The boundary conditions are inflow on the left, outflow on the right and bottom, and reflective on the wall $[\frac{1}{6},1]\times\{0\}$ and $\{1\}\times [-1,0]$. The exact post-shock condition is imposed on $[0,\frac{1}{6}]\times \{0\}$. The top boundary follows the exact motion of a Mach-10 shock, with $(\rho,\mu,\nv,p)=(8, \frac{33}{4}\sin(\frac{\pi}{3}), -\frac{33}{4}\cos(\frac{\pi}{3}), 116.5)$ for $0\le x \le \frac{1}{6}+\frac{1+20t}{\sqrt{3}}$ and $(\rho,\mu,\nv,p)=(1.4,0,0,1)$ elsewhere. The final time is $T = 0.2$. The density contours computed by the HWENO-GR scheme and the HWENO-DR scheme are shown in Fig. \ref{sec3:Mach10shock}. The two schemes produce similar density contours.

	\begin{figure}[!htpb]
		\centering
		\subfigure[HWENO-DR scheme] {\includegraphics[width=8.1cm,height=5.4cm,angle=0]{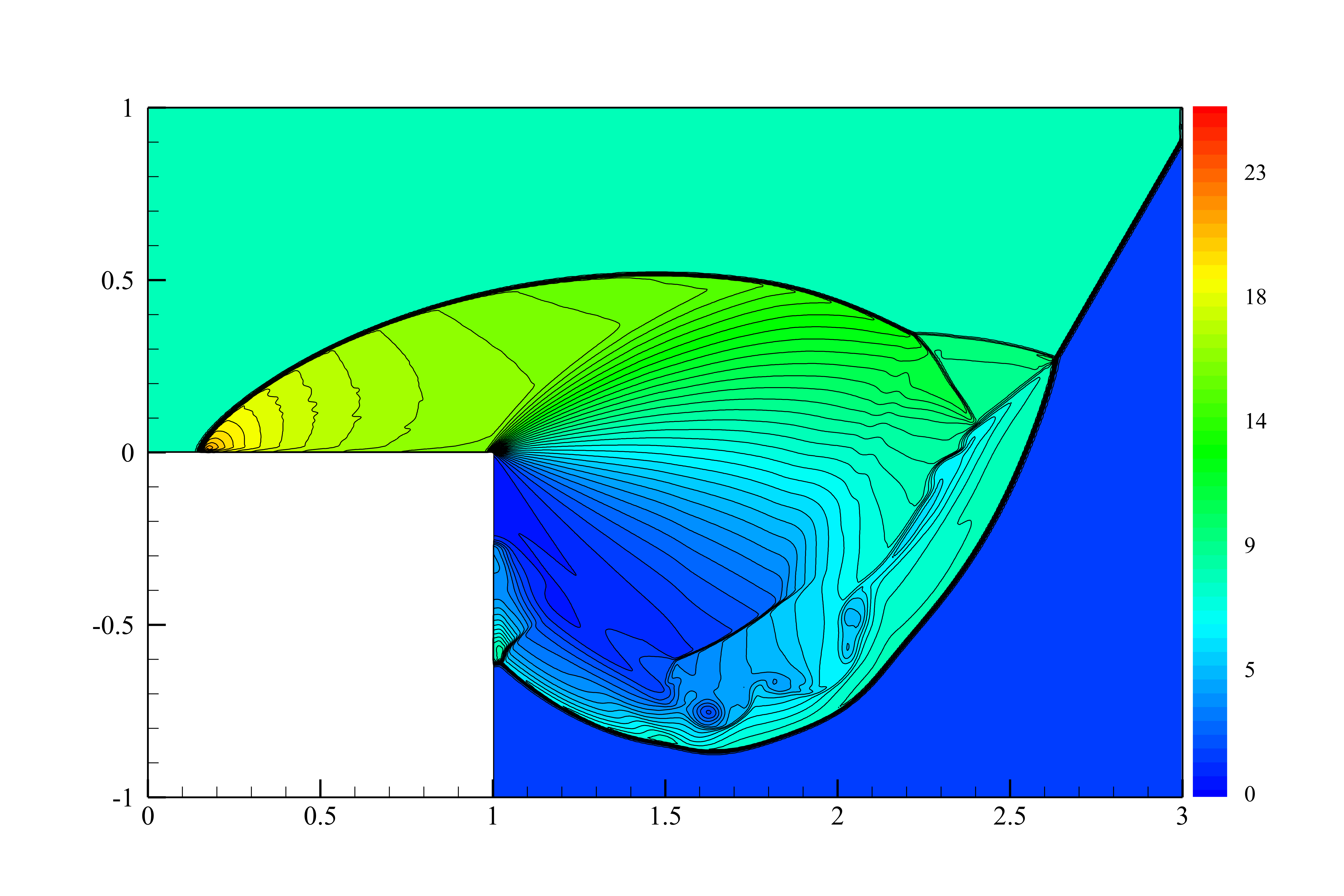}}
		\subfigure[HWENO-GR scheme] {\includegraphics[width=8.1cm,height=5.4cm,angle=0]{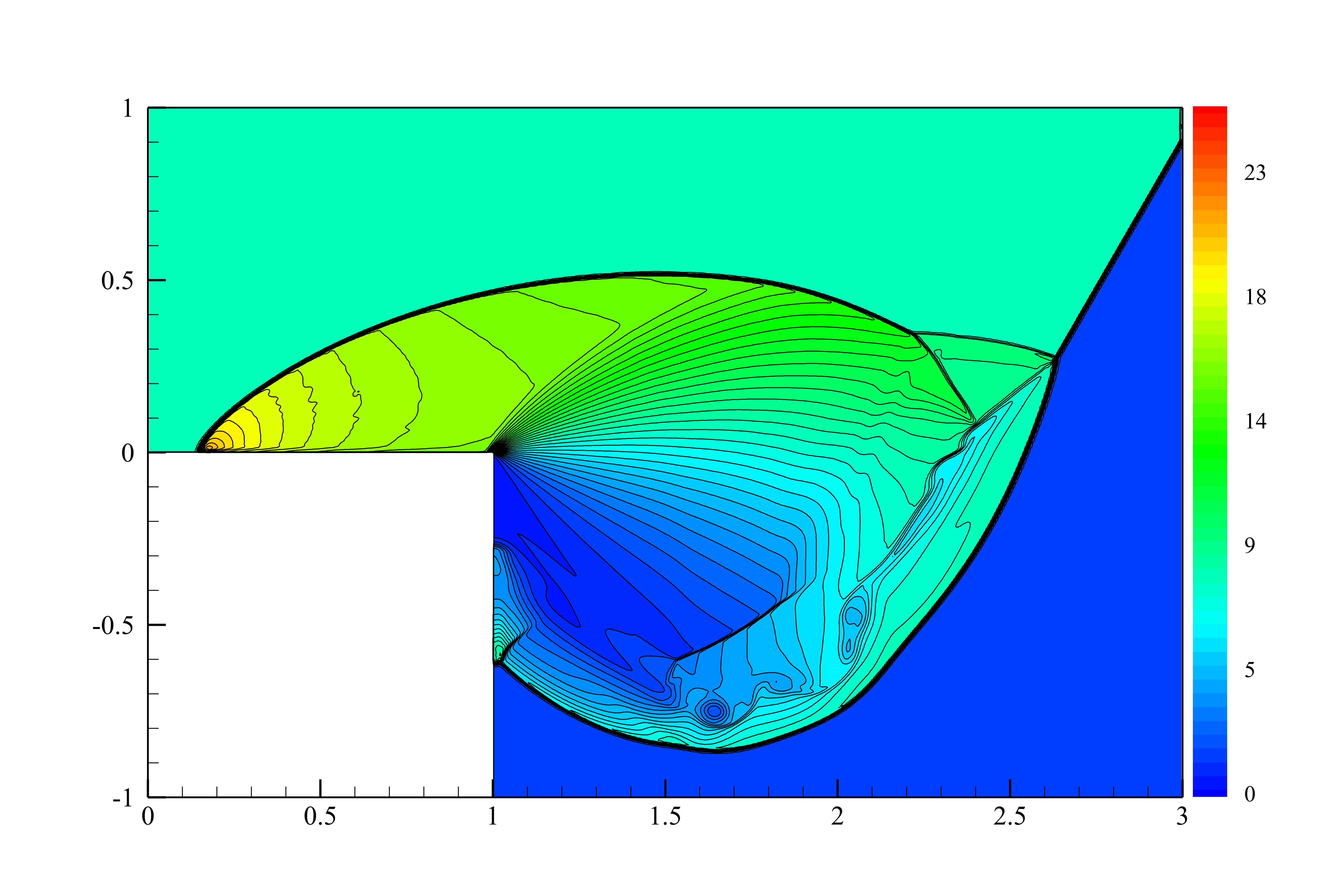}}
		\caption{Example \ref{Example:Mach10shock}.
		Mach-10 shock problem with \(Re=1000\). Contour plots of density with 50 equally spaced lines from 0 to 25. Uniform meshes: $480\times480$.}
		\label{sec3:Mach10shock}
	\end{figure}
\end{example}
\begin{example}\label{Example:ViscousShockMixing}
	We consider a shock wave impinging on a spatially developing mixing layer for the two-dimensional compressible Navier--Stokes equations \cite{PengZhaiNiYongShen2019,SandhamYee1998}. This test examines the resolution of vortical structures generated by the interaction between the shock and the shear layer. The computational domain is $[0,200]\times[-20,20]$. At the inflow, the streamwise velocity and the transverse perturbation are prescribed as
	\begin{equation*}
		\mu(y)=2.5+0.5\tanh(2y),\qquad
		\nv'(y,t)=\sum_{k=1}^{2}a_k\cos\left(2\pi k\frac{t}{T_0}+\phi_k\right)\exp\left(-\frac{y^2}{b}\right),
	\end{equation*}
	where $b=10$, $a_1=a_2=0.05$, $\phi_1=0$, $\phi_2=\pi/2$, and $T_0=\lambda/u_c$ with $\lambda=30$ and $u_c=2.68$. The density and pressure on the two sides of the mixing layer are
	\begin{equation*}
		(\rho,p)=
		\begin{cases}
			(1.6374,0.3327), & y>0,\\
			(0.3626,0.3327), & y<0.
		\end{cases}
	\end{equation*}
	The upper-boundary post-shock state is $(\rho,\mu,\nv,p)^T=(2.1101,2.9709,-0.1367,0.4754)^T$, and a slip-wall condition is imposed at the lower boundary. We take $\gamma=1.4$, $Re=500$, and $\mathrm{Pr}=0.72$. The computation is performed to $T=120$ on a $500\times100$ uniform grid. Figure \ref{sec3:ViscousShockMixing} shows the density contours obtained by the HWENO-GR scheme and the HWENO-DR scheme. The resulting density contours are similar for the two schemes and consistent with those reported in \cite{PengZhaiNiYongShen2019,SandhamYee1998}.

	\begin{figure}[!htpb]
		\centering
		\subfigure[HWENO-DR scheme]{\includegraphics[width=16cm,trim=200 15 200 15,clip,angle=0]{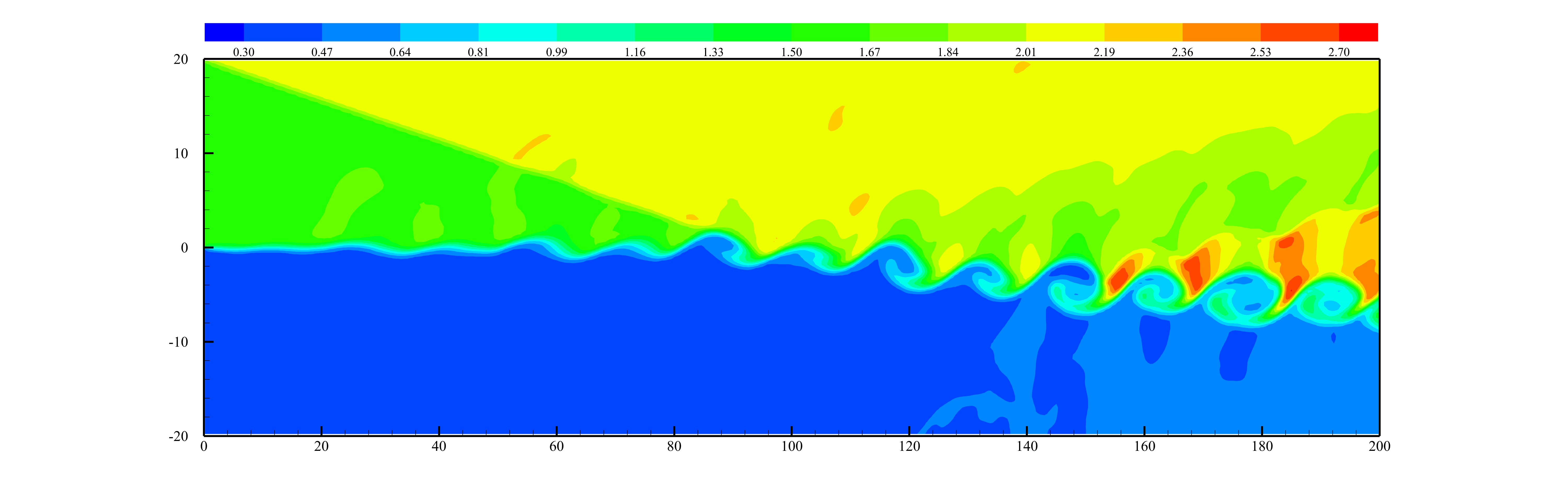}}
		\subfigure[HWENO-GR scheme]{\includegraphics[width=16cm,trim=200 15 200 15,clip,angle=0]{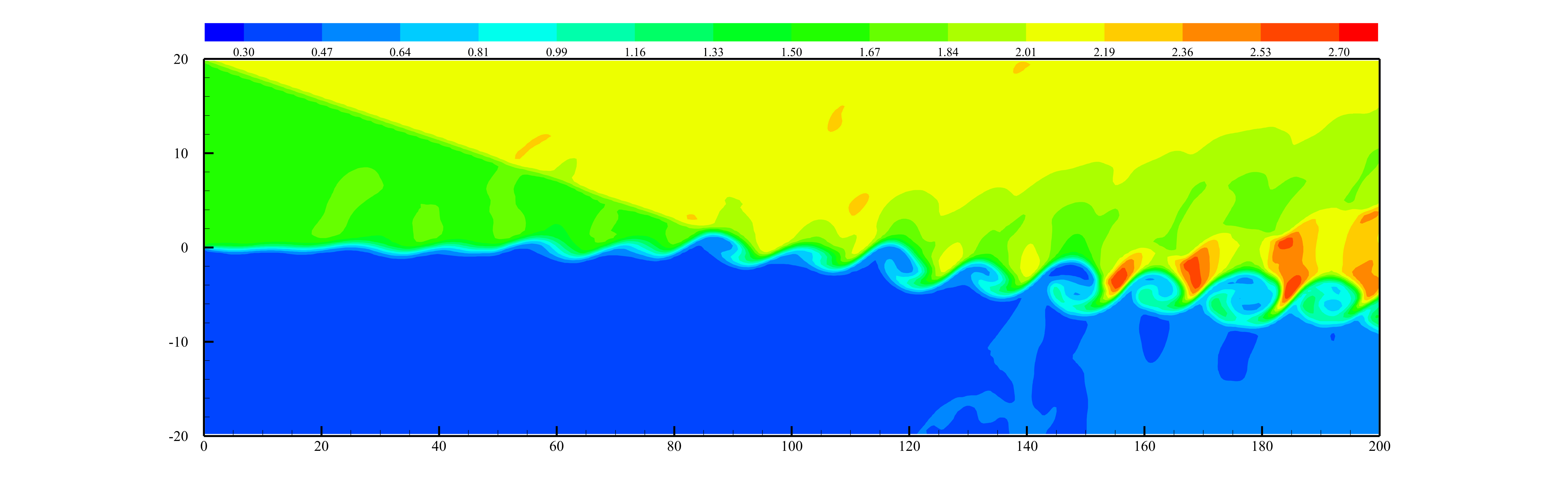}}
		\caption{Example \ref{Example:ViscousShockMixing}. Density contours of the two-dimensional viscous shock--mixing-layer problem at $T=120$. Uniform grid: $500\times100$. 15 equally spaced levels from 0.3 to 2.7.}
		\label{sec3:ViscousShockMixing}
	\end{figure}
\end{example}
\begin{example}\label{Example:ViscousShockTube}
	We consider the two-dimensional viscous shock-tube problem \cite{DaruTenaud2001,SjogreenYee2003,ZhouXuLiu2018}. The interaction of the initial shock, contact discontinuity, and viscous wall layer produces shock reflection, boundary-layer separation, and vortical structures. The computational domain is $[0,1]\times[0,0.5]$, and the diaphragm is initially located at $x=0.5$. With $\gamma=1.4$, the initial condition is
	\begin{equation*}
		(\rho,\mu,\nv,p)=
		\begin{cases}
			(120,0,0,120/\gamma), & x<0.5,\\
			(1.2,0,0,1.2/\gamma), & x\geq 0.5.
		\end{cases}
	\end{equation*}
	{The Prandtl and Reynolds numbers are \(\mathrm{Pr}=0.73\) and \(\mathrm{Re}=200\), respectively, and the dynamic viscosity is held constant.} An adiabatic no-slip condition is imposed on the solid walls, with a symmetry condition at $y=0.5$. The computation is advanced to $T=1$ on an $800\times400$ uniform grid using a CFL number of 0.3. Figure~\ref{sec3:ViscousShockTube} presents the density contours obtained by the HWENO-GR scheme and the HWENO-DR scheme. Both schemes give similar density fields, in agreement with the published results \cite{DaruTenaud2001,SjogreenYee2003,ZhouXuLiu2018,GuermondKronbichlerMaier2022}.

	\begin{figure}[!htpb]
		\centering
		\subfigure[HWENO-DR scheme]{\includegraphics[width=7.5cm,angle=0]{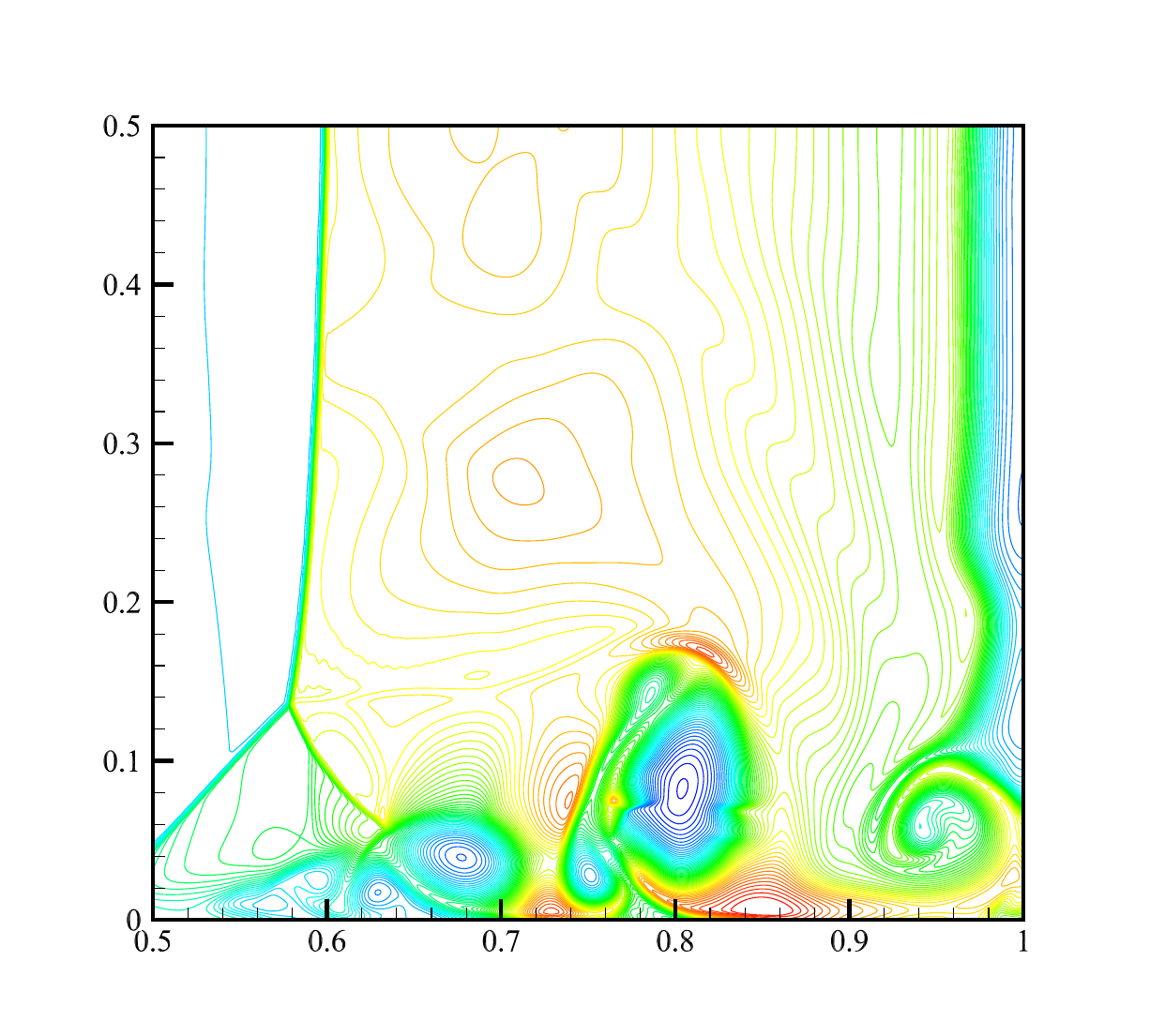}}
		\subfigure[HWENO-GR scheme]{\includegraphics[width=7.5cm,angle=0]{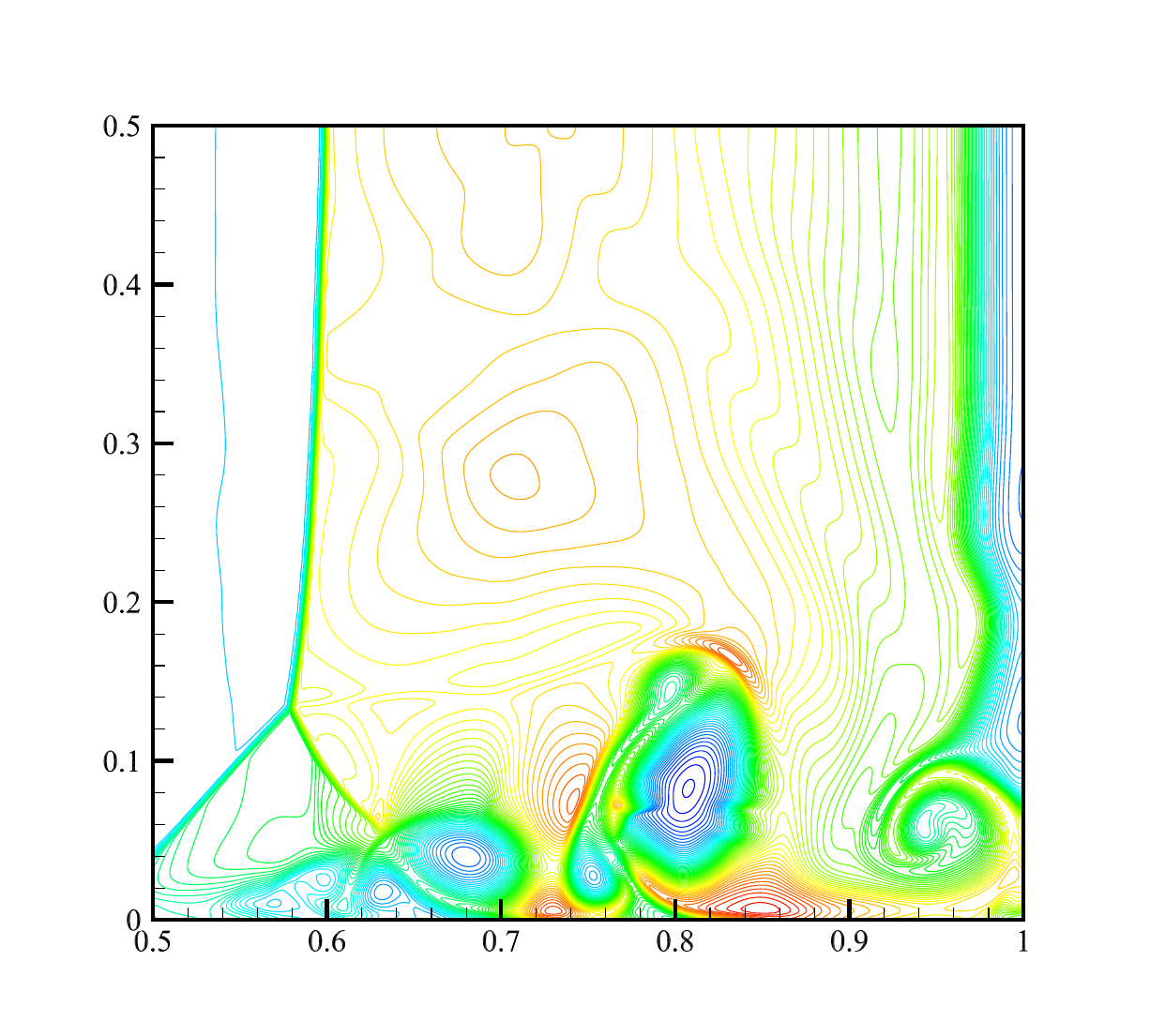}}
		\caption{Example \ref{Example:ViscousShockTube}. Density contours for the two-dimensional viscous shock-tube problem with $Re=200$ at $T=1.0$. Uniform grid: $800\times400$. 60 equally spaced levels from 20 to 120.}
		\label{sec3:ViscousShockTube}
	\end{figure}
\end{example}
\section{Concluding remarks}\label{Sec5:conclusion}
\label{Concluding}

In this paper, we developed an efficient and robust fifth-order finite-volume HWENO scheme with gradient reconstruction for the compressible Navier--Stokes equations on structured meshes. The convective reconstruction uses a HWENO scheme with compact unified stencils, while the main contribution is the compatible gradient reconstruction for the viscous fluxes. The weak-derivative moments of the dissipative variables are formed from HWENO interface traces and reconstructed by the same nonlinear HWENO procedure as the solution variables. Since only the input moments change, the conservative and gradient reconstructions share the same stencil-dependent coefficient matrices and can be implemented with one reconstruction routine. This avoids direct differentiation of reconstructed solution polynomials, does not require a separate derivative-reconstruction implementation, and gives viscous fluxes with the same fifth-order accuracy as the convective fluxes.

The HWENO scheme with unified stencils supplies the compact stencil layout and the high-order modified first-order moments used only in time discretization, while the PP-limiter preserves positive density and pressure. The accuracy tests confirm fifth-order convergence for smooth compressible Navier--Stokes solutions. For the fully nonlinear implementations compared here, the error and CPU-time data are consistent with the reduced reconstruction work. The non-smooth finite-Reynolds-number examples give stable computations, comparable main wave locations, and small differences in local contour structures in several two-dimensional tests.

\section{Acknowledgements}\label{Sec6:Acknowledgements}
The research was partially supported by the National Key R\&D Program of China under Grant Number
2024YFA1012500 and Fujian Provincial Natural Science Foundation of China under Grant Number 2026J009009. The authors declare that they have no financial or non-financial conflicts of interest.

\appendix
\section{HWENO reconstruction details}
\label{app:nonlinear-hweno}

The formulas below describe the nonlinear HWENO reconstruction with unified stencils used in this paper; further details can be found in \cite{FanQiuZhao2024}. Here, the
same reconstruction is applied with different input moments. For the convective fluxes,
the inputs are the conservative cell averages and first-order moments. For the viscous
terms, the inputs are the weak-derivative moments defined in the main text. Thus the
candidate polynomials, smoothness indicators, and nonlinear weights are recomputed from
the chosen input data, while the stencil-dependent formulas are unchanged.

\subsection{One-dimensional nonlinear HWENO reconstruction}
\label{app:nonlinear-hweno1}

\begin{figure}[!htbp]
\centering
\includegraphics{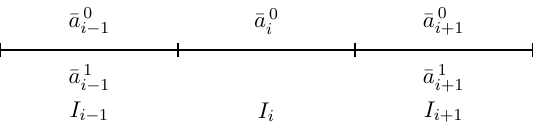}
\caption{One-dimensional stencil for the HWENO reconstruction. The moment data are
attached to cell centers; the entries used by each candidate polynomial are specified by
the moment conditions below.}
\label{fig:app-1d-stencil}
\end{figure}

Let \(\bar a_k^{\,0}\) and \(\bar a_k^{\,1}\) denote the input zeroth-order and first-order
moments on \(I_k\). In the conservative reconstruction,
\((\bar a_k^{\,0},\bar a_k^{\,1})=(\bar u_k,\bar v_k)\);
in the weak-derivative reconstruction,
\((\bar a_k^{\,0},\bar a_k^{\,1})=(\bar r_k^{\,0},\bar r_k^{\,1})\). The stencil
layout is shown in Fig.~\ref{fig:app-1d-stencil}.

The large stencil \(S_0=\{I_{i-1},I_i,I_{i+1}\}\) defines a quartic candidate
\(p_0\) by
\begin{equation}\label{sec:grad-p0-1d}
\begin{aligned}
&\frac1{\dx}\int_{I_k}p_0(x)\,\mathrm{d}x=\bar a_k^{\,0},
&& k=i-1,i,i+1,\\
&\frac1{\dx}\int_{I_k}p_0(x)\frac{x-x_k}{\dx}\,\mathrm{d}x=\bar a_k^{\,1},
&& k=i-1,i+1.
\end{aligned}
\end{equation}
The two linear candidates on \(S_1=\{I_{i-1},I_i\}\) and
\(S_2=\{I_i,I_{i+1}\}\) satisfy
\begin{equation}\label{sec:grad-small-1d}
\frac1{\dx}\int_{I_k}p_1(x)\,\mathrm{d}x=\bar a_k^{\,0},\quad k=i-1,i,
\qquad
\frac1{\dx}\int_{I_k}p_2(x)\,\mathrm{d}x=\bar a_k^{\,0},\quad k=i,i+1.
\end{equation}
The smoothness indicators \(\beta_m\), \(m=0,1,2\), are evaluated from these
candidate polynomials using the scale-invariant formulas of the HWENO scheme with
unified stencils \cite{FanQiuZhao2024}. For positive linear weights
\(\gamma_m\), \(m=0,1,2,\) with \(\sum_{m=0}^{2}\gamma_m=1\), the nonlinear weights are
\begin{equation}\label{sec:grad-weights}
\begin{aligned}
\omega_m=\frac{\widetilde\omega_m}{\sum_{s=0}^{2}\widetilde\omega_s}, \ \ \
\widetilde\omega_m=\gamma_m\left(1+\frac{\tau}{\beta_m+\varepsilon_{\mathrm{WENO}}}\right), \ \ \
\tau=\left(\frac{1}{2}\sum_{s=1}^{2}|\beta_0-\beta_s|\right)^2 .
\end{aligned}
\end{equation}
Here the linear weights are set as \(\gamma_0=199/200\) and \(\gamma_1=\gamma_2=1/400\), and
\(\varepsilon_{\mathrm{WENO}}=10^{-8}\) in all computations as in \cite{FanQiuZhao2024}.

The reconstructed value is then obtained from the nonlinear combination in
\eqref{sec:solution-value-1d} or \eqref{sec:grad-value-1d}, depending only on whether
the input moments are conservative moments or weak-derivative moments.

\subsection{Two-dimensional nonlinear HWENO reconstruction}
\label{app:nonlinear-hweno2}
\begin{figure}[!htbp]
\centering
\includegraphics{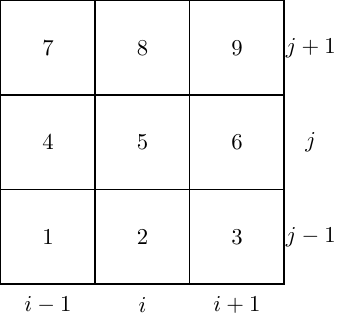}
\caption{Indexing of the \(3\times3\) cell stencil for the two-dimensional HWENO
reconstruction, with \(I_5=I_{i,j}\).}
\label{fig:app-2d-stencil}
\end{figure}

Let \(\bar a_k^{0}\), \(\bar a_k^{\xi}\), and \(\bar a_k^{\eta}\) denote the input
moments on a cell \(I_k\). For conservative variables, these are
\(\bar u_k\), \(\bar v_k\), and \(\bar w_k\). For the
weak derivative in direction \(\ell=x,y\), they are
\(\bar r_k^{\ell,0}\), \(\bar r_k^{\ell,\xi}\), and
\(\bar r_k^{\ell,\eta}\). Relabel the \(3\times3\) cells around \(I_{i,j}\) as
\(I_1,\ldots,I_9\), with \(I_5=I_{i,j}\). The indexing convention is shown in
Fig.~\ref{fig:app-2d-stencil}.

The large-stencil candidate \(p_0\) uses the zeroth-order input moments from all nine
cells and the first-order input moments from the four edge-neighboring cells:
\begin{equation}\label{sec:grad-p0-2d}
\begin{aligned}
&\frac1{\dx\dy}\int_{I_k}p_0(x,y)\,\mathrm{d}x\mathrm{d}y
=\bar a_k^{0},
&& k=1,\ldots,9,\\
&\frac1{\dx\dy}\int_{I_k}p_0(x,y)\frac{x-x_k}{\dx}\,\mathrm{d}x\mathrm{d}y
=\bar a_k^{\xi},
&& k=2,4,6,8,\\
&\frac1{\dx\dy}\int_{I_k}p_0(x,y)\frac{y-y_k}{\dy}\,\mathrm{d}x\mathrm{d}y
=\bar a_k^{\eta},
&& k=2,4,6,8 .
\end{aligned}
\end{equation}
As in the two-dimensional HWENO reconstruction with unified stencils, the nine zeroth-order conditions,
the \(\bar a_k^{\xi}\) conditions for \(k=4,6\), and the \(\bar a_k^{\eta}\)
conditions for \(k=2,8\) are imposed exactly; the remaining first-order moment
conditions are included in a least-squares sense. The four linear candidate polynomials
use only zeroth-order input moments on the compact three-cell substencils:
\begin{equation}\label{sec:grad-small-2d}
\frac1{\dx\dy}\int_{I_k}p_m(x,y)\,\mathrm{d}x\mathrm{d}y
=\bar a_k^{0},
\qquad k\in\mathcal{S}_m,
\end{equation}
where \(\mathcal{S}_1=\{2,4,5\}\), \(\mathcal{S}_2=\{2,5,6\}\),
\(\mathcal{S}_3=\{4,5,8\}\), and \(\mathcal{S}_4=\{5,6,8\}\).
The five smoothness indicators \(\beta_m\), \(m=0,\ldots,4\), are again those of
\cite{FanQiuZhao2024}, evaluated with the same input moments as the candidate
polynomials. For positive linear weights \(\gamma_m\), \(m=0,\ldots,4\), with
\(\sum_{m=0}^4\gamma_m=1\), the nonlinear weights are
\begin{equation}\label{sec:hweno-weights-2d}
\begin{aligned}
\omega_m=\frac{\widetilde\omega_m}{\sum_{s=0}^{4}\widetilde\omega_s}, \ \ \
\widetilde\omega_m=\gamma_m\left(1+\frac{\tau}{\beta_m+\varepsilon_{\mathrm{WENO}}}\right), \ \ \
\tau=\left(\frac{1}{4}\sum_{s=1}^{4}|\beta_0-\beta_s|\right)^2.
\end{aligned}
\end{equation}
Here the linear weights are set as \(\gamma_0=99/100\) and \(\gamma_1=\gamma_2=\gamma_3=\gamma_4=1/400\), and
\(\varepsilon_{\mathrm{WENO}}=10^{-8}\) as in \cite{FanQiuZhao2024}.

The resulting nonlinear combination is evaluated by \eqref{sec:solution-value-2d} for
conservative states and by \eqref{sec:grad-value-2d} for weak-derivative data, where
 the stored polynomial and indicator formulas are therefore shared by the
convective and weak-derivative reconstructions, with only the input moment array being changed.

\end{document}